\documentclass[12pt,leqno]{amsart}
\usepackage{amssymb,mathrsfs,euscript}
\usepackage[hypertex]{hyperref}
\catcode`\@=11
\def\@evenfoot{\rule{0pt}{20pt}[\today] \hfill}
\def\@oddfoot{\rule{0pt}{20pt}\hfill [\today]}

\newtheorem{theorem}{Theorem}
\newtheorem{definition}[theorem]{Definition}
\newtheorem{convention}[theorem]{Convention}

\newtheorem{example}[theorem]{Example}
\newtheorem{lemma}[theorem]{Lemma}
\newtheorem{claim}[theorem]{Claim}

\newtheorem{proposition}[theorem]{Proposition}

\newtheorem{notation}[theorem]{Notation}
\newtheorem*{theoremA}{Theorem~A}
\newtheorem*{theoremB}{Theorem~B}

\def\bad{{\mathcal B}ad}
\def\gls#1#2#3{#1,&\mbox{#2,}&\mbox{ page~\pageref{#3}}}
\def\sphere{{\mathbb S}}\def\ob{{\mathfrak o}}
\def\glsfin#1#2#3{#1,&\mbox{#2,}&\mbox{ page~\pageref{#3}}}
\def\genname{g} \def\Chain{{\EuScript Chain}}
\def\BE{{\EuScript E}}\def\Surjc{{\EuScript X}}
\def\EZ{{\EuScript Z}} \def\Lat{{\EuScript L}} \def\CG{{\EuScript K}}
\def\ss{{\mathbf s}}\def\pp#1{\partial_{\rm prt}^{#1}}
\DeclareMathOperator{\colim}{colim}

\def\Miki{
{
\unitlength=.800000pt
\begin{picture}(40.00,40.00)(0.00,0.00)
\put(20.00,21.00){\makebox(0.00,0.00){$\bullet$}}
\put(20.00,8.00){\makebox(0.00,0.00){$\bullet$}}
\put(20.00,30.00){\makebox(0.00,0.00){$\bullet$}}
{\thinlines
\put(20.00,40.00){\line(0,-1){40.00}}
}
\put(20.00,20.00){\circle{40.00}}
\end{picture}}
}

\def\hFlicek{
{
\unitlength=.800000pt
\begin{picture}(40.00,40.00)(0.00,0.00)
\put(20.00,10.00){\makebox(0.00,0.00){$\bullet$}}
\put(23.00,25.00){\makebox(0.00,0.00){$\bullet$}}
\put(14.00,27.00){\makebox(0.00,0.00){$\bullet$}}
{\thinlines
\put(20.00,40.00){\line(0,-1){3.5}}
\put(20.00,16.00){\line(0,-1){16.00}}
\put(20.00,26.00){\circle{20.00}}
}
\put(20.00,20.00){\circle{40.00}}
\end{picture}}
}

\def\dFlicek{
{
\unitlength=.800000pt
\begin{picture}(40.00,40.00)(0.00,0.00)
\put(20.00,33.00){\makebox(0.00,0.00){$\bullet$}}
\put(25.00,12.00){\makebox(0.00,0.00){$\bullet$}}
\put(16.00,11.00){\makebox(0.00,0.00){$\bullet$}}
{\thinlines
\put(20.00,40.00){\line(0,-1){18.00}}
\put(20.00,2.00){\line(0,-1){2.00}}
\put(20.00,12.00){\circle{20.00}}
}
\put(20.00,20.00){\circle{40.00}}
\end{picture}}
}

\def\Vojacek{
{
\unitlength=.800000pt
\begin{picture}(40.00,40.00)(0.00,0.00)
\put(33.00,22.00){\makebox(0.00,0.00){$\bullet$}}
\put(17.00,18.00){\makebox(0.00,0.00){$\bullet$}}
\put(8.00,21.00){\makebox(0.00,0.00){$\bullet$}}
\put(20.00,20.00){\circle{40.00}}
\end{picture}}
}

\def\lKulda{
{
\unitlength=.800000pt
\begin{picture}(40.00,40.00)(0.00,0.00)
\put(30.00,24.00){\makebox(0.00,0.00){$\bullet$}}
\put(9.00,21.00){\makebox(0.00,0.00){$\bullet$}}
\put(16.00,18.00){\makebox(0.00,0.00){$\bullet$}}
{\thinlines
\put(13.00,20.00){\circle{20.00}}
}
\put(20.00,20.00){\circle{40.00}}
\end{picture}}
}

\def\pKulda{
{
\unitlength=.800000pt
\begin{picture}(40.00,40.00)(0.00,0.00)
\put(24.00,18.00){\makebox(0.00,0.00){$\bullet$}}
\put(33.00,20.00){\makebox(0.00,0.00){$\bullet$}}
\put(8.00,18.00){\makebox(0.00,0.00){$\bullet$}}
{\thinlines
\put(28.00,20.00){\circle{20.00}}
}
\put(20.00,20.00){\circle{40.00}}
\end{picture}}
}

\def\pKrtek{
{
\unitlength=.800000pt
\begin{picture}(40.00,40.00)(0.00,0.00)
\put(9.00,18.00){\makebox(0.00,0.00){$\bullet$}}
\put(30.00,12.00){\makebox(0.00,0.00){$\bullet$}}
\put(30.00,30.00){\makebox(0.00,0.00){$\bullet$}}
{\thinlines
\put(30.00,37.70){\line(0,-1){35.00}}
}
\put(20.00,20.00){\circle{40.00}}
\end{picture}}
}

\def\lKrtek{
{
\unitlength=.800000pt
\begin{picture}(40.00,40.00)(0.00,0.00)
\put(33.00,22.00){\makebox(0.00,0.00){$\bullet$}}
\put(12.00,10.00){\makebox(0.00,0.00){$\bullet$}}
\put(12.00,28.00){\makebox(0.00,0.00){$\bullet$}}
{\thinlines
\put(12.00,39.00){\line(0,-1){37.0}}
}
\put(20.00,20.00){\circle{40.00}}
\end{picture}}
}

\def\lineup#1#2{
{
\unitlength=.800000pt
\begin{picture}(50.00,80.00)(-5,75)
\thicklines
\put(20.00,40.00){\makebox(0.00,0.00){#1}}
\put(20.00,110.00){\makebox(0.00,0.00){#2}}
\put(20.00,20.00){\line(0,-1){15.00}}
\put(20.00,90.00){\line(0,-1){30.00}}
\put(20.00,145.00){\line(0,-1){15.00}}
\end{picture}}
}

\def\Abox#1{
\unitlength=.800000pt
\begin{picture}(0.00,0.00)
\thinlines
\put(0,-65){
\put(248.00,128.00){\makebox(0.00,0.00)[tr]{\scriptsize $#1$}}
\put(0.00,00){\line(1,0){250.00}}
\put(0.00,00){\line(0,1){130.00}}
\put(0.00,130){\line(1,0){250.00}}
\put(250,0){\line(0,1){130.00}}
}
\end{picture}
}

\def\Bbox#1{
\unitlength=.800000pt
\begin{picture}(0.00,0.00)
\thinlines
\put(0,-65){
\put(183.00,128.00){\makebox(0.00,0.00)[tr]{\scriptsize $#1$}}
\put(0.00,00){\line(1,0){185.00}}
\put(0.00,00){\line(0,1){130.00}}
\put(0.00,130){\line(1,0){185.00}}
\put(185,0){\line(0,1){130.00}}
}
\end{picture}
}

\def\hatB{\widehat{\sf B}}
\def\xxx{{\mathbf x}}
\def\velkasvisla{\hbox{$\left|\rule{0pt}{.9em}\right.$}}
\def\muu#1#2{\mu \hskip -.2em \left[\hskip -.2em \rule{0pt}{1.3em} \right.#1 
\hskip -.2em\left|\hskip -.2em
    \rule{0pt}{1.3em} \right. #2
    \hskip -.2em \left. \rule{0pt}{1.3em} \right]  }
\def\jarka#1#2{\left[\hskip -.2em \rule{0pt}{.9em} \right.#1 
       \hskip -.2em\left | \hskip -.2em
    \rule{0pt}{.9em} \right. #2
    \hskip -.2em\left. \rule{0pt}{.9em} \right]  }
\def\jarkA#1#2#3{\left[\hskip -.2em \rule{0pt}{.9em} \right.#1 
       \hskip -.2em\left | \hskip -.2em
    \rule{0pt}{.9em} \right. #2\hskip -.2em\left | \hskip -.2em
    \rule{0pt}{.9em} \right. #3
    \hskip -.2em\left. \rule{0pt}{.9em} \right]  }
\def\Jarka#1#2{\left[\hskip -.2em \rule{0pt}{.9em} \right.#1 
       \hskip -.2em \left|\rule{0pt}{.9em}   \right|  \hskip -.1em
    \rule{0pt}{.9em}  #2
    \hskip -.2em\left. \rule{0pt}{.9em} \right]  }
\def\JArka#1#2{\left[\hskip -.2em \rule{0pt}{.9em} \right.#1 \hskip .3em 
       \hbox{$\hskip -.2em \left|\rule{0pt}{.9em}\right. \hskip -.365em
      \left|\rule{0pt}{.9em} \hskip -.1em  \right|  \hskip -.1em
    \rule{0pt}{.9em}$}  #2
    \hskip -.2em\left. \rule{0pt}{.9em} \right]  }
\def\bfS{{\mathbf S}}\def\Tam{{\mathscr T}}\def\Mar{{\mathscr M}}
\def\bfT{{\mathbf T}}
\def\pomio{\hbox{${\mathbf s}\hskip-.3em%
           \downarrow\hskip-.3em\overline{\EuScript C}$}}
\def\Aff#1{{\mbox{\rm Aff}({\mathbb R}^h)}}\def\dcrit{d_{\rm crit}}
\def\osfF{\hbox{\hskip.1em\raisebox{.8em}{\scriptsize$\circ$}\hskip-.5em{\sf
          F}}}
\def\dirlim{\hbox{\hskip.1em\raisebox{-.5em}
           {\scriptsize $\longrightarrow$}\hskip-1.35em{\rm lim}\hskip .2em}}
\def\sfF{{\sf F}}\def\opal{{\overline{\partial}_{\hskip .1em \rm lin}}}
\def\Free{{\mathbb F}}\def\uFree{\underline{\mathbb F}}
\def\opa{\overline{\partial}} 
\def\End{{\mathcal E \hskip -.1em \it nd}}
\def\opapert{{\overline{\partial}_{\rm prt}}}
\def\papert{{\partial_{\rm prt}}}\def\pareg{{\partial_{\rm reg}}}
\def\pal{{\partial_{\hskip .1em \rm lin}}}\def\Span{{\mbox{\rm Span}}}
\def\uG{\underline {G}}\def\pasing{{\partial_{\rm sng}}}
\def\smod#1#2{\{#1(n)\}_{n \geq #2}}\def\scrF{{\mathscr{F}}}
\def\semidirect{\mbox{$\rule{.15mm}{1.9mm}\hskip-.75mm\times$}}
\def\tria{
{
\unitlength=.55pt
\begin{picture}(20.00,60.00)(0.00,0.00)
\put(20.00,60.00){\line(-1,-6){10.00}}
\put(0.00,60.00){\line(1,0){20.00}}
\put(10.00,0.00){\line(-1,6){10.00}}
\end{picture}}
}

\def\krouzek#1{
\setlength{\unitlength}{#1cm}
\bezier{100}(-2,0)(-2,.5)(-1.73,1)
\bezier{100}(-1.73,1)(-1.48,1.48)(-1,1.73)
\bezier{100}(-1,1.73)(-.5,2)(0,2)
\bezier{100}(2,0)(2,.5)(1.73,1)
\bezier{100}(1.73,1)(1.48,1.48)(1,1.73)
\bezier{100}(1,1.73)(.5,2)(0,2)
\bezier{100}(2,0)(2,-.5)(1.73,-1)
\bezier{100}(1.73,-1)(1.48,-1.48)(1,-1.73)
\bezier{100}(1,-1.73)(.5,-2)(0,-2)
\bezier{100}(-2,0)(-2,-.5)(-1.73,-1)
\bezier{100}(-1.73,-1)(-1.48,-1.48)(-1,-1.73)
\bezier{100}(-1,-1.73)(-.5,-2)(0,-2)
}
\def\calC{{\EuScript C}}\def\scrE{{\mathscr E}}
\def\sx#1{\hbox{$\uparrow \! x_{#1}$}}\def\pabar{{\partial_{\sfB}}}
\def\PTree{{\rm PTree}}\def\Leaf{{\it Leaf}}\def\Int{{ar}}
\def\Ush{{\mathrm {Ush}}}\def\Conf #1#2{{\mathrm{Cnf}(\RRR^{#1},#2)}}
\def\CConf #1#2{{\hskip1em\hbox{\raisebox{.7em}%
                {\scriptsize$\bullet$}\hskip-1.3em
		Cnf}(\RRR^{#1},#2)}}
\def\Ivosek{{\hskip.4em\hbox{\raisebox{.7em}%
                {\scriptsize$\bullet$}\hskip-.85em $[\bfT]$}}}
\def\Tree{{\mathrm {Tree}}}\def\Flag{{\mathrm {Flag}}}
\def\Treereg{{\rm Tree}^{\rm reg}}\def\Treeregh{{\rm Tree}^{{\rm
      reg},h}}
\def\Treeregjedna{{\rm Tree}^{{\rm reg},1}}
\def\uTreereg{\underline{\rm Tree}^{\rm reg}}
\def\uTreeregh{\underline{\rm Tree}^{{\rm reg},h}}
\def\Treeinf{{\mathrm {Tree}}}\def\Vert{{\it Vert}}
\def\uTree{{\underline{\mathrm{Tree}}}}\def\Greg{G^{\rm reg}}
\def\uFlag{{\underline{\mathrm{Fla}}\mathrm{g}}}
\def\uTreeinf{{\underline{\mathrm{Tree}}}}
\def\ZZZ{{\mathbb Z}}\def\TTT{{\mathbb T}}\def\LLL{{\mathbb L}}
\def\RRR{{\mathbb R}}\def\epi{ \twoheadrightarrow}
\def\vlra{{\hbox{$-\hskip-1mm-\hskip-2mm\longrightarrow$}}}
\def\pa{\partial}\def\sgn{{\rm sgn}}
\def\otexp#1#2{#1^{\otimes #2}}\def\sfB{{\sf B}}
\def\cases#1#2#3#4{
                  \left\{
                         \begin{array}{ll}
                           #1,\ &\mbox{#2}
                           \\
                           #3,\ &\mbox{#4}
                          \end{array}
                   \right.
}
\def\tricases#1#2#3#4#5#6{
                  \left\{
                         \begin{array}{ll}
                           #1,\ &\mbox{#2}
                           \\
                           #3,\ &\mbox{#4}
                           \\
                           #5,\ &\mbox{#6}
                          \end{array}
                   \right.
}
\def\ot{\otimes}\def\Mbar{\overline{M}}
\def\id{{1\!\!1}}\def\sfFreg{{\sfF}^{\rm reg}}
\def\rada#1#2{#1,\ldots,#2}
\def\Rada#1#2#3{#1_{#2},\dots,#1_{#3}}
\def\Ass{\mbox{${\mathcal A}${\it ss}}}
\def\Com{\mbox{${\mathcal C}$\hskip -.2mm {\it om}}}

\def\lp{
{\hskip .5mm
\unitlength.1cm
\begin{picture}(3.00,3.00)
\thinlines
  \put(1.5,0){\bezier{30}(0,0)(.5,.5)(1,1)}
  \put(1.5,0){\bezier{30}(0,0)(-.5,.5)(-1,1)}
  \put(2.5,1){\line(0,1){1}}
  \put(.5,1){\bezier{30}(0,0)(-.5,.5)(-1,1)}
  \put(.5,1){\bezier{30}(0,0)(.5,.5)(1,1)}
\end{picture}
\hskip .5mm
}
}

\def\rp{
{\hskip .5mm
\unitlength.1cm
\begin{picture}(3.00,3.00)
\thinlines
  \put(1.5,0){\bezier{30}(0,0)(.5,.5)(1,1)}
  \put(1.5,0){\bezier{30}(0,0)(-.5,.5)(-1,1)}
  \put(.5,1){\line(0,1){1}}
  \put(2.5,1){\bezier{30}(0,0)(-.5,.5)(-1,1)}
  \put(2.5,1){\bezier{30}(0,0)(.5,.5)(1,1)}
\end{picture}
\hskip .5mm
}
}

\def\dvojdum{
{\hskip .5mm
\unitlength.1cm
\begin{picture}(3.00,3.00)
\thinlines
  \put(1,0){\bezier{30}(0,0)(.5,.5)(1,1)}
  \put(1,0){\bezier{30}(0,0)(-.5,.5)(-1,1)}
  \put(0,1){\line(0,1){1}}
  \put(2,1){\line(0,1){1}}
\end{picture}
\hskip .5mm
}
}

\def\levaslozitost{
{\hskip .9mm
\unitlength.1cm
\begin{picture}(4.00,3.00)
\thinlines
  \put(1,0){\bezier{30}(0,0)(.5,.5)(1,1)}
  \put(1,0){\bezier{30}(0,0)(-.5,.5)(-1,1)}
\put(-1,2){  
  \put(1,0){\bezier{30}(0,0)(.5,.5)(1,1)}
  \put(1,0){\bezier{30}(0,0)(-.5,.5)(-1,1)}
}
  \put(0,1){\line(0,1){1}}
  \put(2,1){\line(0,1){2}}
\end{picture}
\hskip .5mm
}
}

\def\trojdum{
{\hskip .5mm
\unitlength.1cm
\begin{picture}(3.00,4.00)
\thinlines
  \put(1,0){\bezier{30}(0,0)(.5,.5)(1,1)}
  \put(1,0){\bezier{30}(0,0)(-.5,.5)(-1,1)}
  \put(0,1){\line(0,1){2}}
  \put(2,1){\line(0,1){2}}
\end{picture}
\hskip .5mm
}
}

\def\Lie{\hbox{{$\mathcal L$}{\it ie\/}}}
\def\strecha{
{\hskip .5mm
\unitlength.1cm
\begin{picture}(3.00,2.00)
\thinlines
  \put(1,0){\bezier{30}(0,0)(.5,.5)(1,1)}
  \put(1,0){\bezier{30}(0,0)(-.5,.5)(-1,1)}
\end{picture}
\hskip .5mm
}
}

\def\trident{
{\hskip .5mm
\unitlength.1cm
\begin{picture}(3.00,3.00)
\thinlines
  \put(1,0){\line(0,1){2}}
  \put(1,1){\bezier{30}(0,0)(.5,.5)(1,1)}
  \put(1,1){\bezier{30}(0,0)(-.5,.5)(-1,1)}
\end{picture}
\hskip .5mm
}
}

\def\assoc{
{\hskip .5mm
\unitlength.1cm
\begin{picture}(3.00,2.00)
\thinlines
  \put(1,0){\line(0,1){1}}
  \put(1,0){\bezier{30}(0,0)(.5,.5)(1,1)}
  \put(1,0){\bezier{30}(0,0)(-.5,.5)(-1,1)}
\end{picture}
\hskip .5mm
}
}

\def\dvojtroj{
{\hskip .5mm
\unitlength.1cm
\begin{picture}(3.00,2.00)
\thinlines
  \put(1,0){\line(0,1){2}}
  \put(0,1){\line(0,1){1}}
  \put(2,1){\line(0,1){1}}
  \put(1,0){\bezier{30}(0,0)(.5,.5)(1,1)}
  \put(1,0){\bezier{30}(0,0)(-.5,.5)(-1,1)}
\end{picture}
\hskip .5mm
}
}

\def\pent{
{\hskip .5mm
\unitlength.1cm
\begin{picture}(3.00,2.00)
\thinlines
  \put(1,0){\bezier{30}(0,0)(.5,.5)(1,1)}
  \put(1,0){\bezier{30}(0,0)(-.5,.5)(-1,1)}
  \put(1,0){\bezier{30}(0,0)(.2,.5)(.4,1)}
  \put(1,0){\bezier{30}(0,0)(-.2,.2)(-.4,1)}
\end{picture}
\hskip .5mm
}
}

\long\def\comment#1\endcomment{}

\title[$E_\infty$-analog of the associahedra]%
{An $E_\infty$-extension of the associahedra  and \\ the
  Tamarkin cell mystery \rule{0em}{1.2em}}
\author[M. MARKL]{M. MARKL}

\thanks{Supported by the grants GA \v CR 201/02/1390 and GA \v CR 201/08/0397.}

\catcode`\@=11
\address{Mathematical Institute of the Academy\\
          \v Zitn\'a 25\\
          115 67 Praha 1\\ The Czech~Republic}
\email{markl@math.cas.cz}
\catcode`\@=13

\begin{document}
\bibliographystyle{plain}

\begin{abstract}
In this note based on the author's communication with {\bf
M.~Batanin}, we study a cofibrant $E_\infty$-operad generated by the
Fox-Neuwirth cells of the configuration space of points in the
Euclidean space.  We show that, below the `critical dimensions' in
which `bad cells' exist, this operad is modeled by the geometry of
the Fulton-MacPherson compactification of this configuration space. We
analyze the Tamarkin bad cell and calculate the differential of the
corresponding generator. We also describe a simpler,
four-dimensional bad cell. We finish the paper by proving an auxiliary
result giving a characterization, over integers, of free Lie algebras.
\end{abstract}

\maketitle

\baselineskip 16pt plus 1pt minus 1pt

\noindent 
{\bf Conventions.}
All algebraic objects will be considered over the ring
of integers $\ZZZ$.  Terminology regarding operads and related
constructions follows~\cite{markl-shnider-stasheff:book}. With few
exceptions, all operads in this paper live in the monoidal 
category $\Chain$ of
non-negatively graded chain complexes of abelian groups. If we consider
cofibrant operads, we refer to the model structure of the category of
operads in $\Chain$ considered
in~\cite[Example~3.3.3]{berger-moerdijk:02}. The central operad of
this paper, $J = (J,\pa)$, will in fact be special cofibrant in a more specific
sense of~\cite[page~143]{markl:ha}.

\vskip .5em 
\noindent 
{\bf Motivations and historical overview.}
Cofibrant resolutions of the operad $\Com$ for
commutative associative algebras play a key r\^ole in
(co)homology theory of commutative or $E_\infty$-algebras,
applications to the
iterated bar construction, infinite loop spaces, and many other areas
of homological algebra and topology. Surprisingly,
not many explicit and combinatorially accessible cofibrant resolutions
are known.

- There is the simplicial {\em Barratt-Eccles\/} 
operad~\cite{barratt-eccles:Top74}
whose chain version $\BE$ was introduced in~\cite{berger-fresse}. The
{\em surjection operad\/} $\Surjc$ was described as a quotient of $\BE$
in~\cite{berger-fresse};  it appeared independently
in~\cite{mcclure-smith}. One also has the {\em Eilenberg-Zilber\/} operad
$\EZ$ of all natural operations $N_*(S) \to N_*(S)^{\otimes n}$ on the
normalized chain complex $N_*(S)$ of a simplicial set $S$,
see~\cite{hinich-schechtman:LNM1289,smirnov:IAN85}. As
proved in~\cite{berger-fresse}, these three operads are related by the
maps
\begin{equation}
\label{dnes_slavim_v_praci_narozeniny}
\BE \stackrel{{\it TR}}\longrightarrow \Surjc 
\stackrel{{\it AW}}\longrightarrow \EZ,  
\end{equation}
where {\it TR} is the table reduction map and {\it AW} abbreviates
Alexander-Whitney.  

Unfortunately, {\em none\/} of the operads
in~(\ref{dnes_slavim_v_praci_narozeniny}) is cofibrant. But, at least,
the components of the operads $\BE(n)$ and $\Surjc(n)$ are free
modules over the symmetric group $\Sigma_n$, 
so one can get cofibrant resolutions $W(\BE)$ and
$W(\Surjc)$ of $\Com$ by applying the cellular version $W(-)$ of
Boardman-Vogt's $W$-construction~\cite{boardman-vogt:73}. The
$W$-construction in fact
coincides with the operadic double bar construction and substantially
inflates the size of the operad to which it is applied.

-- Another chain operads resolving $\Com$ are the chain condensation
of the {\em lattice path operad\/} $|\Lat|$~\cite{bb,batanin-berger-markl}
and the normalization $N(\CG)$ of the complete graph operad
$\CG$~\cite{berger-fresse}. The surjection operad $\Surjc$ is a
suboperad of $|\Lat|$ via the {\em whiskering map\/}, and C.~Berger
informed us that there is a simple zigzag between $|\Lat|$ and
$\CG$. As in the previous item, the operads  $|\Lat|$ and $\CG$ are
only $\Sigma$-free, not cofibrant, so the $W(-)$-functor must be
applied to obtain cofibrant resolutions.

-- In~\cite{fresse:kd}, the existence of cofibrant models of 
$E_n$-operads of the form $B^c(\ss^{-n} E_n^\vee)$, where $B^c(-)$ is
the cobar construction of a cooperad, $\ss^{-n}$ the $n$-fold
operadic desuspension and  $ E_n^\vee$ is the dual cooperad of an
explicit filtration layer of the Barratt-Eccles operad, was proved for
each $n \geq 0$.
The colimit $\colim_nB^c(\ss^{-n} E_n^\vee)$ is a combinatorial
cofibrant model for $\Com$. While its operad structure
is explicit, there is no simple formula for the augmentation to
$\Com$.

The sizes of the cofibrant resolutions $W(\BE)$, $W(\Surjc)$,
$W(N(\CG))$, $W(|\Lat|)$ and $\colim_nB^c(\ss^{-n} E_n^\vee)$ of
$\Com$ recalled above are huge.  There has been, however, a
long-standing candidate for a~small cofibrant chain resolution
proposed by E.~Getzler and J.D.S.~Jones
in~\cite{getzler-jones:preprint}, whose space of generators is of the
size of the Fox-Neuwirth cell decomposition of the configuration
space. It unfortunately turned out that, due to the existence of `bad
cells,' the proposed combinatorial formula for the differential $\pa$
did not work above certain dimensions.

The {\bf purpose of this paper\/} is two-fold. The first one is to
show that the differential $\partial$ can be modified to a correct one
that coincides with the one proposed in~\cite{getzler-jones:preprint}
below the dimensions of the bad cells. 
The second aim is to show that explicit calculations 
of the values of $\partial$ on some bad cells are possible. The methods and
results of the paper are described in more detail in the introduction
below where also the two main results of the paper, Theorems~A and~B,
are formulated.

Our results imply that, in the
applications mentioned in the first paragraph of this subsection,
Getzler-Jones' formula  can be used in dimensions less than the
dimension of the bad cells, and that some explicit results can be
obtained also in dimensions containing bad cells. A closed
formula for $\pa$ is still, however, a challenging open question.  
We would like to note that the 
initial pieces of the cellular Getzler-Jones operad have in fact 
already been calculated in the proceedings~\cite{markl:ws93} 
of the Winter School `Geometry and Physics,' Zd\'\i kov,
Bohemia, January 1993.  

\tableofcontents

\baselineskip 17pt plus 1pt minus 1pt

\section{Introduction}
\label{int}

We describe an operad $J$
that can be viewed as an $E_\infty$-analog of the minimal model
${\mathcal A}_\infty$ of the operad ${\Ass}$ for associative algebras,
governing Stasheff's $A_\infty$-algebras~\cite{markl:zebrulka}.  Our
$J$ lives in the monoidal category of differential graded 
abelian groups. It is of the form $J = (\Free(G),\pa)$ where\label{jarca!}

- $\Free(G)$ is the free operad generated by the graded
  $\Sigma$-module $G = \{G_*(n)\}_{n \geq 2}$ specified in
  Definition~\ref{JArKa}, and

- the differential $\pa$ is the sum $\pa = \pal + \papert$ of the
  `linear' part $\pal$ induced from a differential (denoted by the
  same symbol) $\pal$ on the $\Sigma$-module $G$ introduced in
  Definition~\ref{jarka}, and the `perturbed' part $\papert$ that maps
  $G$ into the decomposables of the free operad $\Free(G)$. Below the
  `critical dimension' in which the `bad' cells exist, explicitly
  specified in Definition~\ref{JarusKa}, $J$ is determined by the cell
  structure of the configuration operad $\sfF$ induced from the
  Fox-Neuwirth decomposition of the configuration space. Moreover,

- the operad $J$ is equipped with a dg-operad homomorphism $\rho : J
  \to \Com$ that makes $J$ a cofibrant resolution of the operad
  $\Com = \End_{\ZZZ}$ for commutative associative algebras.

The $n$th piece $(G_*(n),\pal)$ of the generating $\Sigma$-module $G$
is, for each $n \geq 2$, the colimit of the  (shifted) cellular chain
complexes of the one-point compactifications of the configuration
spaces \label{qaz}$\Conf hn$ of $n$ distinct labeled points in the
$h$-dimensional Euclidean space $\RRR^h$, with the Fox-Neuwirth cell
structure.  Alternative descriptions of the right
$\Sigma_n$-dg-abelian group $(G_*(n),\pal)$ are given in
Section~\ref{JARCA}. Summing up, we prove

\begin{theoremA}
There exist a cofibrant resolution $\rho : J = (\Free(G),\pa) \to
\Com$ such that the linear part $\pal$ of the differential
is as in Definition~\ref{jarka} and the restriction of
$\pa$ to generators below the dimension of bad cells is
determined by the cell structure of the configuration operad \ $\sfF$.
\end{theoremA}

Theorem~A is proved in Subsection~\ref{Jarka_mne_mozna_miluje}. 
Since the image of the canonical embedding ${\sf K} \hookrightarrow
\sfF$ of the Stasheff's associahedron ${\sf K}$ into
the configuration operad $\sfF$ does not contain bad cells 
(see Figure~\ref{Olda1}), the operad $J$ is an
extension of the $A_\infty$-operad ${\mathcal A}_\infty$ of cellular
chains of ${\sf K}$. This explains the title of the paper.

Let $\uG_*(n) \subset G_*(n)$ be a graded abelian group\footnote{As
in~\cite{markl-shnider-stasheff:book}, underlining indicates the
non-$\Sigma$ version of an object.}
generating $G_*(n)$
as a free graded $\Sigma_n$-module -- one such a specific generating space
will be described on page~\pageref{1}. It follows from standard facts that
$J(n)$ is the free $\Sigma_n$-module generated by the $n$th piece
$\uFree(\uG)(n)$ of the free non-$\Sigma$ operad $\uFree(\uG)$\label{jarca!!}
generated by $\uG$. In particular, $(J(n),\pa)$ is, for each $n \geq
1$, a $\Sigma_n$-free resolution of the trivial $\Sigma_n$-module
$\ZZZ = \Com(n)$.

Since $J$ is cofibrant, for an arbitrary dg $E_\infty$-operad ${\mathcal
E}$ (as the Barratt-Eccles operad, surjection operad,
Eilenberg-Zilber operad, see~\cite{berger-fresse,mcclure-smith:JAMS03}, etc.), 
there exist an operadic morphism $J \to {\mathcal E}$ that lifts the
identity endomorphism of the operad $\Com$. 
If the ground ring is a field of characteristic zero, $J$ contains as
its deformation retract the minimal model ${\mathcal C}_\infty$ of the operad
$\Com$~\cite{markl:zebrulka} (operad ${\mathcal C}_\infty$ describes
$C_\infty$, also called commutative or balanced $A_\infty$, algebras).

In Section~\ref{Tamarkin} we analyze two particular bad cells and
calculate the differential of the corresponding generators of $G$. The
first one is the famous $6$-dimensional Tamarkin cell, the second is a
simpler $4$-dimensional bad cell whose existence was a surprise for
us.  The dimension of these two bad cells is precisely the critical
one, i.e.~they are bad cells of the lowest dimension in a given arity.
Although we did not give a general formula, it will be clear that our
construction of the differential applies to any bad cell in the
critical dimension.  Our second main theorem, proved in
Subsection~\ref{b}, says that this formula can always be
extended to a differential on $\Free(G)$.

\begin{theoremB}
Any formula for the differential of (one or more) bad cells in the critical
dimension\footnote{See Definition~\ref{JarusKa}.}
extends to a differential with the properties stated in Theorem~A.
\end{theoremB}

General bad cells are analyzed in Section~\ref{jaRcA}.
There is an obvious question about finding a {\em closed formula\/}
for the differential of Theorem~A that would apply to bad cells of
arbitrarily high dimensions. It is clear from our analysis of bad
cells that this formula should involve cellular approximations of
the images of the source-target conditions.  Theorem~6.1
of~\cite{markl-shnider:TAMS06} shows that any such an approximation is
associative only up to a hierarchy of higher homotopies. All these
homotopies must be build in the formula for the differential as
appropriate `correction terms.'  Writing such a
formula is far beyond our abilities.

\noindent
{\bf Acknowledgement.}  The paper grew up from conversations with
Michael Batanin. I was helped by Benoit Fresse who explained to me the
precise relation between his work and this note. I~am also indebted to
C.~Berger,
Ezra Getzler, Jim Stasheff, Rainer Vogt and Sasha Voronov for reading
the preliminary version of this note and many useful comments.

\noindent 
{\bf Glossary of notation} is given on
page~\pageref{gloss}. 

\section{Trees, barcodes and the $E_\infty$-operad $J = (\Free(G),\pa)$}
\label{JARCA}

In this section we describe the graded $\Sigma$-module $G =
\{G_*(n)\}_{n \geq 2}$ generating the operad $J$, together with the
linear part $\pal$ of the differential. 
Let us start by recalling some definitions
of~\cite{batanin}.  We denote, as usual, by $[k]$ the ordered set $1 <
2 < \cdots < k$.

\begin{definition}
\label{dnes_snad_s_Jaruskou}
Let $h \geq 1$. A {\em tree of height $h$\/} (or {\em tree with $h$
levels\/}, or {\em $h$-tree\/}) is a sequence of order preserving
maps
\begin{equation}
\label{eq:1}
T = [k_h]
\stackrel{\rho_{h-1}}{\vlra} 
[k_{h-1}]
\stackrel{\rho_{h-2}}{\vlra}
\cdots
\stackrel{\rho_{0}}{\vlra} [1].  
\end{equation}
We are not going to consider degenerate trees, so we assume that all
$k_m \geq 1$, for $0 \leq m \leq h$.
\end{definition}

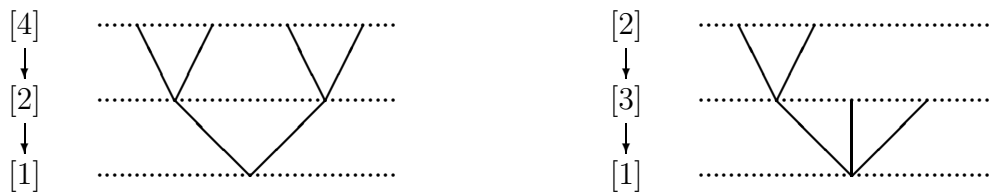
\begin{figure}[ht]
\begin{center}
\unitlength1cm
\begin{picture}(10.00,2.20)(1.00,0.00)
%
\multiput(0,0)(0,1){3}{
\multiput(0,0)(0.1,0){40}{\makebox(0,0){$\cdot$}}
}
\thinlines
\put(-1,1){\makebox(0,0){$[2]$}}
\put(-1,2){\makebox(0,0){$[4]$}}
\put(-1,0){\makebox(0,0){$[1]$}}
\put(-1,1.7){\vector(0,-1){.4}}
\put(-1,0.7){\vector(0,-1){.4}}
\thicklines
\put(2,0){\line(1,1){1}}
\put(2,0){\line(-1,1){1}}
\put(1,1){\line(1,2){.5}}
\put(1,1){\line(-1,2){.5}}
\put(3,1){\line(1,2){.5}}
\put(3,1){\line(-1,2){.5}}
%
\put(8,0)
{
\multiput(0,0)(0,1){3}{
\multiput(0,0)(0.1,0){40}{\makebox(0,0){$\cdot$}}
}
\put(-1,1){\makebox(0,0){$[3]$}}
\put(-1,2){\makebox(0,0){$[2]$}}
\put(-1,0){\makebox(0,0){$[1]$}}
\thinlines
\put(-1,1.7){\vector(0,-1){.4}}
\put(-1,0.7){\vector(0,-1){.4}}
\thicklines
\put(2,0){\line(1,1){1}}
\put(2,0){\line(-1,1){1}}
\put(2,0){\line(0,1){1}}
\put(1,1){\line(1,2){.5}}
\put(1,1){\line(-1,2){.5}}
}
\end{picture}
\end{center}
\caption{\label{fig:1}
Example of trees of height two. The left tree is pruned, the right one
is not. Arity of the left tree is $4$, arity of the right one is $2$.}
\end{figure}

A {\em vertex\/} of height $m$, $0 \leq m \leq h$, is an element of
$[k_m]$. One may imagine that each vertex $i \in [k_m]$, $m \geq 1$,
determines the oriented {\em edge\/}  that starts at $i$ and ends at
the vertex $\rho_{m-1}(i) \in [k_{m-1}]$. With this intuition, one may
indeed interpret objects of Definition~\ref{dnes_snad_s_Jaruskou} as
planar directed trees with vertices arranged at $h+1$ horizontal
lines shown in Figure~\ref{fig:1}.
A {\em leaf of height $m$\/} is
a vertex $i \in [k_m]$ which is not in the image of $\rho_m$. A {\em
tip\/} is a leaf of maximal height $h$. The {\em arity\/} of $T$ is
then the number of tips.  The tree is {\em pruned\/} if all its
leaves are tips. These definitions should be clear from
Figure~\ref{fig:1}.

\begin{figure}[ht]
\begin{center}
\unitlength1cm
\begin{picture}(10.00,3.00)(1.00,0.00)
\multiput(0,0)(0,1){4}{
\multiput(0,0)(0.1,0){40}{\makebox(0,0){$\cdot$}}
}
\thicklines
\put(0,1){
\put(2,0){\line(1,1){1}}
\put(2,0){\line(0,1){1}}
\put(2,0){\line(-1,1){1}}
\put(1,1){\line(1,2){.5}}
\put(1,1){\line(-1,2){.5}}
\put(3,1){\line(1,2){.5}}
\put(3,1){\line(-1,2){.5}}
\multiput(-.03,0)(.01,0){7}{
\put(2,0){\line(0,-1){1}}
}
}
\put(8,-1){
\multiput(0,0)(0,1){3}{
\multiput(0,1)(0.1,0){40}{\makebox(0,0){$\cdot$}}
}
\put(0,1){
\put(2,0){\line(1,1){1}}
\put(2,0){\line(-1,1){1}}
\put(1,1){\line(1,2){.5}}
\put(1,1){\line(-1,2){.5}}
\put(3,1){\line(1,2){.5}}
\put(3,1){\line(-1,2){.5}}
}
}
\end{picture}
\end{center}
\caption{\label{fig:2}
The left tree has a trunk (bold edge) 
and is not pruned. The right tree is its maximal reduced subtree.}
\end{figure}
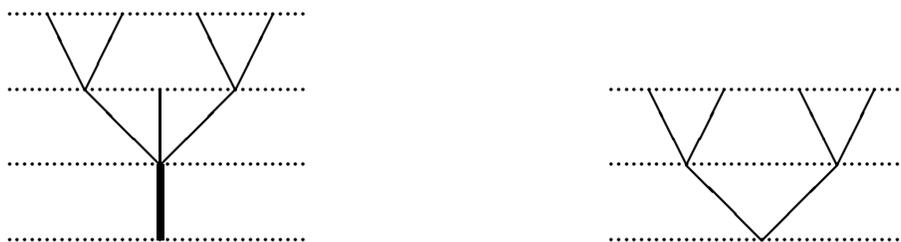

We say that a tree $T$ as in~(\ref{eq:1}) has a {\em trunk\/} if $k_m
= 1$ for some $m \geq 1$. A {\em trunk\/} of $T$ is then everything
`below' $k_m$; see Figure~\ref{fig:2}.

 We say that a tree is {\em
reduced\/}\label{jaRca}\footnote{%
Our terminology is not a standard one -- reduced usually 
means no vertices of arity~$1$.}
if it is pruned and if it has no trunk. Obviously, for each
$T$ there exists a unique maximal reduced subtree $r(T)$ of maximal
height. See again Figure~\ref{fig:2} -- the right tree is obtained
from the left one by first cutting off the trunk and 
then pruning the remaining shrub.
So pruning is for us cutting out branches that do not end in tips as
opposed to what one does in garden, namely cutting of those that stick
up too high.
Finally, for a tree $T$ as in~(\ref{eq:1}) define its {\em
dimension\/} $\dim(T)$ as 
\begin{equation}
\label{eq:2}
\dim(T) := e(T) - h - 1,
\end{equation}
where $e(T)$ is the number of edges and $h$ the height.

The {\em terminal tree\/} $U_h$ is the tree with all 
$k_m=1$. Terminal trees play a very special r\^ole
and, unless stated otherwise, we will not consider them. If
necessary, we set $\dim(U_h) := 0$ (formula~(\ref{eq:2})
would give $\dim(U_h) = -1$). We also define $r(U_h) := U_h$.

\begin{notation}
\label{S_Jarkou_ve_VO}
{\rm 
Denote by \label{prun}$\uTree^h(n) = \bigcup_{i \geq 0}
\uTree^h_{\hskip .2em i}(n)$ the graded set whose $i$th component
consists of pruned $h$-trees of dimension $i$, with $n$ tips, $h,n
\geq 1$.  We also denote \label{445}$\Tree^h(n) = \bigcup_{i \geq 0}
\Tree^h_i(n)$ the graded set of {\em labeled\/} pruned trees of
height $h$.  Elements of $\Tree^h_i(n)$ are couples $\bfT =(T,\ell)$,
where $T \in \uTree^h_{\hskip .2em i}(n)$ is as in~(\ref{eq:1}) and $\ell$ an
isomorphism (labeling) $\ell : [k_h] \stackrel\cong\to [n]$. The
symmetric group $\Sigma_n$ freely acts on $\Tree^h(n)$ by relabeling
the tips.

One has the inclusion $\uTree^h(n) \subset \Tree^h(n)$ given by
numbering the legs of an unlabeled planar tree from the left to the
right.  The subset $\uTree^h(n)$ freely generates
$\Tree^h(n)$ as a right graded $\Sigma_n$-set.

There are the {\em suspensions\/} $\underline{s} : \uTree^h(n)
\hookrightarrow \uTree^{h+1}(n)$ resp.~$s : \Tree^h(n)
\hookrightarrow \Tree^{h+1}(n)$ that adjoin to a (labeled) tree
a trunk of hight one.\label{jaruska-pusinka} 
The graded sets \label{098}$\uTreeinf(n)  := \dirlim \uTree^h(n)$
\label{087}resp.~$\Treeinf(n) := \dirlim \Tree^h(n)$ clearly consist of 
(labeled) reduced trees of an arbitrary height.
}\end{notation}

The first step towards our definition of the operad $J$ is: 

\begin{definition}
\label{JArKa}
The $n$th component $G_*(n)$ of the graded $\Sigma$-module $G_* = \{G_*(n)\}_{n
\geq 2}$ generating the operad $J$ is the free graded abelian group
$\Span(\Treeinf_*(n))$ spanned by the graded $\Sigma_n$-set
$\Treeinf(n)$ of labeled reduced trees with $n$ tips.
\end{definition}

\begin{figure}[ht]
\begin{center}
\unitlength.5cm
\begin{picture}(10.00,22.50)(6.8,-12.00)
\thicklines
\put(0,10){
$\dim = 0$:
\put(1,-.5){
\multiput(0,0)(0,1){2}{
\multiput(0.2,0)(0.2,0){14}{\makebox(0,0){$\cdot$}}
}
\put(1.5,0){\line(1,1){1}}
\put(1.5,0){\line(-1,1){1}}
\put(1.5,-1){\makebox(0,0){$[1|2]$}}
}
}

\put(0,6){
$\dim = 1$:
\put(1,-.5){
  \multiput(0,0)(0,1){2}{
  \multiput(.2,0)(0.2,0){14}{\makebox(0,0){$\cdot$}}
  }
  \put(1.5,0){\line(1,1){1}}
  \put(1.5,0){\line(0,1){1}}
  \put(1.5,0){\line(-1,1){1}}
\put(1.5,-1){\makebox(0,0){$[1|2|3]$}}
}
\put(5,-.5){
  \multiput(0,0)(0,1){3}{
  \multiput(0.2,0)(0.2,0){14}{\makebox(0,0){$\cdot$}}
  }
  \put(1.5,0){\line(1,1){1}}
  \put(1.5,0){\line(-1,1){1}}
  \put(.5,1){\line(0,1){1}}
  \put(2.5,1){\line(0,1){1}}
\put(1.5,-1){\makebox(0,0){$[1||2]$}}
}
}

\put(0,1){
$\dim = 2$:
\put(1,-.5){
  \multiput(0,0)(0,1){2}{
  \multiput(0.2,0)(0.2,0){14}{\makebox(0,0){$\cdot$}}
  }
  \put(1.5,0){\line(1,1){1}}
  \put(1.5,0){\line(-1,3){.3333}}
  \put(1.5,0){\line(1,3){.3333}}
  \put(1.5,0){\line(-1,1){1}}
\put(1.5,-1){\makebox(0,0){$[1|2|3|4]$}}
}
\put(6,-.5){
  \multiput(-1,0)(0,1){3}{
  \multiput(0.2,0)(0.2,0){19}{\makebox(0,0){$\cdot$}}
  }
  \put(1.5,0){\line(1,1){1}}
  \put(1.5,0){\line(-1,1){1}}
  \put(.5,1){\line(-1,1){1}}
  \put(.5,1){\line(1,1){1}}
  \put(2.5,1){\line(0,1){1}}
\put(1,-1){\makebox(0,0){$[1|2||3]$}}
}
\put(10.5,-.5){
  \multiput(-.5,0)(0,1){3}{
  \multiput(.20,0)(0.2,0){19}{\makebox(0,0){$\cdot$}}
  }
  \put(1,0){\line(1,1){1}}
  \put(1,0){\line(-1,1){1}}
  \put(2,1){\line(-1,1){1}}
  \put(2,1){\line(1,1){1}}
  \put(0,1){\line(0,1){1}}
\put(1.5,-1){\makebox(0,0){$[1||2|3]$}}
}
\put(15,-.5){
  \multiput(0,0)(0,1){4}{
  \multiput(0.2,0)(0.2,0){14}{\makebox(0,0){$\cdot$}}
  }
  \put(1.5,0){\line(1,1){1}}
  \put(1.5,0){\line(-1,1){1}}
  \put(.5,1){\line(0,1){2}}
  \put(2.5,1){\line(0,1){2}}
\put(1.5,-1){\makebox(0,0){$[1|||2]$}}
}
}

\put(0,-4.5){
$\dim = 3$:
\put(1,-.5){
\multiput(0,0)(0,1){2}{
\multiput(.20,0)(0.2,0){14}{\makebox(0,0){$\cdot$}}
}
\put(1.5,0){\line(1,1){1}}
\put(1.5,0){\line(-1,1){1}}
\put(1.5,0){\line(0,1){1}}
\put(1.5,0){\line(-1,2){.5}}
\put(1.5,0){\line(1,2){.5}}
\put(1.5,-1){\makebox(0,0){$[1|2|3|4|5]$}}
}

\put(6,-.5){
  \multiput(-1,0)(0,1){3}{
  \multiput(0.2,0)(0.2,0){18}{\makebox(0,0){$\cdot$}}
  }
  \put(1,0){\line(1,1){1}}
  \put(1,0){\line(-1,1){1}}
  \put(2,1){\line(-1,2){.5}}
  \put(2,1){\line(1,2){.5}}
  \put(0,1){\line(-1,2){.5}}
  \put(0,1){\line(1,2){.5}}
\put(1,-1){\makebox(0,0){$[1|2||3|4]$}}
}
\put(11,-.5){
  \multiput(-1,0)(0,1){3}{
  \multiput(0.2,0)(0.2,0){18}{\makebox(0,0){$\cdot$}}
  }
  \put(1.5,0){\line(1,1){1}}
  \put(1.5,0){\line(-1,1){1}}
  \put(.5,1){\line(-1,1){1}}
  \put(.5,1){\line(1,1){1}}
  \put(.5,1){\line(0,1){1}}
  \put(2.5,1){\line(0,1){1}}
\put(1,-1){\makebox(0,0){$[1|2|3||4]$}}
}
\put(15.5,-.5){
  \multiput(-1,0)(0,1){3}{
  \multiput(0.7,0)(0.2,0){18}{\makebox(0,0){$\cdot$}}
  }
  \put(1,0){\line(1,1){1}}
  \put(1,0){\line(-1,1){1}}
  \put(2,1){\line(-1,1){1}}
  \put(2,1){\line(0,1){1}}
  \put(2,1){\line(1,1){1}}
  \put(0,1){\line(0,1){1}}
\put(1.5,-1){\makebox(0,0){$[1||2|3|4]$}}
}
\put(1,-6.5){
  \multiput(0,0)(0,1){3}{
  \multiput(0.2,0)(0.2,0){14}{\makebox(0,0){$\cdot$}}
  }
  \put(1.5,0){\line(1,1){1}}
  \put(1.5,0){\line(0,1){1}}
  \put(1.5,0){\line(-1,1){1}}
  \put(.5,1){\line(0,1){1}}
  \put(1.5,1){\line(0,1){1}}
  \put(2.5,1){\line(0,1){1}}
\put(1.5,-1){\makebox(0,0){$[1||2||3]$}}
}
\put(6,-6.5){
  \multiput(-1,0)(0,1){4}{
  \multiput(0.2,0)(0.2,0){19}{\makebox(0,0){$\cdot$}}
  }
  \put(1.5,0){\line(1,1){1}}
  \put(1.5,0){\line(-1,1){1}}
  \put(.5,1){\line(0,1){1}}
  \put(2.5,1){\line(0,1){2}}
  \put(.5,2){\line(1,1){1}}
  \put(.5,2){\line(-1,1){1}}
\put(1,-1){\makebox(0,0){$[1|2|||3]$}}
}
\put(10,-6.5){
  \multiput(0,0)(0,1){4}{
  \multiput(0.2,0)(0.2,0){19}{\makebox(0,0){$\cdot$}}
  }
  \put(1.5,0){\line(1,1){1}}
  \put(1.5,0){\line(-1,1){1}}
  \put(2.5,1){\line(0,1){1}}
  \put(.5,1){\line(0,1){2}}
  \put(2.5,2){\line(1,1){1}}
  \put(2.5,2){\line(-1,1){1}}
\put(2,-1){\makebox(0,0){$[1|||2|3]$}}
}
\put(15,-6.5){
  \multiput(0,0)(0,1){5}{
  \multiput(0.2,0)(0.2,0){14}{\makebox(0,0){$\cdot$}}
  }
  \put(1.5,0){\line(1,1){1}}
  \put(1.5,0){\line(-1,1){1}}
  \put(.5,1){\line(0,1){3}}
  \put(2.5,1){\line(0,1){3}}
\put(1.5,-1){\makebox(0,0){$[1||||2]$}}
}
}
\end{picture}
\end{center}
\caption{\label{fig:3}
Complete list of reduced trees up to dimension $3$ and their barcodes.}
\end{figure}
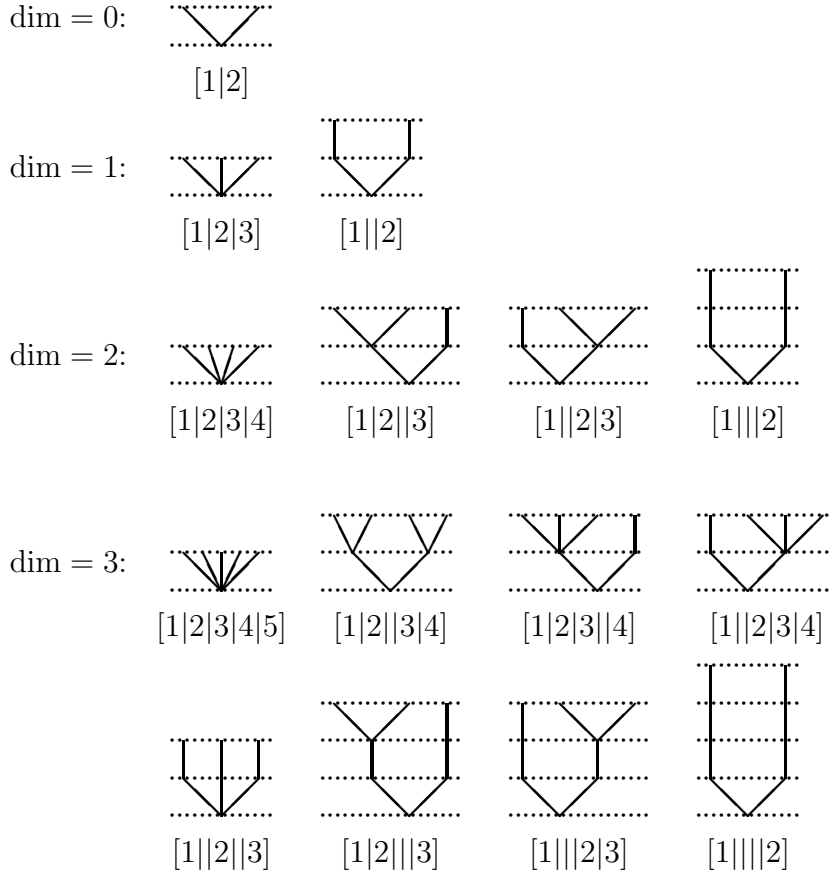

Each $G_*(n)$ is clearly a free $\Sigma_n$-module\label{1}
$\Sigma_n$-generated by the graded abelian group $\uG_*(n) =
\Span(\uTreeinf_*(n))$ spanned by (unlabeled) reduced trees.  A
complete list of unlabeled reduced trees up to dimension $3$ is given
in Figure~\ref{fig:3}; there are exactly $2^d$ reduced trees of
dimension~$d$.

Observe that $G_*(n) = \dirlim G^h_*(n)$, where \label{147}$G^h_*(n) :=
\Span (\Tree^h_*(n))$.
There is a convenient ``barcode'' notation for the reduced labeled
trees (and therefore also the Fox-Neuwirth cells recalled
in Section~\ref{JArca}) introduced in~\cite{getzler-jones:preprint}:

\begin{definition}
\label{S_Jarkou_do_Brna.}
The {\em barcode\/} of a reduced labeled tree is the list of labels of
its tips, separated by the vertical bars whose number equals the
depth of the gaps between the tips.
\end{definition}

Since the tips of an unlabeled tree can be labeled by $1,2,\ldots$
in the increasing order from the left to the right, the barcodes can
be used for unlabeled trees as well. See again
Figure~\ref{fig:3}. 
The height of the corresponding reduced tree is the maximal number of
the adjacent bars, and the dimension is the number of vertical bars
minus $1$.

The shortest way to describe the differential $\pal$ on the
$\Sigma$-module $G = \{G(n)\}_{n \geq 2}$ is to identify this
$\Sigma$-module to a suitable dg-submodule of an iterated bar
construction.  Let \label{jhg}$\sfB(A)$ denote the bar construction of an
associative algebra $A$, i.e.~the tensor algebra \hbox{$\TTT (\uparrow
\hskip -.3em A)$} generated by the suspension \hbox{$\uparrow \hskip
-.3em A$} of the abelian group $A$, with the degree $-1$ differential
$\pa_\sfB$ induced by the multiplication of $A$. It is
classical~\cite[X.12]{maclane:homology} that, if $A$ is commutative,
the shuffle product of the tensor algebra makes $\sfB(A) =
(\hbox{$\TTT (\uparrow \hskip -.3em A)$},\pa_\sfB)$ a commutative
associative algebra, thus the bar construction can be iterated.

Let us denote, for $h \geq 1$, by \label{qwe}$\sfB^h(A)$ the $h$-th
iterate of $\sfB(-)$ applied to $A$. Since the natural inclusion
$\iota^h:\sfB^h_*(A) \hookrightarrow \sfB^{h+1}_*(A)$ is a degree $+1$
map, to have a natural grading on the direct limit we need to regrade
by putting
\label{ty}$\hatB_*^h(A)
:= \downarrow^{h+1} \sfB^h_*(A)$. The induced inclusion
$\hatB^h_*(A) \hookrightarrow \hatB^{h+1}_*(A)$ is
degree-preserving so one may take \label{3r} $\hatB^\infty(A) =
(\hatB^\infty(A),\pa^\infty_\sfB)$, the direct limit $\dirlim
\hatB^h(A)$ with the induced differential.

Consider the free abelian group $V$ spanned by $\Rada x1n$,
interpreted as a commutative algebra with the trivial
multiplication. Denote by $\hatB_{\rada 11}^h(V)$ the sub-dg
abelian group of $\hatB^h(V)$ spanned by monomials that contain
each basic element $\Rada x1n$ exactly once, with the obvious
right $\Sigma_n$-action given by relabeling. 
Finally, let $\hatB_{\rada 11}^\infty(V) := \dirlim \hatB_{\rada 11}^h(V)$.

As observed in~\cite{getzler-jones:preprint} and, in more general
setting, also in~\cite{fresse:iterated}, the graded abelian group
$\hatB_{\rada 11}^\infty(V)$ is isomorphic to the graded abelian group
$G_*(n)$ of Definition~\ref{JArKa}.
\label{T}
The isomorphism $\omega : G(n) \stackrel\cong\to \hatB_{\rada
11}^\infty(V)$ has the following inductive description.

Let $T \in \uTree^h(n)$ and $g_T$ the corresponding  
generator of $G(n)$. If $h = 1$, then $T$
is the $n$-corolla $\star_n$, i.e.~the $1$-tree 
with the barcode $[1|\ldots|n]$. In this
case we put
\[
\omega(g_T) :=
\downarrow^2\left(\sx{1} \ot \cdots  \ot \sx{n}\right) \in
\hatB_{\rada 11}^1(V) \subset \hatB_{\rada11}^\infty(V). 
\]
Assume $\omega(g_S)$ has been defined for all $S \in \uTree^k(n)$ with
$1 \leq k < h$, and has the property that $\omega(g_S)$ actually belongs to
the subspace $\hatB_{\rada 11}^k(V) \subset \hatB_{\rada
11}^\infty(V)$. Let $T \in \uTree^h(n)$. There obviously exists some
$u$, $1 \leq u \leq n$, such that $T$ is obtained by grafting
pruned, not necessarily reduced, $(h-1)$-trees $\Rada T1u$ at the
tips of the $u$-corolla $\star_u$:
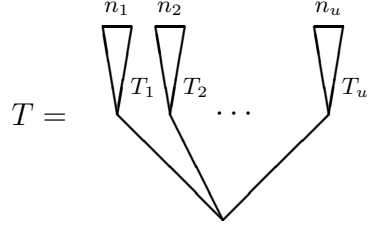
\begin{center}
{
\unitlength=1.000000pt
\begin{picture}(80.00,80.00)(0.00,5.00)
\thicklines
\put(80.00,80.00){\makebox(0.00,0.00){\scriptsize $n_u$}}
\put(20.00,80.00){\makebox(0.00,0.00){\scriptsize $n_2$}}
\put(0.00,80.00){\makebox(0.00,0.00){\scriptsize $n_1$}}
\put(85.00,50){\makebox(0.00,0.00)[l]{\scriptsize $T_u$}}
\put(25.00,50.00){\makebox(0.00,0.00)[l]{\scriptsize $T_2$}}
\put(5.00,50.00){\makebox(0.00,0.00)[l]{\scriptsize $T_1$}}
\put(80.00,40.00){\makebox(0.00,0.00)[b]{$\tria$}}
\put(20.00,40.00){\makebox(0.00,0.00)[b]{$\tria$}}
\put(0.00,40.00){\makebox(0.00,0.00)[b]{$\tria$}}
\put(45.00,40.00){\makebox(0.00,0.00){$\cdots$}}
\put(20.00,40.00){\line(1,-2){20.00}}
\put(40.00,0.00){\line(1,1){40.00}}
\put(0.00,40.00){\line(1,-1){40.00}}
\put(-40.00,40.00){\makebox(0.00,0.00)[l]{$T =$}}
\end{picture}}
\end{center}
where $\Rada n1u \geq 1$ with $n_1+\cdots+ n_u = n$ are the
arities of the trees $\Rada T1u$. For $i$, $1\leq i \leq u$, we denote
\[
V_i := \Span\{x_j;\ n_1 + \cdots + n_{i-1} +1  \leq j \leq n_1 + \cdots +
n_i\}.
\] 
We distinguish two cases.

(a) $n_i \geq 2$. Then let $R_i \in \uTree^{k_i}(n_i)$ be the
maximal reduced subtree of $T_i$. By induction, $\omega(g_{R_i}) \in 
\hatB_{\rada 11}^{k_i}(V_i)$ is defined and we put
$\omega_i \in \hatB_{\rada 11}^{h-1}(V_i)$ the image of $\omega(g_{R_i})$
under the natural inclusion $\hatB_{\rada 11}^{k_i}(V_i)
\hookrightarrow \hatB_{\rada 11}^{h-1}(V_i)$.
 
(b) $n_i = 1$. In this case, let $j := n_1 + \cdots + n_{i-1} +1$ and 
define $\omega_i  \in \hatB_{\rada 11}^{h-1}(V_i)$ the image of
$\downarrow^2 \hskip -.1em(\uparrow \hskip -.3em x_j) \in 
\hatB_{\rada 11}^{1}(V_i)$ under the
natural inclusion  $\hatB_{\rada 11}^{1}(V_i)
\hookrightarrow \hatB_{\rada 11}^{h-1}(V_i)$. 

\noindent 
Observe that, in both cases, $\uparrow^{h+1} \omega_i \in 
\hatB_{\rada 11}^h(V_i)$. Finally, let 
\[
\omega(g_T):=
\downarrow^{h+1}\left(\uparrow^{h+1} \omega_1 \ot \cdots \ot
\uparrow^{h+1} \omega_u\right)  \in \hatB_{\rada 11}^h(V) \subset
\hatB_{\rada 11}^\infty(V). 
\]
For example, 
\begin{equation*}
\begin{split}
\omega(g_\assoc \hskip -.2em) &= \downarrow^2 \left(\sx1\ot\sx2\ot\sx3\right) 
\in \hatB_{\rada11}^1(V) \subset \hatB_{\rada11}^\infty(V),
\\
\omega(g_\rp \hskip .1em) &=  \downarrow^3
\left(\uparrow(\sx1)\ot\uparrow(\sx2\ot\sx3)\right) 
\in \hatB_{\rada11}^2(V) \subset \hatB_{\rada11}^\infty(V), \mbox { \&c.}
\end{split}
\end{equation*}

\begin{definition}
\label{jarka}
The differential $\pal$ on $G_*(n)$ is defined by $\pal := \omega^{-1}
\circ \pabar \circ \omega$. Thus $\pal$ is the unique
differential such that $\omega : (G_*(n),\pal) \to 
(\hatB_{\rada 11}^\infty(V),\pabar)$ is an isomorphism of dg-abelian groups.
\end{definition}

It is easy to see that $(G_*(2),\pal)$ is the cellular chain complex of
the `globular' decomposition of the infinite sphere
${\mathbb{S}^{\infty}}$. The piece $G_*(3)$ contains Stasheff's associator,
two Mac Lane's hexagons (right and left), etc.

\section{Relation to the compactification of the configuration space}
\label{JArca}

Let \label{wsx}$\osfF_h(n) := \Conf hn/\Aff h$ be the moduli space
of configurations of $n$ distinct points in the $h$-dimensional
Euclidean plane $\RRR^h$, modulo the action of the affine group of
dilatations and translations.
Getzler and Jones in~\cite{getzler-jones:preprint} described a compactification
\label{edc}$\sfF_h(n)$ of $\osfF_h(n)$ such that 
the $\Sigma$-space $\sfF_h := \{\sfF_h(n)\}_{n\geq 1}$ 
is an operad in the category of manifolds with
corners, see also~\cite{markl:cf}. 

Will will in fact be interested
in the colimit \label{tgb}$\sfF := \dirlim \sfF_h$ which inherits
an operad structure from its constituents $\sfF_h$. As observed
in~\cite{getzler-jones:preprint}, the Fox-Neuwirth cells (recalled below) of
the open part $\osfF_h(n)$ induce a decomposition of the colimit
$\osfF := \dirlim \osfF_h$ which in turn induces a cell decomposition
of the compactification $\sfF$ compatible with the operad structure.

There was a belief that the operad $J$ can be easily read off from the
combinatorics of this cell decomposition, including a formula for the
differential. This turned out not to be the case because, due to the
existence of `bad' cells, the cell structure of $\sfF$ is not
regular. We identify, in Propositions~\ref{Jarka_byla_u_mne!}
and~\ref{JArkA}, the dimensions in which bad cells exist and analyze
them further in Section~\ref{jaRcA}. We start with:

\begin{definition}
\label{JarusKa}
The {\em critical dimension\/} $\dcrit(n)$ is defined as
\[
\dcrit(n) = \cases{\infty}{if $n=2,3$ and}{n}{if $n \geq 4$.}
\]
\end{definition}

Let use denote by $\Greg_*$ the graded 
sub-$\Sigma$-module of $G$ spanned by the
cells whose dimension is less than the critical one.  
The results of this section are summarized in the
following statement which is a combination of
Propositions~\ref{Jarka_byla_u_mne!},~\ref{JarunA} and the results
of Subsection~\ref{Zitra_zubar}.

\begin{proposition}
The combinatorics of the Fox-Neuwirth decomposition of the
compactification~$\sfF$ determines, by formula~(\ref{hyrim_barvami})
on page~\pageref{hyrim_barvami}, a differential $\pa$ on the free
operad $\Free(\Greg)$.
\end{proposition}

We start by recalling, following
closely~\cite{getzler-jones:preprint}, a correspondence between pruned
trees and flags of pre-orders. This point of view will be useful in
describing decompositions of configuration spaces.

\begin{definition}
A {\em pre-order\/} $\pi$ on a non-empty set $S$ is a reflexive transitive
relation $\leq$ such that if $a,b \in S$, either $a \leq b$ or $b \leq a$.
\end{definition}

A pre-order defines an equivalence $\sim$ on $S$ by $a \sim b$ if and
only if $a \leq b$ and $b \leq a$, and induces a total order on the
quotient $S/\hskip -.2em \sim$. We denote $|\pi|$ the number of
equivalence classes. A pre-order $\pi$ is {\em trivial\/} if $a \leq b$
of all $a,b \in S$ or, equivalently, if $|\pi|=1$.  Pre-orders on $S$
form a poset: $\pi_1 \prec \pi_2$ if $a \leq_2 b$ implies $a \leq_1 b$
for all $a,b \in S$. The maximal elements of this poset are the total
orders of $S$, the unique minimal element is the trivial pre-order.

\begin{definition}
A {\em flag of pre-orders\/} on the set $S$ of height $h \geq 1$ is a
sequence $(\pi_1 \prec \cdots \prec \pi_h)$ of pre-orders on $S$ such
that $\pi_h$ is a total order of $S$. Such a flag is {\em reduced\/}
if $\pi_1$ is not the trivial preorder. 
\end{definition}

Let $\Flag^h(n) = \bigcup_{i \geq 0} \Flag^h_i(n)$ denote the graded
set whose $i$th component is formed by flags of
preorders of height $h$ on the set $\{\rada 1n\}$ satisfying  $i =
\sum_{s=1}^h |\pi_s| - h -1$.  We also denote
$\uFlag^h(n) = \bigcup_{i \geq 0} \uFlag^h_i(n)$ the graded subset of
flags of preorders $(\pi_1 \prec \cdots \prec \pi_h)$ such that
$\pi_h$ is the standard linear order of  $\{\rada 1n\}$.

The {\em suspension\/} $s : \Flag^h(n) \hookrightarrow
\Flag^{h+1}(n)$ resp.~$\underline{s} : \uFlag^h(n) \hookrightarrow 
\uFlag^{h+1}(n)$
extends a given flag from the left by the trivial
preorder. The graded sets $\Flag(n)  := \dirlim \Flag^h(n)$
resp.~$\uFlag(n) := \dirlim \uFlag^h(n)$ consist of 
reduced flags of an arbitrary height.
The next proposition relies on Notation~\ref{S_Jarkou_ve_VO}.

\begin{proposition}
\label{jarina}
For each $n,h \geq 1$, there are natural isomorphisms of graded sets
$\Tree^h(n) \cong \Flag^h(n)$ and $\uTree^h(n) \cong \uFlag^h(n)$
which induce isomorphisms of the colimits $\Tree(n) \cong \Flag(n)$
and $\uTree(n) \cong \uFlag(n)$.
\end{proposition}

\begin{proof}
Let $\bfT =(T,\ell) \in \Tree^h(n)$ be a labeled tree, i.e.~$T \in
\uTree^h(n)$ is as in~(\ref{eq:1}) and $\ell : [k_h] \stackrel\cong\to
[n]$ a labeling. Such a tree $\bfT$ defines a flag of preorders
$(\pi_1 \prec \cdots \prec \pi_h)$ as follows.

For $i,j \in \{\rada 1n\}$ we put $i \leq_h j$ if and only if
$\ell^{-1}(i) \leq \ell^{-1}(j)$. In other words, $\pi_h$ is the image of the
natural order on $\{\rada 1{k_h}\}$ under the isomorphism $\ell$.  For $1
\leq s < h$ we write $i \leq_s j$ if and only if
\[
\rho_{s} \rho_{s+1} \cdots \rho_{h-1}(\ell^{-1}(i)) \leq
\rho_{s} \rho_{s+1} \cdots \rho_{h-1}(\ell^{-1}(j)).
\]
It is easy to prove that the above correspondence is one to one,
induces an isomorphism of the colimits and restricts to isomorphisms
of the `underlined' versions $\uTree^h(n) \cong \uFlag^h(n)$ and
$\uTree(n) \cong \uFlag(n)$. It is also clear that for the flag
$(\pi_1 \prec \cdots \prec \pi_h)$ corresponding to a tree $\bfT =
(T,\ell)$ one has $e(T) = \sum_{s=1}^h|\pi_s|$, therefore,
by~(\ref{eq:2}),
\[
\dim(T) = \textstyle\sum_{s=1}^h|\pi_s| - h -1,
\]
so the above isomorphism are compatible with the gradings.
\end{proof}

\begin{convention}
\label{Dnes_snad_Jarca_prijde_ke_mne.}
{\rm
Given Proposition~\ref{jarina}, we will make no difference between
pruned trees and the corresponding flags of preorders. Thus, for
instance, a boldfaced $\bfT$ will denote both a labeled tree
$(T,\ell)$ and the corresponding flag $(\pi_1 \prec \cdots \prec \pi_h)$.
}
\end{convention}

Recall that $\Conf hn$ denotes the configuration space of $n$ distinct
labeled points $\Rada p1n$ in the $h$-dimensional affine space
$\RRR^h$, $n,h \geq 1$. It is an $hn$-dimensional oriented smooth
non-compact manifold whose points are monomorphisms $f
: \{\rada 1n\} \to \RRR^h$ given as $f(k) := p_k$, $1 \leq k \leq
n$. For such an $f$ and $1\leq s \leq h$, denote by $f_s$ the
composition of $f$ with the projection $\RRR^h \epi \RRR^s$, $(\Rada
x1h) \mapsto (\Rada x1s)$, to the first $s$ coordinates. We finally
denote $\pi^f_s$ the preorder on the set $\{\rada 1n\}$ given by the
pullback of the lexicographic order of $\RRR^s$ via $f_s$. In this
way, each monomorphism $f : \{\rada 1n\} \to \RRR^h \in \Conf hn$
determines a flag of preorders $\bfT^f = (\pi^f_1 \prec \cdots \prec
\pi^f_h)$.

Conversely, for a given tree\footnote{We are already using
Convention~\ref{Dnes_snad_Jarca_prijde_ke_mne.}.} 
$\bfT = (T,\ell) \in \Tree^h(n)$ define
\label{yhn}$[\bfT] := \{f \in \Conf hn;\ \bfT^f = \bfT\}$.
It is clear that $\Conf hn$ is the disjoint union
\begin{equation}
\label{FNe}
\Conf hn = 
\textstyle\bigcup_{\bfT \in \Tree^h(n)} [\bfT].
\end{equation}
Each $[\bfT]$ is an open  ball of dimension $e(T) =
\sum_{s=1}^h|\pi_s|$, therefore a tree $\bfT \in \Tree_i^h(n)$
determines a cell of dimension $i+h+1$. For $h=2$,~(\ref{FNe}) describes
the classical {\em Fox-Neuwirth\/}
decomposition~\cite{fox-neuwirth:MathScand62} generalized
in~\cite{getzler-jones:preprint} to arbitrary $h \geq 2$ .
One may assign an orientation
to $[\bfT]$ taking first the coordinates
$x_1$ of the equivalence class $\pi_1$ in the increasing order, next the
coordinates $x_2$ of the equivalence class $\pi_2$, also in the
increasing order, \&c. An example can be found
on page~\pageref{ukl} of Section~\ref{free_Lie}.

As we already indicated, we will need the {\em
moduli space\/}  $\osfF_h(n) := \Conf hn/\Aff h$, where the affine group
$\Aff h = \RRR^h \semidirect \RRR^\times_+$ acts by
translations and dilatations in the obvious manner. We denote the quotient of
the cell $[\bfT]$ modulo $\Aff h$ by \label{plm}$\mu[\bfT]$. It
is clear that $\mu[\bfT]$ is an open ball of
dimension $e(T)-h-1$. This explains formula~(\ref{eq:2}) for the
dimension of a tree. 
One has the disjoint decomposition
\begin{equation}
\label{Jarca_na_chalupe}
\osfF_h(n) = 
\textstyle\bigcup_{\bfT \in \Tree^h(n)} \mu[\bfT].
\end{equation}

Let us denote $\osfF_h$ the $\Sigma$-space $\osfF_h = \smod
\osfF2$. Fulton and MacPherson construct
in~\cite{fulton-macpherson:anm94} a compactification $\sfF_h = \smod
{\sfF_h}1$ of $\osfF_h$ such that $\sfF_h(n)$ is, for $n \geq 2$, a
smooth manifold with corners containing $\osfF_h(n)$ as its unique
open stratum. The $\Sigma$-space $\sfF_h$ is obtained by gluing the
free operad $\Free(\osfF_h)$. In particular,
decomposition~(\ref{Jarca_na_chalupe}) of $\osfF_h$ induces a
decomposition of the free operad $\Free(\osfF_h)$ which in turn
induces a CW-structure of $\sfF_h$ via the gluing map $\Free(\osfF_h)
\epi \sfF_h$.

This implies that the cells of $\sfF_h$ are indexed by the free set-operad
$\Free(\Tree^h)$.\label{Jarka_na_chalupe.} 
Since the pieces $\Tree^h(n)$ of the generating
$\Sigma$-set $\Tree^h = \smod {\Tree^h}2$ are freely
$\Sigma_n$-generated by the subset $\uTree^h(n) \subset \Tree^h(n)$,
the natural inclusion $\uFree(\uTree^h)
\hookrightarrow \Free(\Tree^h)$ induces, for each $n \geq 2$, the
isomorphism of right $\Sigma_n$-sets
\[
\Free(\Tree^h)(n) \cong {\rm Ind\/}_{\rm {\mathbf 1}_n}^{\Sigma_n}
\uFree(\uTree^h)(n) = \uFree(\uTree^h)(n)\times \Sigma_n,
\]
where ${\mathbf 1}_n$ denotes the trivial representation of $\Sigma_n$
and $\uFree(-)$ the free non-$\Sigma$ operad functor. We will abbreviate the above
display by
\begin{equation}
\label{Jarca!}
\Free(\Tree^h) \cong \uFree(\uTree^h)\times \Sigma.
\end{equation}
It follows from~(\ref{Jarca!}) and the structure theorem for free 
operads~\cite[Section~II.1.9]{markl-shnider-stasheff:book} that
\begin{equation}
\label{X}
\Free(\Tree^h)(n) =  
\textstyle\bigcup_{\tau \in \PTree(n)}\tau(\uTree^h),\ n \geq 2,
\end{equation}
where $\PTree(n)$ denotes the set of {\em planar rooted\/} trees whose each
vertex has at least two input edges\footnote{This is the usual meaning
of being reduced, compare the footnote on page~\pageref{jaRca}.}, with leaves
labeled by an isomorphism $\omega : \Leaf(\tau) \to \{n\}$,
see~\cite[Sect.~II.1.5]{markl-shnider-stasheff:book} for the
terminology and notation.
The trees in $\PTree(n)$ are different than the trees
in Section~\ref{JARCA} in that they do not have levels. 
The set $\tau(\uTree^h)$ in~(\ref{X}) is the cartesian product
\[
\tau(\uTree^h) := 
{\mbox {\LARGE $\times$}}_{v \in \Vert(\tau)} \uTree^h(\Int(v)),
\]
where $\Vert(\tau)$ is the set of vertices of $\tau$ and $\Int(v)$ the
number of input edges (= the arity) of a vertex $v$.
Informally,~(\ref{X}) means that the cells of $\sfF_h$ are indexed by
planar leaf-labeled rooted trees whose vertices are decorated by the
graded set $\uTree^h$ of pruned non-labeled $h$-trees. 

The inclusion $\RRR^h \hookrightarrow \RRR^{h+1}$, $(\Rada x1h)
\mapsto (0,\Rada x1h)$, induces an inclusion of $\Sigma$-spaces 
$\osfF_h \hookrightarrow \osfF_{h+1}$ so one can take the colimit
$\osfF := \dirlim \osfF_h$. Decomposition~(\ref{Jarca_na_chalupe})
induces the decomposition

\vskip -1.9em
\begin{equation}
\label{Jarca_na_chalupe.}
\osfF(n) = 
\textstyle\bigcup_{\bfT \in \Tree(n)} \mu[\bfT]
\end{equation}
with the cells indexed by {\em reduced\/} trees. The colimit $\sfF =
\dirlim \sfF_h$ is again obtained by
gluing the free operad $\Free(\osfF)$, so~(\ref{Jarca_na_chalupe.})
gives a decomposition of $\sfF$ with cells indexed by the free
set-operad~$\Free(\Tree)$.

At this stage we need to extend Definition~\ref{JarusKa} of the
critical dimension for finite $h$
by 
\begin{equation}
\label{Vochomurka}
\dcrit^h(n) := \cases\infty{if $n = 2,3$, or $n = 4,5$ and $h\leq 2$, or
$n \geq 6$ and $h=1$,}n{in the remaining cases.} 
\end{equation}
Clearly, $\dcrit(n) = \lim_{h \to \infty} \dcrit^h(n)$.
Let \label{qw}$\Treeregh_*$ be, for $h \geq 1$, 
the graded $\Sigma$-subset
of the graded $\Sigma$-set $\Tree^h_*$ consisting of reduced
trees of dimension less that the critical one, i.e.~the graded
$\Sigma$-set such that
\[
\Treeregh_i(n) := 
\cases{\Tree^h_i(n)}{if $i < \dcrit^h(n)$, and}
\emptyset{if $i \geq \dcrit^h(n)$.\rule{0em}{1.1em}}
\] 
Observe that $\Treeregjedna_* = \Tree^1_*$. We will also need the direct
limit \label{bv}$\Treereg_* := \dirlim \Treeregh_*$. Clearly
$\Treereg_*(n) = \Tree_*(n)$ if $n =2,3$ while, for $n \geq 4$,
\[
\Treereg_i(n) = 
\cases{\Tree_i(n)}{if $i < n$, and}\emptyset{if $i \geq n$.}
\] 

We will call, just for the purposes of this section, the trees in
$\Treeregh$ or $\Treereg$ the {\em regular\/} trees.  We also denote
by $\uTreeregh$ (resp.~ $\uTreereg$) the $\Sigma$-subset of $\Treeregh$
(resp.~$\Treereg$) of unlabeled regular trees,
i.e.~$\uTreeregh := \Treeregh \cap \uTree$ (resp.~$\uTreereg :=
\Treereg \cap \uTree$).

\begin{definition}
\label{rt}
The {\em regular skeleton\/} $\sfFreg$ of $\sfF$ 
is the union  of the cells of
the CW-complex $\sfF$ indexed by the suboperad
$\Free(\Treereg) \subset \Free(\Tree)$.
The {\em regular skeleton\/} $\sfFreg_h$ of $\sfF_h$ is the
intersection $\sfFreg \cap \sfF_h$.
\end{definition}

For $\Free(\Treereg)$ 
we have a formula similar to~(\ref{X}), i.e.
\[
\Free(\Treereg)(n) =  
\textstyle\bigcup_{\tau \in \PTree(n)}\tau(\uTreereg),\ n \geq 2,
\]
thus the cells of $\sfFreg$ are indexed by
planar rooted labeled trees with
vertices decorated by unlabeled reduced regular trees.
It is clear that $\sfFreg_h$ is the union of the cells indexed by the
suboperad $\Free(\Treeregh)$. It is equally obvious that 
the sub-$\Sigma$-modules $\sfFreg_h \subset
\sfF_h$ and $\sfFreg \subset
\sfF$ are suboperads and that $\sfFreg_1 = \sfF_1$. 
Let us recall the following standard

\begin{definition}
\label{JarunKa}
A CW-complex is {\em regular\/} if
(i) the attaching maps are homeomorphisms and
(ii) the boundary of each cell is a union of cells.
\end{definition}

We have the following correction
to~\cite[Lemma~5.11]{getzler-jones:preprint} whose proof we
postpone to Section~\ref{jaRcA}.

\begin{proposition}
\label{Jarka_byla_u_mne!}
The CW-structures of $\sfF_h$, $h \geq 1$, and of $\sfF$ 
are compatible with the operad
structures and the symmetric group acts freely on the cells. The
spaces $\sfF(n)$ are regular cell complexes if and only if $n=2$ or $3$. 
The complexes $\sfF_h(n)$
are regular if and only if 
\begin{itemize}
\item[(i)] 
$n = 2,3$ and $h$ arbitrary, or
\item[(ii)] 
$n\leq 5$ and $h \leq 2$, or
\item[(iii)] 
$n$ arbitrary and $h=1$.
\end{itemize}

The spaces $\sfF_h(n)$ and $\sfF(n)$ satisfy condition~(i) of
Definition~\ref{JarunKa} for arbitrary~$n$ and $h$.
The CW-subcomplexes $\sfFreg_h\subset \sfF_h$ and $\sfFreg\subset
\sfF$ are regular, with the cell structures 
compatible with the operad
structures and the symmetric group acting freely on the cells.
\end{proposition}

Observe that, by the definition of the critical dimension,
Proposition~\ref{Jarka_byla_u_mne!} says that the spaces $\sfF_h(n)$
(resp.~ $\sfF(n)$) are regular if and only if $\dcrit^h(n) = \infty$
(resp.~$\dcrit(n) = \infty$).  We call cells violating (ii) of
Definition~\ref{JarunKa} the {\em bad\/} cells. The following
statement provides a `coordinate-free' definition of the regular
skeleta.

\begin{proposition}
\label{JArkA}
For each $n \geq 4$, $h \geq 3$ or $n \geq 6$, $h \geq 2$, 
there exist a bad cell $e_n^h$ in the
open stratum $\osfF_h(n)$ whose dimension $\dim(e_n^h)$ equals
$\dcrit^h(n)$.
Likewise, for each $n \geq 4$ there exists a bad cell $e_n \subset
\osfF(n)$ such that $\dim(e_n) = \dcrit(n)$. 
The regular skeleta are therefore the maximal regular subcomplexes closed under
the operad structure.
\end{proposition}

\begin{figure}[ht]
\begin{center}
{
\unitlength=1.3pt
\begin{picture}(100.00,120.00)(0.00,0.00)
\put(50.00,120.00){\makebox(0.00,0.00){$\vdots$}}
\put(60.00,120.00){\makebox(0.00,0.00){$\vdots$}}
\put(70.00,120.00){\makebox(0.00,0.00){$\vdots$}}
\put(80.00,120.00){\makebox(0.00,0.00){$\vdots$}}
\put(10.00,120.00){\makebox(0.00,0.00){$\vdots$}}
\put(90.00,120.00){\makebox(0.00,0.00){$\vdots$}}
\put(90.00,110.00){\makebox(0.00,0.00){\scriptsize $\times$}}
\put(80.00,100.00){\makebox(0.00,0.00){\scriptsize $\times$}}
\put(80.00,110.00){\makebox(0.00,0.00){\scriptsize $\times$}}
\put(70.00,110.00){\makebox(0.00,0.00){\scriptsize $\times$}}
\put(70.00,100.00){\makebox(0.00,0.00){\scriptsize $\times$}}
\put(70.00,90.00){\makebox(0.00,0.00){\scriptsize $\times$}}
\put(60.00,80.00){\makebox(0.00,0.00){\scriptsize $\times$}}
\put(60.00,90.00){\makebox(0.00,0.00){\scriptsize $\times$}}
\put(60.00,100.00){\makebox(0.00,0.00){\scriptsize $\times$}}
\put(60.00,110.00){\makebox(0.00,0.00){\scriptsize $\times$}}
\put(50.00,70.00){\makebox(0.00,0.00){\scriptsize $\times$}}
\put(50.00,80.00){\makebox(0.00,0.00){\scriptsize $\times$}}
\put(50.00,90.00){\makebox(0.00,0.00){\scriptsize $\times$}}
\put(50.00,100.00){\makebox(0.00,0.00){\scriptsize $\times$}}
\put(50.00,110.00){\makebox(0.00,0.00){\scriptsize $\times$}}
\put(40.00,60.00){\makebox(0.00,0.00){\scriptsize $\times$}}
\put(40.00,70.00){\makebox(0.00,0.00){\scriptsize $\times$}}
\put(40.00,80.00){\makebox(0.00,0.00){\scriptsize $\times$}}
\put(40.00,90.00){\makebox(0.00,0.00){\scriptsize $\times$}}
\put(40.00,100.00){\makebox(0.00,0.00){\scriptsize $\times$}}
\put(40.00,110.00){\makebox(0.00,0.00){\scriptsize $\times$}}
\put(30.00,50.00){\makebox(0.00,0.00){\scriptsize $\times$}}
\put(30.00,60.00){\makebox(0.00,0.00){\scriptsize $\times$}}
\put(30.00,70.00){\makebox(0.00,0.00){\scriptsize $\times$}}
\put(30.00,80.00){\makebox(0.00,0.00){\scriptsize $\times$}}
\put(30.00,90.00){\makebox(0.00,0.00){\scriptsize $\times$}}
\put(30.00,100.00){\makebox(0.00,0.00){\scriptsize $\times$}}
\put(30.00,110.00){\makebox(0.00,0.00){\scriptsize $\times$}}
\put(-20,0){
\put(0.00,120.00){\makebox(0.00,0.00)[r]{dimension}}
\put(0.00,110.00){\makebox(0.00,0.00){\scriptsize$10$}}
\put(0.00,100.00){\makebox(0.00,0.00){\scriptsize$9$}}
\put(0.00,90.00){\makebox(0.00,0.00){\scriptsize$8$}}
\put(0.00,80.00){\makebox(0.00,0.00){\scriptsize$7$}}
\put(0.00,70.00){\makebox(0.00,0.00){\scriptsize$6$}}
\put(0.00,60.00){\makebox(0.00,0.00){\scriptsize$5$}}
\put(0.00,50.00){\makebox(0.00,0.00){\scriptsize$4$}}
\put(0.00,40.00){\makebox(0.00,0.00){\scriptsize$3$}}
\put(0.00,30.00){\makebox(0.00,0.00){\scriptsize$2$}}
\put(0.00,20.00){\makebox(0.00,0.00){\scriptsize$1$}}
\put(0.00,10.00){\makebox(0.00,0.00){\scriptsize$0$}}
}
\put(100.00,0.00){\makebox(0.00,0.00){\scriptsize $11$}}
\put(113.00,10.00){\makebox(0.00,0.00)[l]{arity}}
\put(90.00,0.00){\makebox(0.00,0.00){\scriptsize$10$}}
\put(80.00,0.00){\makebox(0.00,0.00){\scriptsize$9$}}
\put(70.00,0.00){\makebox(0.00,0.00){\scriptsize$8$}}
\put(60.00,0.00){\makebox(0.00,0.00){\scriptsize$7$}}
\put(50.00,0.00){\makebox(0.00,0.00){\scriptsize$6$}}
\put(40.00,0.00){\makebox(0.00,0.00){\scriptsize$5$}}
\put(30.00,0.00){\makebox(0.00,0.00){\scriptsize$4$}}
\put(20.00,0.00){\makebox(0.00,0.00){\scriptsize$3$}}
\put(10.00,0.00){\makebox(0.00,0.00){\scriptsize$2$}}
\put(0.00,0.00){\makebox(0.00,0.00){\scriptsize$1$}}
\put(-10.00,0.00){\makebox(0.00,0.00){\scriptsize$0$}}
\put(40.00,120.00){\makebox(0.00,0.00){$\vdots$}}
\put(30.00,120.00){\makebox(0.00,0.00){$\vdots$}}
\put(20.00,120.00){\makebox(0.00,0.00){$\vdots$}}
\put(90.00,90.00){\makebox(0.00,0.00){$\bullet$}}
\put(80.00,80.00){\makebox(0.00,0.00){$\bullet$}}
\put(70.00,70.00){\makebox(0.00,0.00){$\bullet$}}
\put(60.00,60.00){\makebox(0.00,0.00){$\bullet$}}
\put(50.00,50.00){\makebox(0.00,0.00){$\bullet$}}
\put(40.00,40.00){\makebox(0.00,0.00){$\bullet$}}
\put(30.00,30.00){\makebox(0.00,0.00){$\bullet$}}
\put(20.00,20.00){\makebox(0.00,0.00){$\bullet$}}
\put(10.00,10.00){\makebox(0.00,0.00){$\bullet$}}
\multiput(10,10)(0,10){11}{\makebox(0.00,0.00){$\bullet$}}
\multiput(20,20)(0,10){10}{\makebox(0.00,0.00){$\bullet$}}
\multiput(30,40)(10,10){7}{\makebox(0.00,0.00){$\bullet$}}
\put(7,30){\makebox(0.00,0.00)[r]{\scriptsize $h\!\! = \!\!3$}}
\put(7,40){\makebox(0.00,0.00)[r]{\scriptsize $h\!\! = \!\!4$}}
\put(0,50){\makebox(0.00,0.00)[b]{$\vdots$}}
\put(10,30){\line(1,3){27}}
\put(10,40){\line(1,4){17}}
\thicklines
\put(-10.00,5.00){\vector(0,1){115.00}}
\put(-15.00,10.00){\vector(1,0){125.00}}
\multiput(0,10)(10,0){11}{\makebox(0.00,0.00){\rule{.7pt}{5pt}}}
\multiput(-10,10)(0,10){11}{\makebox(0.00,0.00){\rule{5pt}{.7pt}}}
\end{picture}}
\end{center}
\caption{\label{Olda} 
The bad (marked $\times$) and regular (marked
$\bullet$) cells of $\hbox{$\protect\osfF$}$, and of
$\hbox{$\protect\osfF_h$}$ for $h \geq 3$.
The lines marked $h=3,4$ show the dimension of $\sfF_h(n)$.
}
\end{figure}
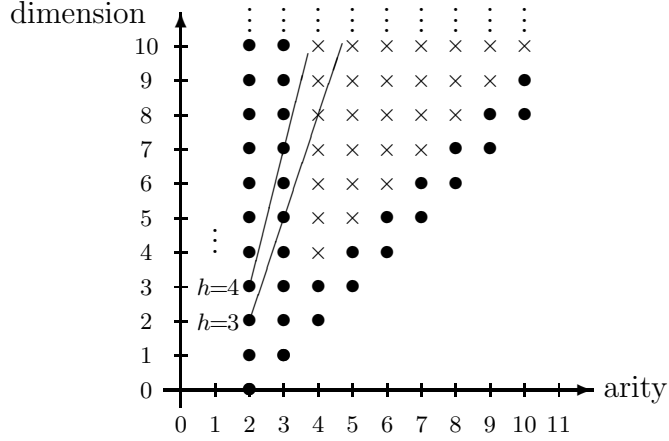

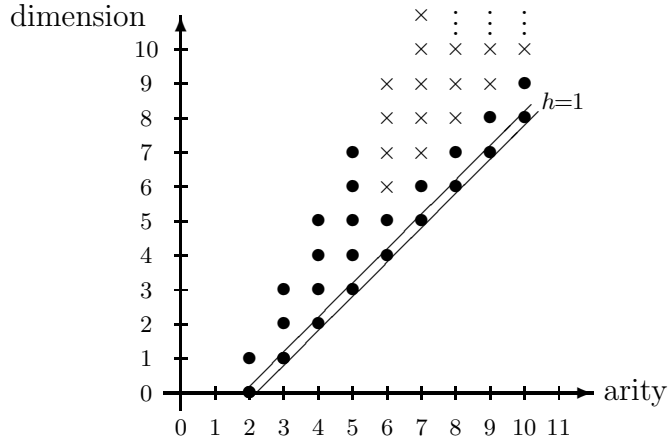
\begin{figure}[ht]
\begin{center}
{
\unitlength=1.3pt
\begin{picture}(100.00,120.00)(0.00,0.00)
\put(70.00,120.00){\makebox(0.00,0.00){$\vdots$}}
\put(80.00,120.00){\makebox(0.00,0.00){$\vdots$}}
\put(90.00,120.00){\makebox(0.00,0.00){$\vdots$}}
\put(90.00,110.00){\makebox(0.00,0.00){\scriptsize $\times$}}
\put(80.00,100.00){\makebox(0.00,0.00){\scriptsize $\times$}}
\put(80.00,110.00){\makebox(0.00,0.00){\scriptsize $\times$}}
\put(70.00,110.00){\makebox(0.00,0.00){\scriptsize $\times$}}
\put(70.00,100.00){\makebox(0.00,0.00){\scriptsize $\times$}}
\put(70.00,90.00){\makebox(0.00,0.00){\scriptsize $\times$}}
\put(60.00,80.00){\makebox(0.00,0.00){\scriptsize $\times$}}
\put(60.00,90.00){\makebox(0.00,0.00){\scriptsize $\times$}}
\put(60.00,100.00){\makebox(0.00,0.00){\scriptsize $\times$}}
\put(60.00,110.00){\makebox(0.00,0.00){\scriptsize $\times$}}
\put(60.00,120.00){\makebox(0.00,0.00){\scriptsize $\times$}}
\put(50.00,70.00){\makebox(0.00,0.00){\scriptsize $\times$}}
\put(50.00,80.00){\makebox(0.00,0.00){\scriptsize $\times$}}
\put(50.00,90.00){\makebox(0.00,0.00){\scriptsize $\times$}}
\put(50.00,100.00){\makebox(0.00,0.00){\scriptsize $\times$}}
\put(40.00,60.00){\makebox(0.00,0.00){$\bullet$}}
\put(40.00,70.00){\makebox(0.00,0.00){$\bullet$}}
\put(40.00,80.00){\makebox(0.00,0.00){$\bullet$}}
\put(30.00,50.00){\makebox(0.00,0.00){$\bullet$}}
\put(30.00,60.00){\makebox(0.00,0.00){$\bullet$}}
\put(-20,0){
\put(0.00,120.00){\makebox(0.00,0.00)[r]{dimension}}
\put(0.00,110.00){\makebox(0.00,0.00){\scriptsize$10$}}
\put(0.00,100.00){\makebox(0.00,0.00){\scriptsize$9$}}
\put(0.00,90.00){\makebox(0.00,0.00){\scriptsize$8$}}
\put(0.00,80.00){\makebox(0.00,0.00){\scriptsize$7$}}
\put(0.00,70.00){\makebox(0.00,0.00){\scriptsize$6$}}
\put(0.00,60.00){\makebox(0.00,0.00){\scriptsize$5$}}
\put(0.00,50.00){\makebox(0.00,0.00){\scriptsize$4$}}
\put(0.00,40.00){\makebox(0.00,0.00){\scriptsize$3$}}
\put(0.00,30.00){\makebox(0.00,0.00){\scriptsize$2$}}
\put(0.00,20.00){\makebox(0.00,0.00){\scriptsize$1$}}
\put(0.00,10.00){\makebox(0.00,0.00){\scriptsize$0$}}
}
\put(100.00,0.00){\makebox(0.00,0.00){\scriptsize $11$}}
\put(113.00,10.00){\makebox(0.00,0.00)[l]{arity}}
\put(90.00,0.00){\makebox(0.00,0.00){\scriptsize$10$}}
\put(80.00,0.00){\makebox(0.00,0.00){\scriptsize$9$}}
\put(70.00,0.00){\makebox(0.00,0.00){\scriptsize$8$}}
\put(60.00,0.00){\makebox(0.00,0.00){\scriptsize$7$}}
\put(50.00,0.00){\makebox(0.00,0.00){\scriptsize$6$}}
\put(40.00,0.00){\makebox(0.00,0.00){\scriptsize$5$}}
\put(30.00,0.00){\makebox(0.00,0.00){\scriptsize$4$}}
\put(20.00,0.00){\makebox(0.00,0.00){\scriptsize$3$}}
\put(10.00,0.00){\makebox(0.00,0.00){\scriptsize$2$}}
\put(0.00,0.00){\makebox(0.00,0.00){\scriptsize$1$}}
\put(-10.00,0.00){\makebox(0.00,0.00){\scriptsize$0$}}
\put(90.00,90.00){\makebox(0.00,0.00){$\bullet$}}
\put(80.00,80.00){\makebox(0.00,0.00){$\bullet$}}
\put(70.00,70.00){\makebox(0.00,0.00){$\bullet$}}
\put(60.00,60.00){\makebox(0.00,0.00){$\bullet$}}
\put(50.00,50.00){\makebox(0.00,0.00){$\bullet$}}
\put(40.00,40.00){\makebox(0.00,0.00){$\bullet$}}
\put(30.00,30.00){\makebox(0.00,0.00){$\bullet$}}
\put(20.00,20.00){\makebox(0.00,0.00){$\bullet$}}
\put(10.00,10.00){\makebox(0.00,0.00){$\bullet$}}
\put(11.00,9.00){\line(1,1){83.00}}
\put(9.00,11.00){\line(1,1){83.00}}
\multiput(10,10)(0,10){2}{\makebox(0.00,0.00){$\bullet$}}
\multiput(20,20)(0,10){3}{\makebox(0.00,0.00){$\bullet$}}
\multiput(30,40)(10,10){7}{\makebox(0.00,0.00){$\bullet$}}
\put(95,95){\makebox(0.00,0.00)[l]{\scriptsize $h\!\! = \!\!1$}}
\thicklines
\put(-10.00,5.00){\vector(0,1){115.00}}
\put(-15.00,10.00){\vector(1,0){125.00}}
\multiput(0,10)(10,0){11}{\makebox(0.00,0.00){\rule{.7pt}{5pt}}}
\multiput(-10,10)(0,10){11}{\makebox(0.00,0.00){\rule{5pt}{.7pt}}}
\end{picture}}
\end{center}
\caption{\label{Olda1} 
The bad (marked $\times$) and regular (marked
$\bullet$) cells of 
$\hbox{$\protect\osfF_2$}$.
The double line
represents the embedding of the top-dimensional cell of the associahedron 
${\sf K}$. 
}
\end{figure}

The first bad cell was found by D.~Tamarkin. We will call
this particular bad cell the Tamarkin cell and recall its definition
in Section~\ref{Tamarkin} in which we also prove
Propositions~\ref{Jarka_byla_u_mne!} and~\ref{JArkA}. The
case $h=1$ is special; $\sfF_1$ is the {\em Stasheff's
operad\/} of the associahedra~\cite{jds:hahI} which indeed forms a regular cell
complex. The dimensions/arities in which the bad cells of the open
strata sit are shown
in Figures~\ref{Olda} and~\ref{Olda1}.

Let us define the increasing filtration $\scrF_*(n) = \cdots
\scrF_0(n) \subset \scrF_1(n)\subset \scrF_2(n) \cdots$ of $\sfFreg(n)$
by
\[
\scrF_p(n) := \bigcup\{e \mbox { a cell of $\sfFreg(n)$};\ \dim(e) \leq p\}.
\]
Since the cells of $\sfFreg$ are indexed by the free operad $\Free(\Treereg)$,
this filtration is manifestly operadic, i.e.~each $\scrF_p(n)$ is
$\Sigma_n$-invariant and $\scrF_p(m) \circ_i \scrF_q(n) \subset
\scrF_{p+q}(m+n-1)$ for $m,n \geq 1$, $1 \leq i \leq m$. 
Each layer $(\scrE^r_{**}(n),\pa^r)$ of the induced
spectral sequence determines a dg-operad $\scrE^r :=
\{(\scrE^r_*(n),\pa^r)\}_{n \geq 1}$ with $\scrE^r_*(n) :=
\bigoplus_{* = p+q} \scrE^r_{pq}(n)$. 

The dg-operad $\scrE^2$ will be of a particular importance for us.
As usual, the abelian group $\scrE^2_{pq}(n)$ equals the reduced
homology $\overline{H}_{p+q}(\scrF_p(n)/\scrF_{p-1}(n))$, so 
the regularity of the CW-structure established in
Proposition~\ref{Jarka_byla_u_mne!} implies that
\[
\scrE^2_{pq}(n) = \cases{\Span(\mbox{set of $p$-dimensional
  cells of $\sfF(n)$})}{if $q=0$ and}0{if
  $q \not=0$.}
\]
In particular, $\scrE^2_{p0}$ contains, 
for each $\bfT \in \Treereg_p(n)$, 
the generator $c_\bfT$ corresponding to the cell~$\mu[\bfT]$.

Recall that $\Greg_* = \Span(\Treereg_*)$ is the graded
$\Sigma$-submodule of the $\Sigma$-module $G_*$ from
Definition~\ref{JArKa} spanned by the generators $g_\bfT$ indexed by
regular trees $\bfT\in \Treereg$.  We have a natural map of graded
operads $j: \Free(\Greg) \to \scrE^2$ given by $j(g_\bfT) := c_\bfT$
for $\bfT \in \Treereg$.

\begin{proposition}
\label{JarunA}
There is a unique differential $\pa$ on the free operad $\Free(\Greg)$
such that the map $j:  (\Free(\Greg),\pa) \to (\scrE^2,\pa^2)$ is a map of
dg-operads. Moreover, $\pa$ is the sum $\pal + \papert$ where $\pal$
is as in Definition~\ref{jarka}.
\end{proposition}

\begin{proof}
It is clear from the description of the cell structure of $\sfFreg$
via the free set-operad $\Free(\Treereg)$ that the map $j$ is an
isomorphism of graded operads, which implies the existence and
uniqueness of the differential $\pa$. The fact that $\pa$ constructed
in this way is a perturbation of the linear part $\pal$ of
Definition~\ref{jarka} will follow from explicit
calculations given below and Proposition~\ref{V_9:30_sraz_s_Jarkou_u_Pakousu.}.
\end{proof}

\noindent 
\subsection{The differential $\pa$ in sub-critical dimensions}
\label{Zitra_zubar}
Proposition~\ref{JarunA} translates the description of the
differential operad ${\Free(\Greg},\pa)$ into the
standard task of calculating the second term of the spectral sequence
associated to the regular cell complex $\sfFreg$.  Given an
$i$-dimensional cell $e$ of $\sfFreg$, one needs first to identify
cells forming the boundary of $e$. The differential of the generator
corresponding to $e$ then is then the sum of the generators
corresponding to $(i-1)$-dimensional cells in the boundary of $e$,
with the signs determined by the orientations.

In our particular case, the compatibility of the differential with the
operad structure implies that it suffices to describe the boundaries
of the cells $\mu[T]$ corresponding to the operadic generators in
$\Greg$, indexed by unlabeled regular reduced trees $T \in
\uTreereg$.  This was in fact already done in~\cite{batanin}, so we
only need to recall the necessary notions. Let us recollect the
notation first.

\begin{notation}
\label{Jarka_zavolala_sama_od_sebe!}
{\rm 
We introduced the following objects indexed by trees $\bfT \in
\Tree$ (resp.~the unlabeled versions $T \in \uTree$): the
corresponding generator $g_\bfT$ (resp.~$g_T$) of $G = \Span(\Tree)$
(resp.~of $\uG = \Span(\uTree)$), the Fox-Neuwirth cell $\mu[\bfT]$
(resp.~$\mu[T]$) of $\osfF$, and $c_\bfT$ (resp.~$c_T$) -- the
corresponding generator of $\scrE^2$. We will also denote by $E_T$ the
corresponding generator of the free set-operad $\Free(\Tree)$ and, for
an element $C \in \Free(\Tree)$, by $\mu[C]$ the corresponding cell of
$\sfF$.  
}\end{notation}

Let us return to our task of describing the differential
$\pa$. According to~\cite[Definition~2.2]{batanin}, a~{\em morphism of
$h$-trees\/}
\[
T = [k_h]
\stackrel{\rho_{h-1}}{\vlra} 
[k_{h-1}]
\stackrel{\rho_{h-2}}{\vlra}
\cdots
\stackrel{\rho_{0}}{\vlra} [1]
\] 
and 
\[
S = [s_h]
\stackrel{\xi_{h-1}}{\vlra} 
[s_{h-1}]
\stackrel{\xi_{h-2}}{\vlra}
\cdots
\stackrel{\xi_{0}}{\vlra} [1]
\]
is given by a sequence $\sigma = (\sigma_h,\ldots,\sigma_0)$ of not
necessary order preserving maps
$\sigma_m : [k_m] \to [s_m]$, $0 \leq m \leq h$, with the property 
that for each $m$
and each $j \in [k_{m-1}]$, the restriction of $\sigma_m$ to
$\rho^{-1}_{m-1}(j)$ preserves the induced 
order.\footnote{We believe
the same implicit notation for a permutation
and a~morphism of trees will not confuse the reader.}

Let $T \in \uTree$ be a reduced unlabeled $h$-tree as above. We will
consider maps $\sigma : T \to S$ of $h$-trees such that
\begin{itemize}
\item[(i)]
the tree $S$ is pruned, but possibly with a trunk, and
\item[(ii)]
the map $\sigma$ induces an epimorphism of tips, that is, $\sigma_h : [k_h]
\to [s_h]$ is onto.
\end{itemize}

We will call such a map $\sigma$ a {\em face\/} of the tree
$T$. Observe that $\sigma$ is determined by the values
$\sigma_h(i)$, $i \in [k_h]$. Let us explain how a  
face $\sigma$ determines a cell of $\sfF$ in
the boundary of $\mu[T]$.
We need first to describe, following again M.~Batanin's~\cite{batanin}, 
faces $\sigma$ in terms of fibers.

\begin{figure}[ht]
\begin{center}
\unitlength.65cm
\begin{picture}(10.00,2.50)
\thicklines
\put(0,0){
  \multiput(-1,0)(0,1){3}{
  \multiput(0.2,0)(0.1,0){37}{\makebox(0,0){$\cdot$}}
  }
\multiput(-.05,0)(.01,0){11}{
  \put(1.5,0){\line(1,1){1}}
  \put(1.5,0){\line(-1,1){1}}
  \put(.5,1){\line(1,1){1}}
}
  \put(2.5,1){\line(0,1){1}}
  \put(.5,1){\line(-1,1){1}}
\put(-0.5,2.2){\makebox(0,0)[b]{\scriptsize$1$}}
\put(1.5,2.2){\makebox(0,0)[b]{\scriptsize$2$}}
\put(2.5,2.2){\makebox(0,0)[b]{\scriptsize$3$}}
}
\put(5,0){
\multiput(0,0)(0,1){3}{
\multiput(0.2,0)(0.1,0){27}{\makebox(0,0){$\cdot$}}
}
\put(1.5,1){\line(-1,1){1}}
\multiput(-.03,0)(.01,0){7}{
\put(1.5,1){\line(1,1){1}}
\put(1.5,0){\line(0,1){1}}
}
\put(0.5,2.2){\makebox(0,0)[b]{\scriptsize$1$}}
\put(2.5,2.2){\makebox(0,0)[b]{\scriptsize$2$}}
}
\put(4,1){\makebox(0,0){$\stackrel{\sigma}{\longrightarrow}$}}
\end{picture}
\end{center}
\caption{\label{fig:4}
Fiber $F_2$ (shown in bold lines) of the map $\sigma : T \to S$ given
by $\sigma_2(1) = \sigma_2(3) = 1$ and $\sigma_2(2) = 2$.}
\end{figure}
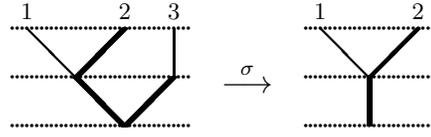

Let $\sigma : T \to S$ be a face of $T$ as above. 
For each tip $j \in [s_h]$, let $S_j$ be the path in $S$
connecting $\{j\}$ with the root of $S$. Then the {\em $j$th 
fiber of $\sigma$\/} is the subtree $F_j : = \sigma^{-1}(S_j)$ of
$T$. We believe that Figure~\ref{fig:4} elucidates this definition. 

\begin{figure}[ht]
\begin{center}
\unitlength.65cm
\begin{picture}(10.00,3.50)(3.5,.5)
\thicklines
\put(-2,0){
\put(1.5,1){\line(-1,1){1}}
\put(1.5,1){\line(1,1){1}}
\put(1.5,0){\line(0,1){1}}
\put(-1,2.2){\makebox(0,0)[t]{$D_1$:}}
\put(-.44,2.2){
\unitlength.4cm
  \put(1.5,0){\line(1,1){1}}
  \put(1.5,0){\line(-1,1){1}}
  \put(.5,1){\line(0,1){1}}
  \put(.5,2.2){\makebox(0,0)[b]{\scriptsize 1}}
  \put(2.5,1){\line(0,1){1}}
  \put(2.5,2.2){\makebox(0,0)[b]{\scriptsize 3}}
}
\put(1.57,2.2){
\unitlength.4cm
  \put(1.5,0){\line(1,1){1}}
  \put(1.5,0){\line(-1,1){1}}
  \put(.5,1){\line(0,1){1}}
  \put(.5,2.2){\makebox(0,0)[b]{\scriptsize 2}}
}
}
\put(5,0){
\put(1.5,1){\line(-1,1){1}}
\put(1.5,1){\line(1,1){1}}
\put(1.5,0){\line(0,1){1}}
\put(-2,2.2){\makebox(0,0)[t]{$D_2$:}}
\put(-.44,2.2){
\unitlength.4cm
  \put(1.5,0){\line(1,1){1}}
  \put(1.5,0){\line(-1,1){1}}
  \put(.5,1){\line(1,1){1}}
  \put(.5,1){\line(-1,1){1}}
  \put(-.5,2.2){\makebox(0,0)[b]{\scriptsize 1}}
  \put(1.5,2.2){\makebox(0,0)[b]{\scriptsize 2}}
}
\put(1.57,2.2){
\unitlength.4cm
  \put(1.5,0){\line(1,1){1}}
  \put(1.5,0){\line(-1,1){1}}
  \put(2.5,1){\line(0,1){1}}
  \put(2.5,2.2){\makebox(0,0)[b]{\scriptsize 3}}
}
}

\put(11,0){
\put(1.5,1){\line(-1,1){1}}
\put(1.5,1){\line(1,1){1}}
\put(1.5,1){\line(0,1){1}}
\put(1.5,0){\line(0,1){1}}
\put(-1,2.2){\makebox(0,0)[t]{$D_3$:}}
\put(0,2.2){
\unitlength.2cm
  \put(1.5,0){\bezier{30}(0,0)(.5,.5)(1,1)}
  \put(1.5,0){\bezier{30}(0,0)(-.5,.5)(-1,1)}
  \put(.5,1){\line(0,1){3}}
  \put(.5,4.4){\makebox(0,0)[b]{\scriptsize 1}}
}
\put(1,2.2){
\unitlength.2cm
  \put(1.5,0){\bezier{30}(0,0)(.5,.5)(1,1)}
  \put(1.5,0){\bezier{30}(0,0)(-.5,.5)(-1,1)}
  \put(2.5,1){\line(0,1){3}}
  \put(2.5,4.4){\makebox(0,0)[b]{\scriptsize 3}}
}
\put(2,2.2){
\unitlength.2cm
  \put(1.5,0){\bezier{30}(0,0)(.5,.5)(1,1)}
  \put(1.5,0){\bezier{30}(0,0)(-.5,.5)(-1,1)}
  \put(.5,1){\line(0,1){3}}
  \put(.5,4.4){\makebox(0,0)[b]{\scriptsize 2}}
}
}
\put(17.5,0){
\put(1.5,0){\line(-1,1){1}}
\put(1.5,0){\line(1,1){1}}
\put(.5,1){\line(0,1){1}}
\put(2.5,1){\line(0,1){1}}
\put(-1.5,2.2){\makebox(0,0)[t]{$D_4$:}}
\put(-.44,2.2){
\unitlength.4cm
  \put(1.5,0){\line(1,1){1}}
  \put(1.5,0){\line(-1,1){1}}
  \put(.5,1){\line(1,1){1}}
  \put(.5,1){\line(-1,1){1}}
  \put(-.5,2.2){\makebox(0,0)[b]{\scriptsize 1}}
  \put(1.5,2.2){\makebox(0,0)[b]{\scriptsize 2}}
}
\put(1.57,2.2){
\unitlength.4cm
  \put(1.5,0){\line(0,1){2}}
  \put(1.5,2.2){\makebox(0,0)[b]{\scriptsize 3}}
}
}
\end{picture}
\end{center}
\caption{\label{fig:5}%
Examples of fiber diagrams.}
\end{figure}
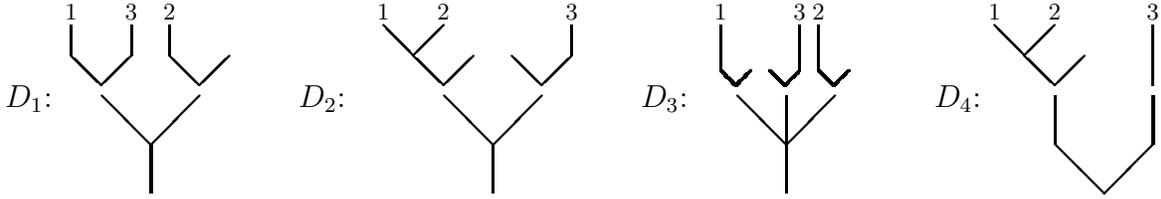

Each such a $\sigma : T \to S$ is characterized by its {\em fiber
diagram\/}, obtained by drawing fibers $F_j$ over the corresponding tips of
$S$. Some examples of fiber diagrams can be found in
Figure~\ref{fig:5}.

Diagram $D_1$ is the fiber diagram
of face $\sigma$ from Figure~\ref{fig:4}. Diagram $D_2$ is the
fiber diagram of the same trees as in Figure~\ref{fig:4}, but with
$\sigma$ determined by  $\sigma_2(1) = \sigma_2(2) = 1$ and $\sigma_2(3) = 2$.
Diagram $D_3$ is the diagram of the map $\sigma : \lp \to \trident$
given by $\sigma_2(1) = 1$,  $\sigma_2(2) = 3$ and $\sigma_2(3) =
2$. Diagram $D_4$ describes the map $\sigma : \lp \to \dvojdum$ given by
$\sigma_2(1) = \sigma_2(2) = 1$ and $\sigma_2(3) = 2$.

We sometimes decorate the tips of fibers by numbers that
indicate to which tip of $S$ they are mapped, see again Figure~\ref{fig:5}.
Other examples of fiber diagrams can be found in 
Figures~XV and~XVI of~\cite{batanin}.

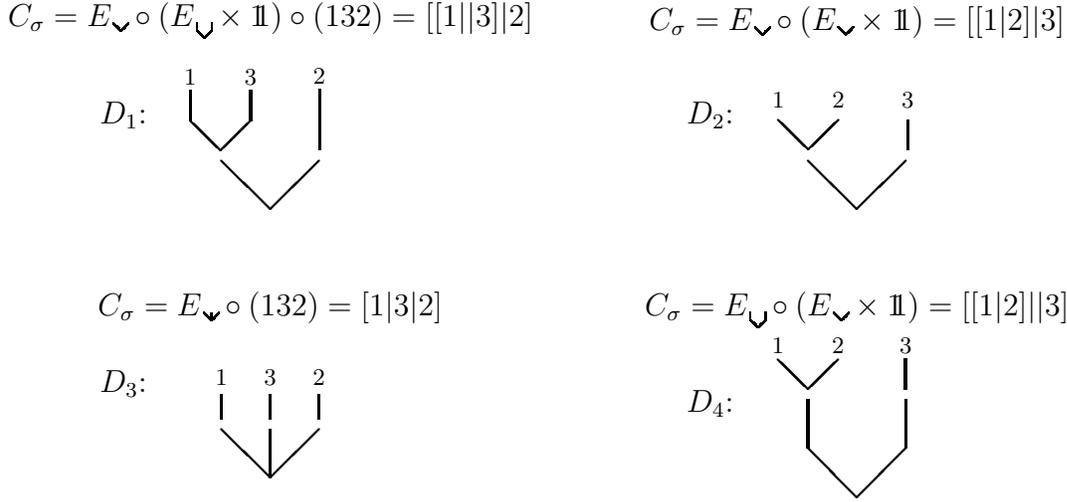
\begin{figure}[ht]
\begin{center}
\unitlength.65cm
\begin{picture}(10.00,10)(3,.5)
\thicklines
\put(0,6){
\put(1.5,3.5){\makebox(0,0)[b]%
{$C_\sigma = E_\strecha \hskip-5pt \circ (E_\dvojdum  
\hskip-5pt  \times \id)\circ(132) = [[1||3]|2]$}
}

\put(1.5,0){\line(-1,1){1}}
\put(1.5,0){\line(1,1){1}}
\put(-1.5,2.2){\makebox(0,0)[t]{$D_1$:}}
\put(-.44,1.2){
\unitlength.4cm
  \put(1.5,0){\line(1,1){1}}
  \put(1.5,0){\line(-1,1){1}}
  \put(0.5,1){\line(0,1){1}}
  \put(2.5,1){\line(0,1){1}}
  \put(.5,2.2){\makebox(0,0)[b]{\scriptsize 1}}
  \put(2.5,2.2){\makebox(0,0)[b]{\scriptsize 3}}
}
\put(2.2,.6){
\unitlength.4cm
  \put(.5,1){\line(0,1){2}}
  \put(.5,3.2){\makebox(0,0)[b]{\scriptsize 2}}
}
}
\put(12,6){
\put(1.5,3.5)%
   {\makebox(0,0)[b]%
           {$C_\sigma = E_\strecha  \hskip-5pt \circ (E_\strecha  
            \hskip-5pt \times \id) = [[1|2]|3]$}}
\put(1.5,0){\line(-1,1){1}}
\put(1.5,0){\line(1,1){1}}
\put(-1.5,2.2){\makebox(0,0)[t]{$D_2$:}}
\put(.2,0.6){
\unitlength.4cm
  \put(.5,1){\line(1,1){1}}
  \put(.5,1){\line(-1,1){1}}
  \put(-.5,2.4){\makebox(0,0)[b]{\scriptsize 1}}
  \put(1.5,2.4){\makebox(0,0)[b]{\scriptsize 2}}
}
\put(1,.6){
\unitlength.4cm
  \put(2.5,1){\line(0,1){1}}
  \put(2.5,2.4){\makebox(0,0)[b]{\scriptsize 3}}
}
}
\put(0,.5){
\put(1.5,3.8){\makebox(0,0)[t]{$C_\sigma = E_\assoc  
          \hskip-5pt \circ (132) = [1|3|2]$}}
\put(1.5,0){\line(-1,1){1}}
\put(1.5,0){\line(1,1){1}}
\put(1.5,0){\line(0,1){1}}
\put(-1.5,2.2){\makebox(0,0)[t]{$D_3$:}}
\put(.5,1.2){
\unitlength.65cm
  \put(0,0){\line(0,1){.5}}
  \put(0,.7){\makebox(0,0)[b]{\scriptsize 1}}
  \put(1,0){\line(0,1){.5}}
  \put(1,.7){\makebox(0,0)[b]{\scriptsize 3}}
  \put(2,0){\line(0,1){.5}}
  \put(2,.7){\makebox(0,0)[b]{\scriptsize 2}}
}
}
\put(12,.1){
\put(1.5,4.2){\makebox(0,0)[t]%
         {$C_\sigma = E_\dvojdum \hskip-5pt 
          \circ (E_\strecha \hskip-5pt \times \id) = [[1|2]||3]$}}
\put(1.5,0){\line(-1,1){1}}
\put(1.5,0){\line(1,1){1}}
\put(.5,1){\line(0,1){1}}
\put(2.5,1){\line(0,1){1}}
\put(-1.5,2.2){\makebox(0,0)[t]{$D_4$:}}
\put(.5,1.6){
\unitlength.4cm
  \put(0,1){\line(1,1){1}}
  \put(0,1){\line(-1,1){1}}
  \put(-1,2.2){\makebox(0,0)[b]{\scriptsize 1}}
  \put(1,2.2){\makebox(0,0)[b]{\scriptsize 2}}
}
\put(1.57,1.6){
\unitlength.4cm
  \put(1.5,1){\line(0,1){1}}
  \put(1.5,2.2){\makebox(0,0)[b]{\scriptsize 3}}
}
}
\end{picture}
\end{center}
\caption{\label{fig:6}%
Reduced fiber diagrams and elements $C_\sigma
\in \Free(\Tree) \cong \uFree(\uTree) \times \Sigma$ they determine. 
The symbol $\id$ denotes the identity
and $(132) \in \Sigma_3$ the permutation $(1,2,3) \mapsto
(1,3,2)$. For $T \in \uTree$, $E_T$ is the corresponding
generator of $\uFree(\uTree)$.}
\end{figure}

We are ready to describe the element $C_\sigma \in \Free(\Tree)$
indexing the cell $\mu[\sigma]:=\mu[C_\sigma]$ of
$\sfF$ \label{dnes_oslici}corresponding
to the face $\sigma$. We take the fiber diagram of $\sigma$ and
replace all trees of this diagram by their maximal reduced
subtrees. We obtain a tree-shaped diagram of reduced trees in $\Tree$
which are, by definition, the generators of the free operad set operad
$\Free(\Tree)$. The terminal tree represents the identity $\id \in
\Free(\Tree)(1)$. We then interpret this {\em reduced fiber diagram\/}
as the indicated composition of elements in the free operad
$\Free(\Tree)$ using the direct limit of isomorphisms~(\ref{Jarca!})
\[
\Free(\Tree) \cong \uFree(\uTree)\times \Sigma.
\] 
Let us denote this composition by $C_\sigma$.  We
believe that the construction of $C_\sigma$ is clear from
Figure~\ref{fig:6} which relies on
Notation~\ref{Jarka_zavolala_sama_od_sebe!}.

It will be convenient to extend the barcode notation of
Definition~\ref{S_Jarkou_do_Brna.}  to elements of the free operad
$\Free(\Tree)$. For example, the {\em extended barcode\/}
$[[1||3]|2]$ of the element $E_\strecha \hskip-5pt \circ (E_\dvojdum
\hskip-5pt \times \id)\circ(132)$ is obtained by inserting the 
barcode $[1||2]$ for the tree $E_\dvojdum$ into the first position in the
barcode $[1|2]$ for $E_\strecha$ and permuting the labels according
to the permutation $(132)$.
See Figure~\ref{fig:6} for more
examples of the extended barcodes.

The last step is counting $\deg(C_\sigma)$ by adding up the degrees of
generators that constitute~$C_\sigma$. For example, in
Figure~\ref{fig:6} all $C_\sigma$'s are of degree $1$ except the one
corresponding to $D_2$ which is of degree $0$. Let $\iota :
\Free(\Tree) \hookrightarrow \Free(G)$ be the monomorphism induced by
the inclusion $\Tree \hookrightarrow G = \Span(\Tree)$, $\bfT \mapsto
g_\bfT$, of graded $\Sigma$-sets. For $T \in \uTree$ and $g_T$ the
corresponding generator of $G$, put
\begin{equation}
\label{hyrim_barvami}
\pa(g_T) := \textstyle\sum_\sigma \pm \iota(C_\sigma),
\end{equation}
with the sum taken over faces $\sigma$ of $T$ such that
$\dim(C_\sigma) = \dim(T)-1$. The signs are determined by the orientation of
the cells. While it is possible to determine the signs for
each particular $g_T$, we do not know a reasonable general formula.

\begin{proposition}
\label{V_9:30_sraz_s_Jarkou_u_Pakousu.}
Formula~(\ref{hyrim_barvami}) extends to a differential on
$\Free(\Greg) \subset \Free(G)$ having the form $\pa = \pal + \papert$,
where $\pal$ is as in Definition~\ref{jarka}.
\end{proposition}

\begin{proof}
The first part of the proposition follows from the fact
that~(\ref{hyrim_barvami}) calculates the cellular differential of the
regular cell complex $\sfFreg$. Let us prove that the linear part 
of~(\ref{hyrim_barvami}) coincides with $\pal$.

It is obvious that, for $T \in \uTreereg$, $\pal(g_T)$ is given by the
sum~(\ref{hyrim_barvami}) restricted to the faces $\sigma : T \to S$
with trivial reduced fibers. Equivalently, we restrict to $\sigma$'s
that induce isomorphisms $\sigma_h : [k_h] \cong [s_h]$ of the
tips. Batanin calls, in~\cite{batanin:conf}, such maps $\sigma$ {\em
quasibijections\/}. The reduced fiber diagram of a quasiisomorphism
$\sigma$ is simply $S$ with the tips labeled by the permutation
$\sigma_h$, see $D_3$ in Figure~\ref{fig:6} for an example.  If we
denote this labeled tree by $\bfS_\sigma := (S,\sigma_h)$, then
$C_\sigma = \bfS_\sigma \in \Tree \subset \Free(\Tree)$.  Condition
$\dim(S) = \dim(T)-1$ means that the tree $S$ has one edge less than
$T$. Each such $S$ is obtained from $T$ by the following procedure.

Assume that $T$ is as in~(\ref{eq:1}) and choose 
$1 \leq m < h$ such that there exists $u \in [k_{m}]$ satisfying
$\rho_{m-1}(u) = \rho_{m-1}(u+1)$. Let $\Rada b1s$
(resp.~$\Rada b{s+1}{s+t}$) be the branches of $T$ over $u$
(resp.~$u+1$). 
By a {\em branch\/} over $u$ we mean a subtree 
determined by a vertex $\tilde u \in [k_{m+1}]$ satisfying
$\rho_m(\tilde u) = u$. The corresponding branch is the maximal subtree
of $T$ of height $h-m$ whose trunk is the edge connecting $\tilde u$
and $u$. Branches over $u+1$ are defined analogously. The situation is
shown in Figure~\ref{Vezu_Jarce_ruzicku.}.

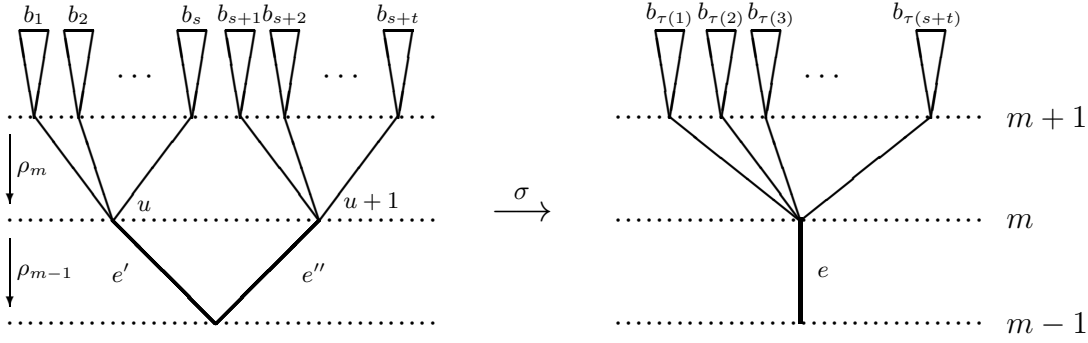
\begin{figure}[ht]
\begin{center}
\thicklines
\unitlength=1.300000pt
\begin{picture}(150.00,90.00)(0.00,0.00)
\put(-100,0){
\put(123.00,90.00){\makebox(0.00,0.00){\scriptsize $b_{s+t}$}}
\put(90.00,90.00){\makebox(0.00,0.00){\scriptsize $b_{s+2}$}}
\put(77.00,90.00){\makebox(0.00,0.00){\scriptsize $b_{s+1}$}}
\put(63.00,90.00){\makebox(0.00,0.00){\scriptsize $b_s$}}
\put(30.00,90.00){\makebox(0.00,0.00){\scriptsize $b_2$}}
\put(17.00,90.00){\makebox(0.00,0.00){\scriptsize $b_1$}}
\put(107.00,72.00){\makebox(0.00,0.00){$\cdots$}}
\put(47.00,72.00){\makebox(0.00,0.00){$\cdots$}}
\put(123.00,60.00){\makebox(0.00,0.00)[b]{\tria}}
\put(90.00,60.00){\makebox(0.00,0.00)[b]{\tria}}
\put(77.00,60.00){\makebox(0.00,0.00)[b]{\tria}}
\put(63.00,60.00){\makebox(0.00,0.00)[b]{\tria}}
\put(30.00,60.00){\makebox(0.00,0.00)[b]{\tria}}
\put(17.00,60.00){\makebox(0.00,0.00)[b]{\tria}}
\multiput(10,0)(0,30){3}{
\multiput(0,0)(3,0){42}{\makebox(0,0){$\cdot$}}
}
\put(150.00,30.00){\makebox(0.00,0.00)[lb]{{\large$\stackrel\sigma\longrightarrow$}}}
\put(107.00,33.00){\makebox(0.00,0.00)[lb]{\scriptsize $u+1$}}
\put(47.00,33.00){\makebox(0.00,0.00)[lb]{\scriptsize $u$}}
\put(100.00,30.00){\line(3,4){23.00}}
\put(100.00,30.00){\line(-1,3){10.00}}
\put(100.00,30.00){\line(-3,4){23.00}}
\put(40.00,30.00){\line(3,4){23.00}}
\put(40.00,30.00){\line(-1,3){10.00}}
{\thinlines
\put(10.00,25.00){\vector(0,-1){20.00}}
\put(10.00,55.00){\vector(0,-1){20.00}}
\put(12.00,15.00){\makebox(0.00,0.00)[l]{\scriptsize $\rho_{m-1}$}}
\put(12.00,45.00){\makebox(0.00,0.00)[l]{\scriptsize $\rho_{m}$}}
}
\put(45.00,15.00){\makebox(0.00,0.00)[r]{\scriptsize $e'$}}
\put(95.00,15.00){\makebox(0.00,0.00)[l]{\scriptsize $e''$}}
\put(40.00,30.00){\line(-3,4){23.00}}
\multiput(66.70,0)(.02,0){20}{
{\line(-1,1){30.00}}}
\multiput(66.70,0)(.02,0){20}{
{\line(1,1){30.00}}}
}
\put(100,0){
\unitlength=1.300000pt
\put(77.00,90.00){\makebox(0.00,0.00){\scriptsize $b_{\tau(s+t)}$}}
\put(30.00,90.00){\makebox(0.00,0.00){\ \scriptsize $b_{\tau(3)}$}}
\put(17.00,90.00){\makebox(0.00,0.00){\scriptsize $b_{\tau(2)}$}}
\put(2.00,90.00){\makebox(0.00,0.00){\scriptsize $b_{\tau(1)}$}}
\put(47.00,72.00){\makebox(0.00,0.00){$\cdots$}}
\put(78.00,60.00){\makebox(0.00,0.00)[b]{\tria}}
\put(30.00,60.00){\makebox(0.00,0.00)[b]{\tria}}
\put(17.00,60.00){\makebox(0.00,0.00)[b]{\tria}}
\put(2.00,60.00){\makebox(0.00,0.00)[b]{\tria}}
\multiput(-13,0)(0,30){3}{
\multiput(0,0)(3,0){36}{\makebox(0,0){$\cdot$}}
}
\put(100.00,30.00){\makebox(0.00,0.00)[l]{$m$}}
\put(100.00,60.00){\makebox(0.00,0.00)[l]{$m+1$}}
\put(100.00,0.00){\makebox(0.00,0.00)[l]{$m-1$}}
\put(40.00,30.00){\line(5,4){38.00}}
\put(40.00,30.00){\line(-1,3){10.00}}
\put(40.00,30.00){\line(-3,4){23.00}}
\put(40.00,30.00){\line(-5,4){38.00}}
\put(45.00,15.00){\makebox(0.00,0.00)[l]{\scriptsize $e$}}
\multiput(39.60,0)(.04,0){20}{\line(0,1){31.00}}
}
\end{picture}
\end{center}
\caption{\label{Vezu_Jarce_ruzicku.}
The relevant parts of the trees $T$ (left) and $S$ (right). 
The `fat' edge $e$ is obtained
by identifying $e'$ and $e''$.}
\end{figure}

Choose finally an $(s,t)$-unshuffle $\tau \in \Sigma_{s,t}$ and denote
by $\bfS$ the labeled tree obtained from $T$ by identifying the edge
$e'$ starting from $u$ with the edge $e''$ starting at $u+1$, and
permuting the branches $\Rada b1s,\Rada b{s+1}{s+t}$ according the
shuffle $\tau$, see again Figure~\ref{Vezu_Jarce_ruzicku.}.  Let
$\sigma : T \to \bfS$ be the projection.  It is clear that all
codimension-one faces $\sigma$ of $T$ are of this form and that
$g_\bfS \in G$ corresponds, under the isomorphism $\omega : G \to
\sfB_{\rada 11}^\infty(V)$ defined on page~\pageref{T}, to a component
of the top-level differential in $\sfB^{m+1}(V) \subset
\sfB^\infty(V)$ applied to $\omega(g_T)$.
\end{proof}

\noindent
\subsection{Proof of Theorem~A}
\label{Jarka_mne_mozna_miluje}  
Recall that the {\em cobar construction\/} on a coaugmented cooperad
$\calC$ is the dg-operad $\Omega(\calC)$ of the form $\Omega(\calC) =
(\Free(\pomio),\pa_\Omega)$, where $\overline{\calC}$ is the
co-augmentation co-ideal of $\calC$ and $\pomio$ the $\Sigma$-module
defined~by
\[
\pomio(n) := \sgn_n \otimes \uparrow^{n-2} \overline{\calC}(n),\ n
\geq 2,
\]
the product of the signum representation and the suspension
of $\overline{\calC}$ iterated 
$(n-2)$-times. The differential $\pa_\Omega$ is induced in the standard
manner from the structure operations of the cooperad $\calC$, 
see~\cite[Definition~II.3.9]{markl-shnider-stasheff:book}.

Let $\Lie = \{\Lie(n)\}_{n \geq 1}$ be the operad for Lie
algebras. It is well-known that each component $\Lie(n)$ of this
operad is a finite-dimensional free abelian group.  Denote by $\Lie' =
\{\Lie(n)'\}_{n \geq 1}$ the component-wise linear dual of $\Lie$ with
the induced cooperad structure. It follows from Theorem~6.7, Fact 6.2
and Proposition 5.2.12 of~\cite{fresse:CM04} that the natural morphism
\begin{equation}
\label{Za_chvili_volam_Jarce.}
\alpha :\Omega(\Lie') \to \Com
\end{equation}
of dg-operads is a homology isomorphism over $\ZZZ$.  This can also be
expressed by saying that the operads $\Com$ and $\Lie$ are Koszul dual
to each other, and Koszul over $\ZZZ$.  The last step is based on
Lemma~\ref{Zitra_do_VO} below, formulated in generality that allows to
use it in our proof of Theorem~B as well. For related results, see
Propositions~5.2.13 and~3.2.6 of~\cite{fresse:CM04}. Before we
formulate the lemma, we explain what precisely we mean by a
perturbation.

\begin{definition}
\label{zase_strasny_smutek}
Let $(M,\pal) = \{(M(n),\pal)\}_{n\geq 2}$ be a dg-$\Sigma$-module. A
{\em perturbation\/} of $\pal$ is a degree $-1$ map $\papert : M \to
\Free^{\geq 2}(M)$ from $M$ to the decomposables in the free operad
$\Free(M)$ such that $\pa := \pal + \papert : M \to \Free(M)$ extends
to a differential of \ $\Free(M)$.
\end{definition}

Let us state the assumptions of the lemma.  As everywhere in the
paper, we will use the same symbols for $\Sigma$-module maps $M \to
\Free(M)$ and for their extensions to derivations of~$\Free(M)$.
Suppose we are given a dg-operad $(\Free(E),\vartheta)$, for some
$\Sigma$-module $E = \{E(n)\}_{n \geq 2}$, such that $E(n)$ is
concentrated in degree $n-2$ and $\vartheta(E)$ consist of
decomposable elements in $\Free(E)$.  Assume we also have, for each $n
\geq 2$, a $\Sigma_n$-projective resolution $\gamma_n : (M(n),\pal) \to
(E(n),0)$ such that $M(n)$ is trivial in degrees $< n-2$. We denote by
$\gamma : \Free(M) \to \Free(E)$ the map of free operads induced by
$\{\gamma_n\}_{n \geq 2}$.

Suppose there are, for each $n \geq 4$, $\Sigma_n$-modules
$A(n)$ and $B(n)$ such that $M_n(n) = A(n) \oplus B(n)$ ($A(n)$ and
$B(n)$ are then projective, too). Let us denote by $\Mbar =  
\{\Mbar(n)\}_{n \geq 2}$ the graded $\Sigma$-sub-module of $M$ defined by
\[
\Mbar_i(n) := \tricases{M_i(n)}{if $n < 4$ or $n \geq  4$ and $i<n$,}
                         {A(n)}{if $n \geq 4$ and $i=n$, and}
                         {0}{in the remaining cases.}
\]
Finally, suppose we are given a perturbation $\opapert$ of the
restricted differential $\opal := \pal|_{\Mbar}$ such that 
the restricted map $\gamma : (\Free(\Mbar),\opa) \to
(\Free(E),\vartheta)$, where $\opa := \opal + \opapert$, is 
a morphism of dg-operads.

\begin{lemma}
\label{Zitra_do_VO}
In the above situation, there exists a perturbation $\papert$ of the
differential $\pal$ on $\Free(M)$ that extends the given restricted
perturbation
$\opapert$ such that $\gamma : (\Free(M),\pa) \to (\Free(E),\vartheta)$,
where $\pa := \pal + \papert$, is a map of dg-operads inducing an
isomorphism of homology.
\end{lemma}

\begin{proof}
We start by stating three simple facts. 

(i) From degree reasons, for an arbitrary extension $\papert$ of the
    perturbation $\opapert$, the extended differential $\pa = \pal +
    \papert$ will always commute with the homomorphism $\gamma :
    (\Free(M),\pa) \to (\Free(E),\vartheta)$. We therefore do not need
    to check the compatibility with $\gamma$.

(ii) By an easy spectral sequence argument, for every extension $\pa$
     as in (i), the map $\gamma: (\Free(M),\pa) \to
     (\Free(E),\vartheta)$ is a homology isomorphism, i.e.~the last
     property of the lemma is automatic.

(3) Since the free operad functor is a tensor functor,
    $H_*(\Free(M),\pal) \cong \Free(H_*(M,\pal)) \cong \Free(E)$. Again
    from simple degree reasons, $H_i(\Free^{\geq 2}(M),\pal)(n) = 0$ for
    $i \geq n-2$.

By (i), the only equation for the extended differential $\pa = \pal +
\papert$ which
has to be checked is $\pa^2=0$, that is
\[
(\pal + \papert)(\pal + \papert) = \pal\pal + \pal\papert +
\papert\pal + \papert \papert = 0
\]
which, since $\pal\pal = 0$, reduces to
\[
\pal\papert +
\papert\pal + \papert \papert = 0.
\]
Decomposing $\papert = \sum_{k \geq 2} \pp k$ (locally finite
sum), in which $\pp k : M \to \Free^k(M)$ is the component that sends
the generators $M$ to the sub-$\Sigma$-module $\Free^k(M) \subset
\Free(M)$ 
spanned by trees with $k$
vertices, one can expand the above equation into the system
\begin{equation}
\label{5}
\pal \pp l + \pp l \pal + \textstyle\sum_{i+j = l+1} \pp i\pp j =0
\end{equation}
which has to be satisfied for each $l \geq 2$. The extension $\papert$
will be constructed by induction on $l$, degree $i$ and arity $n$
using standard obstruction theory.

Let us start by extending $\opapert$ to $B(4)$. Let $k \geq 2$ and
suppose we already constructed, for each $2 \leq l < k$, 
$\pp l : B(4) \to \Free^l(M)(4)$ satisfying~(\ref{5}). Let us denote
\[
\ob^k := -  \pp k \pal - \textstyle\sum_{i+j = k+1} \pp i\pp j : B(4)
\to  \Free^k_2(M)(4).
\] 
Observe that the term $\pp k \pal$ in the above display has already
been defined. It is simple to verify that $\pal \ob^k =
0$, that is, $\ob^k$ is in fact a map from $B(4)$ to
$2$-dimensional  $\pal$-cycles in
$\Free^{\geq 2}(M)(4)$. 
It follows from the $\pal$-acyclicity~(iii) and the projectivity
of $B(4)$ that there exists a map $\pp k :  B(4) \to \Free^k_3(M)(4)$
such that $\ob^k = \pal \pp k$ which is the same as
\[
\pal \pp k + \pp k \pal + \textstyle\sum_{i+j = k+1} \pp i\pp j =0.
\] 
Repeating this process we extend $\opapert$ onto $B(4)$. In exactly
the same way, we extend $\opapert$ onto $M(4)$. 
The construction of an extension of $\opapert$ onto $M(n)$ for $n >
4$ is the same.  
\end{proof}

The existence of the operad $J = (\Free(G),\pa)$, $\pa =\pal +
\papert$, as in Theorem~A now follows from Lemma~\ref{Zitra_do_VO} in
which we take $E$ the $\Sigma$-module whose $n$th component equals $
\sgn_n \otimes \uparrow^{n-2}\! \Lie(n)'$ and $\vartheta$ the cobar
differential, i.e.~$(\Free(E),\vartheta) := \Omega(\Lie')$.  As
$(M,\pal)$ we take the $\Sigma$-module $(G,\pal)$ and set
$A(n)$ to be trivial for each $n \geq 4$, so $M = \Mbar$. In place of $\opa =
\opal + \opapert$ we take the differential from
Proposition~\ref{JarunA}.  Since the $\Sigma$-module $M$ is, by
Theorem~\ref{Stale_chodim_s_Jarkou}, a component-wise $\Sigma$-free
resolution of the collection $E$ defined above, the assumptions of the
lemma are fulfilled. 

The operad $J$ resolves $\Com$ via the composition
\[
J = (\Free(G),\pa) \stackrel\gamma\to \Omega(\Lie') \stackrel\alpha\to \Com,
\]
of the map $\gamma$ of Lemma~\ref{Zitra_do_VO} and $\alpha$
in~(\ref{Za_chvili_volam_Jarce.}).  It remains to prove that $J$ is
cofibrant.  To this end we show that there exists a totally ordered set $\Lambda$
such that $G$ decomposes as $G = \bigoplus_{\lambda \in
\Lambda}G_\lambda$ and
\begin{equation}
\label{Opicak_Fuk}
\pa(G_\lambda) \subset \Free(G)_{< \lambda},
\mbox { for each $\lambda \in  \Lambda$,}
\end{equation}
where $\Free(E)_{< \lambda}$ denotes the suboperad of $\Free(E)$
generated by $\bigoplus_{\lambda' < \lambda} E_{\lambda'}$.  It
follows from~\cite[Lemma~20]{markl:ha}
that~(\ref{Opicak_Fuk}) guarantees the lifting property of $J$ with
respect to trivial fibrations, so $J$
is cofibrant in the standard model structure of the category of
operads~\cite[Example~3.3.3]{berger-moerdijk:02}.

Let $\Lambda := {\mathbb N} \times {\mathbb N}$ be the cartesian
product of two copies of natural numbers with the lexicographic order
and, for $(i,n) \in \Lambda$, $G_{(i,n)}:= G_i(n)$ (the subspace of
arity $n$ and degree $i$). While obviously $\pal ( G_i(n)) \subset
G_{i-1}(n)$, it follows from simple combinatorics that $\papert(G(n))$
consists of compositions of elements of arities $< n$.  This
establishes~(\ref{Opicak_Fuk}) for the sum $\pa = \pal +
\papert$ and finishes our proof of Theorem~A.

\noindent 
\subsection{Some formulas}
\label{jarUnka}
In this subsection we calculate the differential of some
low-dimensional generators of the operad $J = \Free(G)$.
Recall that, for a tree $T \in \uTree(n)$, we denoted by $\genname_T$ the
corresponding generator of $G(n)$. We will denote 
a \label{Dnes_mi_trhali_zub!}permutation
$\sigma : \{\rada 1n\} \to \{\rada 1n\}  \in
\Sigma_n$ by the $n$-tuple $(\Rada \sigma 1n)$,
where $\sigma_i := \sigma^{-1}(i)$ for $1\leq i \leq n$.

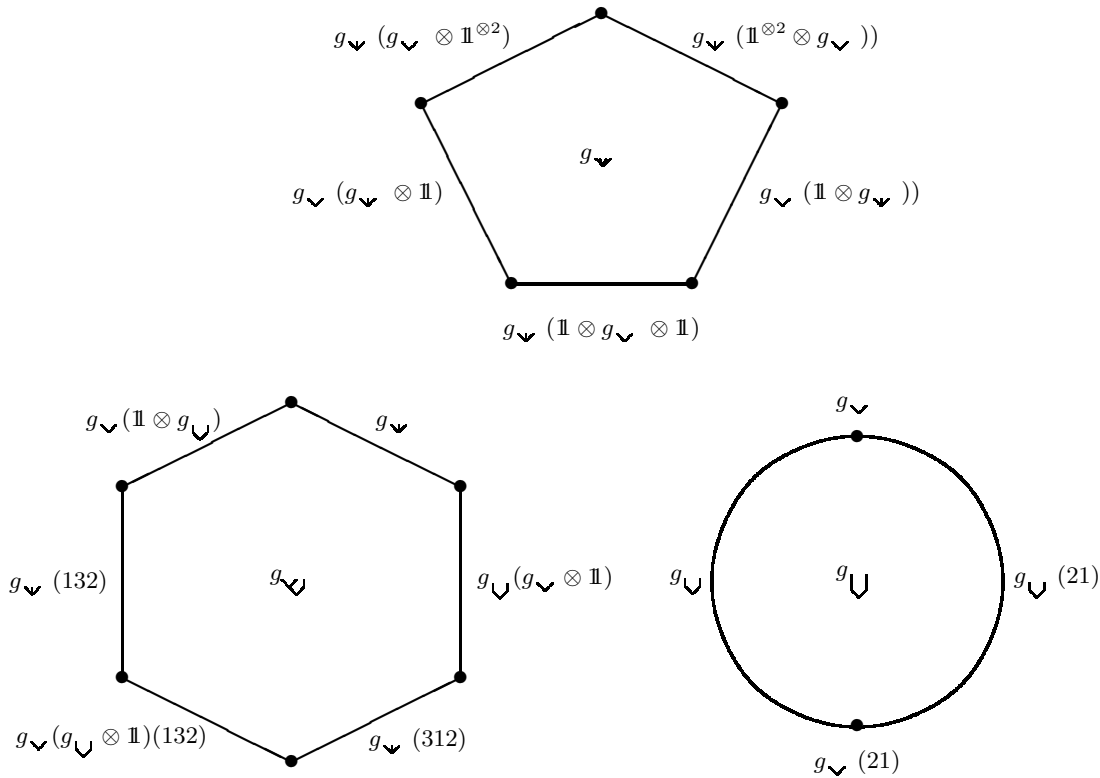
\begin{figure}[ht]
\begin{center}
{
\thicklines
\unitlength=.8pt
\begin{picture}(210.00,360.00)(0.00,0.00)
\put(-150,0){
\put(45.00,160.00){\makebox(0.00,0.00){%
       \scriptsize$\genname_\strecha \hskip -1mm (\id \ot \genname_\dvojdum \hskip -1mm  )$}}
\put(160.00,160.00){\makebox(0.00,0.00){\scriptsize$\genname_\assoc$}}
\put(230.00,85.00){\makebox(0.00,0.00){%
       \scriptsize$\genname_\dvojdum\hskip -1mm (\genname_\strecha \hskip -1mm \ot \id)$}}
\put(170.00,10.00){\makebox(0.00,0.00){\scriptsize$\genname_\assoc(312)$}}
\put(25.00,10.00){\makebox(0.00,0.00){%
    \scriptsize$\genname_\strecha\hskip -1mm (\genname_\dvojdum\hskip -1mm  \ot \id)(132)$}}
\put(0.00,85.00){\makebox(0.00,0.00){\scriptsize$\genname_\assoc (132)$}}
\put(110.00,85.00){\makebox(0.00,0.00){\scriptsize$\genname_\lp$}}
\put(110.00,170){\makebox(0.00,0.00){$\bullet$}}
\put(190.00,130.00){\makebox(0.00,0.00){$\bullet$}}
\put(190.00,40.00){\makebox(0.00,0.00){$\bullet$}}
\put(110.00,0.00){\makebox(0.00,0.00){$\bullet$}}
\put(30.00,40.00){\makebox(0.00,0.00){$\bullet$}}
\put(30.00,130.00){\makebox(0.00,0.00){$\bullet$}}
\put(110.00,170.00){\line(-2,-1){80.00}}
\put(190.00,130.00){\line(-2,1){80.00}}
\put(190.00,40.00){\line(0,1){90.00}}
\put(110.00,0.00){\line(2,1){80.00}}
\put(30.00,40.00){\line(2,-1){80.00}}
\put(30.00,130.00){\line(0,-1){90.00}}
}
\unitlength=.8cm
\put(1,0){
\put(7,3){\krouzek{.968}}
\put(7,3){\makebox(0,0){\scriptsize $\genname_\trojdum$}}
\put(7,5.4){\makebox(0,0){$\bullet$}}
\put(7,.6){\makebox(0,0){$\bullet$}}
\put(7,5.8){\makebox(0,0)[b]{\scriptsize $\genname_\strecha$}}
\put(7,.2){\makebox(0,0)[t]{\scriptsize $\genname_\strecha(21)$}}
\put(4.6,3){\makebox(0,0)[r]{\scriptsize $\genname_\dvojdum$}}
\put(9.6,3){\makebox(0,0)[l]{\scriptsize $\genname_\dvojdum(21)$}}
}
\unitlength 1.2cm
\put(2.5,5.3){
%
\put(-1,0){\line(1,0){2}}
\put(1,0){\line(1,2){1}}
\put(0,3){\line(2,-1){2}}
\put(-2,2){\line(2,1){2}}
\put(-2,2){\line(1,-2){1}}
%
\put(1,0){\makebox(0,0){$\bullet$}}
\put(2,2){\makebox(0,0){$\bullet$}}
\put(0,3){\makebox(0,0){$\bullet$}}
\put(-2,2){\makebox(0,0){$\bullet$}}
\put(-1,0){\makebox(0,0){$\bullet$}}
\put(0,-.5){\makebox(0,0){{%
         \scriptsize $\genname_\assoc(\id \ot \genname_\strecha \ot \id)$}}}
\put(-1.75,1){\makebox(0,0)[r]{{%
         \scriptsize $\genname_\strecha(\genname_\assoc \ot \id)$}}}
\put(1.75,1){\makebox(0,0)[l]{{%
         \scriptsize $\genname_\strecha(\id \otimes \genname_\assoc))$}}}
\put(-1,2.75){\makebox(0,0)[r]{{%
         \scriptsize $\genname_\assoc(\genname_\strecha \ot \otexp \id2)$}}}
\put(1,2.75){\makebox(0,0)[l]{{%
         \scriptsize $\genname_\assoc(\otexp \id2 \ot \genname_\strecha))$}}}
\put(0,1.4){\makebox(0,0){{%
         \scriptsize $\genname_\pent$}}}
}
\end{picture}}
\end{center}
\caption{\label{fig:7}%
Pentagon, left hexagon and disk.}
\end{figure}

The degree $0$ generator $\genname_\strecha \hskip -.2em$ 
is mapped to the commutative associative multiplication in
$\Com$. Of course, $\pa(\genname_\strecha \hskip -.2em) = 0$.
The degree one generator $\genname_\assoc$ is the `associator' and
\[
\pa(\genname_\assoc \hskip -.2em) 
= \genname_\strecha \hskip -.2em \circ (\genname_\strecha \hskip -.2em\ot \id 
- \id \ot \genname_\strecha\hskip -.2em).
\]
The second degree one generator $\genname_\dvojdum \hskip -.2em$ 
represents the homotopy
for the commutativity~of~$\genname_\strecha \hskip -.2em$:
\[
\pa(\genname_\dvojdum \hskip -.2em) = \genname_\strecha \hskip -.2em(\id - (21)).
\]
The degree two generator $\genname_\pent$ is the Stasheff/Mac Lane pentagon and
we all know the formula
\[
\pa(\genname_\pent \hskip -.2em) =
\genname_\assoc \hskip -1.5mm (\genname_\strecha \hskip -1.5mm \ot \otexp {\id}2) -
\genname_\assoc\hskip -1.5mm (\id \ot \genname_\strecha \hskip -1.5mm \ot \id)+
\genname_\assoc \hskip -1.5mm (\otexp {\id}2 \ot \genname_\strecha \hskip -1.5mm)
- \genname_\strecha\hskip -1.5mm (\genname_\assoc \hskip -1.5mm \ot \id) - 
\genname_\strecha \hskip -1.5mm (\id \ot \genname_\assoc \hskip -1.5mm)
\]
from kindergarten, see Figure~\ref{fig:7}. The
degree two generator $\genname_\lp$ is the left 
hexagon whose differential is given by
\[
\pa(\genname_\lp \hskip -.2em) = \genname_\assoc \hskip -1.5mm \circ (\id - (132) + (312)) 
           - \genname_\dvojdum \hskip -1.5mm\circ (\genname_\strecha \hskip -1.5mm \ot \id)
           + \genname_\strecha\hskip -1.5mm \circ (\id \ot \genname_\dvojdum \hskip -1mm) 
       + \genname_\strecha \hskip -1.5mm\circ (\genname_\dvojdum \hskip -1.5mm \ot \id)(132),
\]
see Figure~\ref{fig:7}.
The formula for the right hexagon $\genname_\rp$ is 
a similar:
\[
\pa(\genname_\rp) = \genname_\assoc \hskip -1.5mm \circ (\id - (213) + (231)) 
       + \genname_\dvojdum \hskip -1.5mm\circ (\id \ot \genname_\strecha \hskip -1.5mm)
   - \genname_\strecha\hskip -1.5mm \circ (\id \ot \genname_\dvojdum \hskip -1mm)(213)
      -\genname_\strecha \hskip -1.5mm\circ (\genname_\dvojdum \hskip -1.5mm \ot \id).
\]
The last degree two generator $\genname_\trojdum \hskip -.2em$ is the
homotopy for the anticommutativity of $\genname_\dvojdum  \hskip -.2em$:
\[
\pa(\genname_\trojdum  \hskip -.2em) = \genname_\dvojdum \hskip -.2em(\id +(21)),
\]
see again Figure~\ref{fig:7}.
The differential of the degree three generator $\genname_\dvojtroj  \hskip -.2em$ is
\[
\pa(\genname_\dvojtroj) = \genname_\lp \circ (\id - (213)) - \genname_\rp \circ (\id - (132))
                 - \genname_\dvojdum  \hskip -.2em \circ (\genname_\dvojdum  \hskip
		 -.2em\ot \id - \id \ot \genname_\dvojdum  \hskip -.2em).
\]
Let us give also a formula for $\partial(\genname_\levaslozitost \hskip -1.5mm)$:
\[
\pa(\genname_\levaslozitost \hskip -2mm) = \genname_\lp - \genname_\rp \circ (312) 
       + \genname_\trojdum \hskip -1mm \circ (\genname_\strecha\hskip -1mm \ot \id)
       - \genname_\strecha  \hskip -1mm\circ (\id \ot \genname_\trojdum \hskip -1mm)
       - \genname_\strecha  \hskip -1mm \circ (\genname_\trojdum \hskip -1mm \ot \id)(132).
\]

\noindent 
\subsection{$E_\infty$-algebras} 
As explained in the Introduction, algebras over the differential
graded operad $J = (J,\pa)$ are particular realizations of {\em
$E_\infty$-algebras\/}. A structure of this $E_\infty$-algebra on a
dg-abelian group $V = (V,d)$ is given
by multilinear maps $\mu_T : \otexp Vn \to V$ indexed by
reduced trees. The degree of $\mu_T$ equals $\dim(T)$ and the arity
$n$ equals the number of the tips of~$T$. The axioms could be read off from the
formulas for differential $\pa$  given in
Subsection~\ref{jarUnka}. One gets
\begin{subequations}
\allowdisplaybreaks
\begin{eqnarray}
\label{dnes_vecer_v_8_sraz_s_Jarunkou}
\delta \mu_\strecha (a,b) 
&=&
0,
\\
\label{hom-ass}
\delta \mu_\assoc(a,b,c)
&=&
\mu_\strecha(\mu_\strecha(a,b),c) - \mu_\strecha(a,\mu_\strecha(b,c)),
\\
\label{cup-jedna}
\delta \mu_\dvojdum(a,b)
&=&
\mu_\strecha(a,b) - \mu_\strecha(b,a),
\\
\label{hexag}
\delta \mu_\lp(a,b,c)
&=&
\mu_\assoc(a,b,c) - \mu_\assoc(a,c,b) + \mu_\assoc(c,a,b) -
\mu_\dvojdum(\mu_\strecha(a,b),c) +
\\
\nonumber 
&&\hskip 4mm
(-1)^{\deg(a)} 
\mu_\strecha(a,\mu_\dvojdum(b,c))
+\mu_\strecha(\mu_\dvojdum(a,c),b),
\\
\label{hexagr}
\delta \mu_\rp(a,b,c)
&=&
\mu_\assoc(a,b,c) - \mu_\assoc(b,a,c) + \mu_\assoc(b,c,a) -
\mu_\strecha(\mu_\dvojdum(a,b),c), 
\\
\nonumber 
&&\hskip 4mm
(-1)^{\deg(a)}\mu_\dvojdum(a,\mu_\strecha(b,c)) -(-1)^{\deg(b)} 
\mu_\strecha(b,\mu_\dvojdum(a,c))
\\
\label{pentag}
\delta \mu_\pent(a,b,c,d)
&=&
\mu_\assoc(\mu_\strecha(a,b),c,d) -
\mu_\assoc(a,\mu_\strecha(b,c),d)+ \mu_\assoc(a,b,\mu_\strecha(c,d))
\\ 
\nonumber 
&&\hskip 4mm -
\mu_\strecha(\mu_\assoc(a,b,c),d) - (-1)^{\deg(a)} 
\mu_\strecha(a,\mu_\assoc(b,c,d)),
\\
\label{cup-dva}
\delta \mu_\trojdum(a,b)
&=&
\mu_\dvojdum(a,b) + \mu_\dvojdum(b,a),
\\
\nonumber 
\delta \mu_\dvojtroj(a,b,c) 
&=&
\nonumber 
\mu_\lp(a,b,c) - \mu_\lp(b,a,c) - \mu_\rp(a,b,c) + \mu_\rp(a,c,a)+
\\
\nonumber 
&&\hskip 4mm       - \mu_\dvojdum (\mu_\dvojdum(a,b),c) - 
                  (-1)^{\deg(a)}  \mu_\dvojdum (a,\mu_\dvojdum(b,c)),
\\
\nonumber 
\delta \mu_\levaslozitost(a,b,c) 
&=& 
\mu_\lp(a,b,c) - \mu_\rp(c,a,b) + \mu_\trojdum(\mu_\strecha(a,b),c) -
\\
\nonumber 
&&\hskip 4mm - \mu_\strecha (a,\mu_\trojdum(b,c))
- \mu_\strecha (\mu_\trojdum(a,c),b), \mbox { \&c.}
\end{eqnarray}
\end{subequations}

In the above formulas, $a$, $b$, $c$, $d$ are homogeneous elements
of $V$, and $\delta$ the induced differential in
the endomorphism complex of $V = (V,d)$. For example
\[
\delta \mu_\strecha (a,b) := d \mu_\strecha (a,b) -  \mu_\strecha
(da,b) - (-1)^{\deg(a)} \mu_\strecha (a,db).
\]

Some of the above axioms have obvious
interpretations. Axiom~(\ref{hom-ass}) says that $\mu_\assoc$ is a
(chain) homotopy for the multiplication $\mu_\strecha$,
axiom~(\ref{cup-jedna}) means that $\mu_\dvojdum$ is $\smile_1$ for
$\mu_\strecha$ and~(\ref{cup-dva}) means that $\mu_\trojdum$ is
$\smile_2$ for $\mu_\strecha$. Axioms~(\ref{hexag}) and~(\ref{hexagr})
are algebraic versions of the left and right
hexagons. Axiom~(\ref{pentag}) is an algebraic version of the
pentagon.

More generally, if $\mu_n := \mu_{\star_n}$ with $\star_n \in \uTree^1(n)$ 
the $n$-corolla with the
barcode $[1|\cdots|n]$, then $(V,d,\mu_2,\mu_3,\ldots)$ is
an $A_\infty$-algebra, with~(\ref{dnes_vecer_v_8_sraz_s_Jarunkou}),
(\ref{hom-ass}) and~(\ref{pentag}) Axiom~(1)
of~\cite{markl:JPAA92} for $n=2,3$ and~$4$.  This justifies calling
the operad $J$ an extension of Stasheff's operad.
Axioms~(\ref{dnes_vecer_v_8_sraz_s_Jarunkou})--(\ref{cup-dva}) were
already obtained in Example~4.8 of the 1996 paper~\cite{markl:JPAA96}.

\section{The Tamarkin cell mystery}
\label{Tamarkin}

The first example of a cell violating the regularity of the CW-complex
$\sfF$ was found by Dimitri Tamarkin.  It is a $6$-dimensional cell
$\Tam \subset \sfF(6)$ which actually lives in the subcomplex
$\sfF_2(6)$ of compactified configurations of six points in
$\RRR^2$. Surprisingly, there exist even a simpler, $4$-dimensional
`bad' cell $\Mar \subset \sfF(4)$ living in the subcomplex $\sfF_3(4)$
of compactified configurations of $4$ points in $\RRR^3$.  It will be
clear from Section~\ref{jaRcA} that $\Mar$ has the smallest possible
dimension.  In this section we analyze the above two cells and show
how to construct the differential of the corresponding generators of $G$.
Theorem~B proved at the end of this section shows that any such a
`partial' differential extends to a `global' one. 

\subsection{The Tamarkin cell}
\label{Jarca!!}
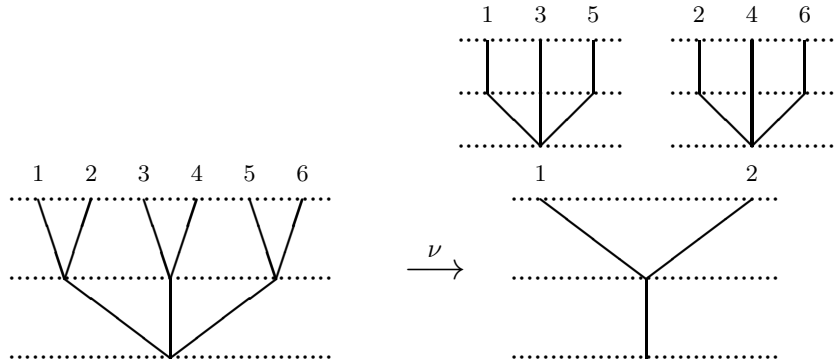
\begin{figure}[ht]
\begin{center}
{
\thicklines
\unitlength=1.000000pt
\begin{picture}(270.00,130.00)(0.00,0.00)
\put(130.00,30.00){\makebox(0,0)[b]{\large $\stackrel\nu\longrightarrow$}}
\put(-30,0){
\multiput(0,0)(0,30){3}{
\multiput(0,0)(3,0){41}{\makebox(0,0){$\cdot$}}
}
\put(60.00,30.00){\line(0,-1){30.00}}
\put(60.00,0.00){\line(4,3){40.00}}
\put(20.00,30.00){\line(4,-3){40.00}}
\put(100.00,30.00){\line(1,3){10.00}}
\put(90.00,60.00){\line(1,-3){10.00}}
\put(60.00,30.00){\line(1,3){10.00}}
\put(50.00,60.00){\line(1,-3){10.00}}
\put(20.00,30.00){\line(1,3){10.00}}
\put(10.00,60.00){\line(1,-3){10.00}}
\put(110.00,70.00){\makebox(0.00,0.00){\scriptsize $6$}}
\put(90.00,70.00){\makebox(0.00,0.00){\scriptsize $5$}}
\put(70.00,70.00){\makebox(0.00,0.00){\scriptsize $4$}}
\put(50.00,70.00){\makebox(0.00,0.00){\scriptsize $3$}}
\put(30.00,70.00){\makebox(0.00,0.00){\scriptsize $2$}}
\put(10.00,70.00){\makebox(0.00,0.00){\scriptsize $1$}}
}
\multiput(220,80)(0,20){3}{
\multiput(0,0)(3,0){21}{\makebox(0,0){$\cdot$}}
}
\multiput(140,80)(0,20){3}{
\multiput(0,0)(3,0){21}{\makebox(0,0){$\cdot$}}
}
\multiput(160,0)(0,30){3}{
\multiput(0,0)(3,0){34}{\makebox(0,0){$\cdot$}}
}
\put(190.00,130.00){\makebox(0.00,0.00){\scriptsize $5$}}
\put(170.00,130.00){\makebox(0.00,0.00){\scriptsize $3$}}
\put(150.00,130.00){\makebox(0.00,0.00){\scriptsize $1$}}
\put(230.00,130.00){\makebox(0.00,0.00){\scriptsize $2$}}
\put(250.00,130.00){\makebox(0.00,0.00){\scriptsize $4$}}
\put(270.00,130.00){\makebox(0.00,0.00){\scriptsize $6$}}
\put(170.00,70.00){\makebox(0.00,0.00){\scriptsize $1$}}
\put(250.00,70.00){\makebox(0.00,0.00){\scriptsize $2$}}
\put(270.00,100.00){\line(0,1){20.00}}
\put(250.00,80.00){\line(1,1){20.00}}
\put(250.00,120.00){\line(0,-1){40.00}}
\put(230.00,100.00){\line(0,1){20.00}}
\put(250.00,80.00){\line(-1,1){20.00}}
\put(190.00,100.00){\line(0,1){20.00}}
\put(170.00,80.00){\line(1,1){20.00}}
\put(170.00,120.00){\line(0,-1){40.00}}
\put(150.00,100.00){\line(0,1){20.00}}
\put(170.00,80.00){\line(-1,1){20.00}}
\put(210.00,30.00){\line(0,-1){30.00}}
\put(210.00,30.00){\line(4,3){40.00}}
\put(170.00,60.00){\line(4,-3){40.00}}
\end{picture}}
\end{center}
\caption{\label{Tam}
The Tamarkin tree $T$ (left), the tree $S$ (right
bottom), the map $\nu : T \to S$ and its fiber diagram.}
\end{figure}

\begin{figure}[ht]
\begin{center}
{
\unitlength=1.000000pt
\begin{picture}(190.00,160.00)(0.00,0.00)

\put(30.00,130.00){\makebox(0.00,0.00){\scriptsize $\mu[\nu]$}}
\put(1.00,40.00){\makebox(0.00,0.00){\large $\bullet$}}
\put(60.00,110.00){\makebox(0.00,0.00){\large $\bullet$}}
\put(110.00,70.00){\makebox(0.00,0.00)[lb]{\scriptsize $\Tam  =
    \mu[T]$}}
\put(20,60){\makebox(0.00,0.00)[lt]{\scriptsize$\pasing\Tam$}}
\put(120,10){\makebox(0.00,0.00)[t]{\scriptsize$\pareg\Tam$}}
\put(40,160){\makebox(0.00,0.00)[b]{\scriptsize$U'$}}
\put(40,0){\makebox(0.00,0.00)[t]{\scriptsize$U''$}}
\multiput(1.00,40.00)(1.6,1.89){38}{\makebox(0.00,0.00){$\cdot$}}
\put(70.00,60.00){\line(-1,5){10.00}}
\put(52.20,15.50){\line(2,5){18.00}}
\thicklines
\put(50.00,10.00){\line(2,5){4.00}}
\put(20.00,0.00){\line(3,1){30.00}}
\put(0.00,40.00){\line(1,-2){20.00}}
\put(0.00,100.00){\line(0,-1){60.00}}
\put(20.00,160.00){\line(-1,-3){20.00}}
\put(50.00,150.00){\line(-3,1){30.00}}
\put(60.00,110.00){\line(-1,4){10.00}}
\put(110.00,110.00){\line(-1,0){50.00}}
\put(190.00,90.00){\line(-4,1){80.00}}
\put(170.00,30.00){\line(1,3){20.00}}
\put(90.00,10.00){\line(4,1){80.00}}
\put(0.00,40.00){\line(3,-1){90.00}}
\end{picture}}
\end{center}
\caption{\label{!Jaruska}
The Tamarkin cell, the cell $\mu[\nu]$ and the intersection 
$\pasing \Tam := \pa \Tam \cap \mu[\nu]$.}
\end{figure}
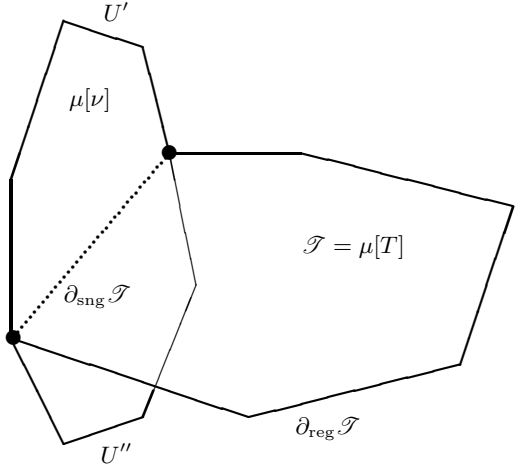

The Tamarkin cell \label{tam}$\Tam$ 
corresponds to the tree $T$ with
the barcode $[1|2||3|4||5|6]$ shown in Figure~\ref{Tam} (left). 
Consider the tree $S := [1|2]$ and the map $\nu : T \to S$ that
sends the tips of $T$ labeled 1,3,5 (reps.~labeled 2,4,6) to the tip
1 (resp.~tip 2) of $S$, see Figure~\ref{Tam}. We easily read off from
the fiber diagram of $\nu$ that $\mu[\nu] =
\mu\jarka{[1||3||5]\hskip .1em}{[2||4||6]}$. Simple degree count shows that
$\dim(\mu[\nu]) = 6$ i.e.~it is the same as the dimension of the
Tamarkin cell $\Tam$! The explanation is that the face $\mu[\nu]$
is not a subset of the boundary $\pa \Tam$ of $\Tam$ but that
$\pa \Tam$ intersect $\mu[\nu]$ in a $5$-dimensional
subspace which is not a union of $5$-dimensional cells. 
This violates Definition~\ref{JarunKa}(ii). An ``ideological'' picture of this
situation is shown in Figure~\ref{!Jaruska}.

\begin{figure}[ht]
\begin{center}
{
\unitlength=.600000pt
\begin{picture}(140.00,150.00)(90.00,20.00)
\put(-100,0){
\put(120.00,30.00){\makebox(0.00,0.00){\scriptsize $5$}}
\put(120.00,150.00){\makebox(0.00,0.00){\scriptsize $6$}}
\put(70.00,80.00){\makebox(0.00,0.00){\scriptsize $3$}}
\put(70.00,110.00){\makebox(0.00,0.00){\scriptsize $4$}}
\put(30.00,70.00){\makebox(0.00,0.00){\scriptsize $1$}}
\put(30.00,130.00){\makebox(0.00,0.00){\scriptsize $2$}}
\put(110.00,20.00){\makebox(0.00,0.00){$\bullet$}}
\put(110.00,140.00){\makebox(0.00,0.00){$\bullet$}}
\put(60.00,70.00){\makebox(0.00,0.00){$\bullet$}}
\put(60.00,100.00){\makebox(0.00,0.00){$\bullet$}}
\put(20.00,60.00){\makebox(0.00,0.00){$\bullet$}}
\put(20.00,120.00){\makebox(0.00,0.00){$\bullet$}}
\put(110.00,150.00){\line(0,-1){150.00}}
\put(60.00,150.00){\line(0,-1){150.00}}
\put(20.00,150.00){\line(0,-1){150.00}}
\thicklines
\put(-10.00,40.00){\vector(1,0){180.00}}
}
\put(115,20)
{
\thinlines
\put(60.00,30.00){\line(0,-1){30.00}}
\put(60.00,100.00){\line(0,-1){30.00}}
\put(60.00,160.00){\line(0,-1){20.00}}
\put(60.00,50.00){\circle{40.00}}
\put(60.00,120.00){\circle{40.00}}
\put(70.00,40.00){\makebox(0.00,0.00){$\bullet$}}
\put(62.00,60.00){\makebox(0.00,0.00){$\bullet$}}
\put(50.00,50.00){\makebox(0.00,0.00){$\bullet$}}
\put(70.00,130.00){\makebox(0.00,0.00){$\bullet$}}
\put(55.00,110.00){\makebox(0.00,0.00){$\bullet$}}
\put(48.00,120.00){\makebox(0.00,0.00){$\bullet$}}
\put(87.00,127.00){\makebox(0.00,0.00)[lb]{\scriptsize ${\mathbf x}_u$}}
\put(87.00,57.00){\makebox(0.00,0.00)[lb]{\scriptsize ${\mathbf x}_b$}}
\thicklines
\put(0.00,20.00){\vector(1,0){140.00}}
}
\put(300,20)
{
\thinlines
\multiput(50.00,160.00)(10,0){3}{
\multiput(0,0)(0,-5){32}{\makebox(0,0){$\cdot$}}}
\put(60.00,50.00){\circle{40.00}}
\put(60.00,120.00){\circle{40.00}}
\put(70.00,40.00){\makebox(0.00,0.00){$\bullet$}}
\put(60.00,60.00){\makebox(0.00,0.00){$\bullet$}}
\put(50.00,50.00){\makebox(0.00,0.00){$\bullet$}}
\put(70.00,130.00){\makebox(0.00,0.00){$\bullet$}}
\put(60.00,110.00){\makebox(0.00,0.00){$\bullet$}}
\put(50.00,120.00){\makebox(0.00,0.00){$\bullet$}}
\put(87.00,127.00){\makebox(0.00,0.00)[lb]{\scriptsize ${\mathbf x}_u$}}
\put(87.00,57.00){\makebox(0.00,0.00)[lb]{\scriptsize ${\mathbf x}_b$}}

\thicklines
\put(0.00,20.00){\vector(1,0){140.00}}
}
\end{picture}}
\end{center}
\caption{
\label{snad_mne_miluje}
A generic point of the Tamarkin cell  $\mu[1|2||3|4||5|6]$
(left), of the cell $\nu = 
\mu\protect\jarka{[1||3||5]\hskip .1em}{[2||4||6]}$ (middle) and
of their intersection (right).}
\end{figure}
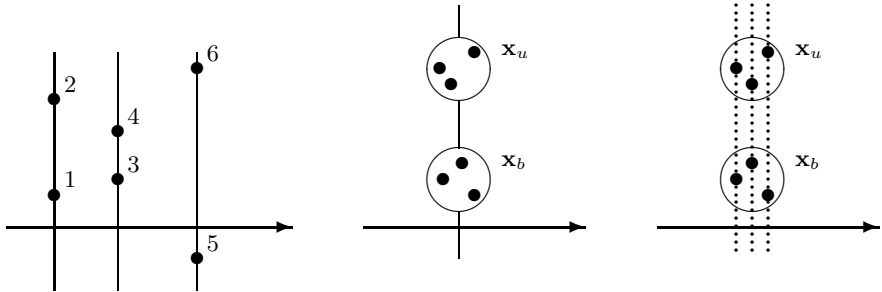

Let us analyze this phenomenon in detail. If we denote by  $c :
\RRR^2 \to \RRR$ the projection to the first coordinate, then 
each point ${\mathbf x} = (\Rada x16)$ in the interior of $\Tam$ satisfies
\begin{equation}
\label{!Jarka}
\frac{c(x_3)-c(x_1)}{c(x_5) - c(x_3)} = \frac{c(x_4)-c(x_2)}{c(x_6) - c(x_4)}.
\end{equation}
Observe that both sides are invariant under the affine group action.
The same condition is satisfied also by the points in the intersection
$\pa \Tam \cap \mu[\nu]$. More precisely, $\mu[\nu]$ consists of two
`microscopic' configurations ${\mathbf x}_u$ (resp.~${\mathbf x}_b$)
of points $(x_2,x_4,x_6)$ (resp.~$(x_1,x_2,x_3)$) in $\sfF(3)$
arranged at the vertical line (Figure~\ref{snad_mne_miluje}
middle). Since the points in $\pa \Tam \cap \mu[\nu]$ are the limits of
the points in the interior of $\Tam$, the configurations ${\mathbf
x}_u$ and ${\mathbf x}_b$ are still tied by~(\ref{!Jarka}), see
Figure~\ref{snad_mne_miluje} (right). Therefore the intersection $\pa
\Tam \cap \mu[\nu]$ is a codimension-one subspace of $\mu[\nu]$.
Loosely speaking, when the point ${\mathbf x} \in \Tam$ moves towards
the boundary, it still `remembers' that its coordinates were lined up
at tree vertical lines parametrized by a point in ${\sf K}(3) =
\sfF_1(3)$.  This is a particular instance the `source-target'
condition of in a globular category~\cite{batanin:globular}, see also
formula~(\ref{Pozitri_jedu_za_Jaruskou.}) in Section~\ref{jaRcA}.

In the rest of this subsection $T$, $S$
and $\nu : T \to S$ will have the same meaning as above.
Let us try to determine the value $\pa(g_T)$ of the differential on
the generator $g_T \in G(6)$ corresponding to the Tamarkin tree.
Inspired by~(\ref{hyrim_barvami}), define
\[
\pareg(g_T) := \textstyle\sum_\sigma \pm \iota(C_\sigma),
\]
with the sum taken over all `regular' faces $\sigma$ of $T$,
i.e.~faces $\sigma$ such that $\dim(C_\sigma) = \dim(T)-1$. In
Figure~\ref{!Jaruska}, the union of these faces is denoted
$\pareg \Tam$. Since $\pareg(g_T)$ is of dimension $5 < \dcrit(6)$, the
value $\pa (\pareg (g_T))$ is determined by calculations of
Subsection~\ref{Zitra_zubar} and equals the sum of
elements of $\Free(G)$ corresponding to the $4$-dimensional cells in
the intersection $\pareg \Tam \cap \pasing \Tam$ marked by two bullets
$\bullet$ in Figure~\ref{snad_mne_miluje}. In particular, $\pa (\pareg
(g_T)) \not= 0$. We shall find a `counterterm' $\pasing(g_T)$ such
that $\pa(\pareg (g_T)) = - \pa(\pasing (g_T))$ and put 
\[
\pa(g_T) := \pareg (g_T) +\pasing (g_T).
\]

The idea of finding such a couterterm is clear from
Figure~\ref{!Jaruska}; $\pasing (g_T)$ shall correspond to an
union of $5$-dimensional cells $U$ such that $\pa U = \pareg \Tam \cap
\pasing \Tam$. In the ideological Figure~\ref{snad_mne_miluje}, there
are two such unions, the `upper' $U'$ and the `lower' $U''$ which
indicates that the choice of $U$ need not be unique.  

The first step is to identify $4$-dimensional cells in 
the intersection $\pareg \Tam \cap
\pasing \Tam$. The extended barcodes of these cells
will be of the form $[b_1|b_2]$ for some  barcodes
$b_1,b_2$. To shorten the formulas, we use an `additive'
notation for the corresponding cells, so that $\mu[b_1' \cup b_1''|
b_2] =\mu [b_1'|b_2] \cup \mu[b_1''|b_2]$,~\&c. With this notation, one
easily expresses $\pareg \Tam \cap \pasing$ as the union of
$4$-dimensional cells:
\begin{eqnarray*}
\pareg \Tam \cap \pasing \Tam &=&
\muu{
\bigcup_{\tau \in \Sigma_{1,3,5}}[\tau_1|\tau_3|\tau_5] 
\cup
\bigcup_{\tau \in \Sigma_{1,3,5}'}\jarka{\tau_1}{[\tau_3||\tau_5]} 
\cup
\bigcup_{\tau \in \Sigma_{1,3,5}''}\jarka{[\tau_1||\tau_3]\hskip .1em}{\tau_5} 
}{[2||4||6]} 
\\
&\cup&
\muu{[1||3||5]\hskip .2em}{
\bigcup_{\tau \in \Sigma_{2,4,6}}[\tau_2|\tau_4|\tau_6] 
\cup
\bigcup_{\tau \in \Sigma_{2,4,6}'}\jarka{\tau_2}{[\tau_4||\tau_6]} 
\cup
\bigcup_{\tau \in \Sigma_{2,4,6}''}\jarka{[\tau_2||\tau_4]\hskip .1em}{\tau_6} 
}
\\
&\cup&
\muu{\Jarka{[1||3]\hskip .1em}5 \cup
  \bigcup_{\tau\in\Sigma_{1,3}}[\tau_1|\tau_3||5] \hskip .2em}
{\Jarka{[2||4]}6 \cup
  \bigcup_{\tau\in\Sigma_{2,4}}[\tau_2|\tau_4||6]}
\\
&\cup&
\muu{\Jarka1{[3||5]} \cup 
   \bigcup_{\tau\in\Sigma_{3,5}}[1||\tau_3|\tau_5]\hskip .2em}
{\Jarka {2\hskip .1em}{[4||6]} \cup
  \bigcup_{\tau\in\Sigma_{4,6}}[2\hskip .1em||\tau_4|\tau_6]}.
\end{eqnarray*}
In the above display, $\Sigma_{1,3,5}$ is the group of permutations of
the set $\{1,3,5\}$, and the symbols $\Sigma_{2,4,6}$, $\Sigma_{1,3}$,
$\Sigma_{2,4}$, $\Sigma_{3,5}$ and $\Sigma_{4,6}$ have the
obvious  similar meanings. Moreover,
\begin{align*}
&\Sigma'_{1,3,5} := \{ \tau \in \Sigma_{1,3,5};\ \tau_3 < \tau_5\},\ 
\Sigma''_{1,3,5} := \{ \tau \in \Sigma_{1,3,5};\ \tau_1 < \tau_3\},
\\
&\Sigma'_{2,4,6} := \{ \tau \in \Sigma_{2,4,6};\ \tau_4 < \tau_6\}
\mbox { and } 
\Sigma''_{2,3,6} := \{ \tau \in \Sigma_{2,4,6};\ \tau_2 < \tau_4\}.
\end{align*}
The structure of the right hand side should be clear from
Figure~\ref{!Jarka!} which shows, without specifying the labels,
generic points of the corresponding configurations. 
\begin{figure}[ht]
\begin{eqnarray*}
\raisebox{-65pt}{\rule{0pt}{0pt}}
&
\Bbox1
\lineup\Miki\Vojacek  \cup \lineup\hFlicek\Vojacek  \cup
\lineup\dFlicek\Vojacek   
\cup 
\Bbox2\lineup\Vojacek\Miki  \cup \lineup\Vojacek\hFlicek  \cup
\lineup\Vojacek\dFlicek&
\\
\raisebox{-65pt}{\rule{0pt}{0pt}}
&
\Abox3
\lineup\lKulda\lKulda \cup \lineup\lKrtek\lKulda
\cup \lineup\lKulda\lKrtek  \cup \lineup\lKrtek\lKrtek
\cup
\Abox4
\lineup\pKulda\pKulda \cup \lineup\pKrtek\pKulda
\cup \lineup\pKulda\pKrtek  \cup \lineup\pKrtek\pKrtek&
\end{eqnarray*}
\caption{\label{!Jarka!}
4-dimensional cells in the intersection $\pareg \Tam \cap
  \pasing \Tam$.} 
\end{figure}

\begin{figure}[ht]
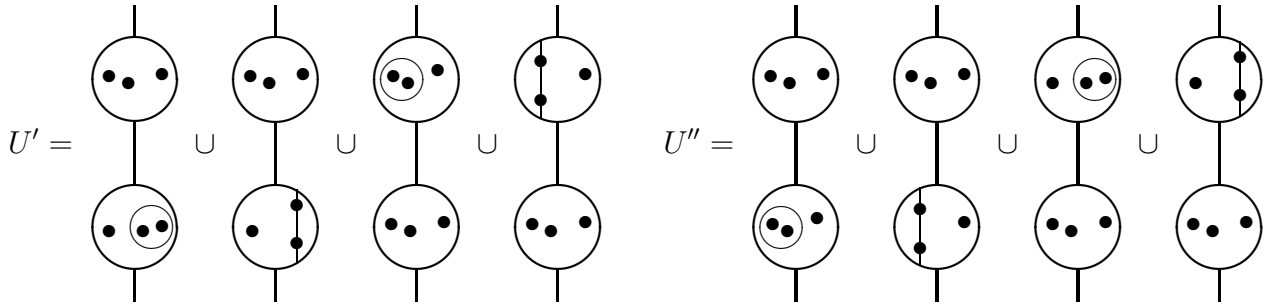

\[
\raisebox{-65pt}{\rule{0pt}{0pt}}
U' = \lineup\pKulda\Vojacek \cup \lineup\pKrtek\Vojacek \cup 
\lineup\Vojacek\lKulda \cup \lineup\Vojacek\lKrtek 
\hskip 20pt
U'' = \lineup\lKulda\Vojacek \cup \lineup\lKrtek\Vojacek \cup 
\lineup\Vojacek\pKulda \cup \lineup\Vojacek\pKrtek 
\]
\caption{\label{Jaruska-pusinka}
Two possible choices of the set $U$.}
\end{figure}

The four boxes of
this figure correspond to the four lines of the display.
One of the possible choices for the set $U$ is then
\[
U' := 
\muu
{\Jarka{1\hskip .1em}{[3||5]} \cup
  \bigcup_{\tau\in\Sigma_{3,5}}[1||\tau_3|\tau_5]\hskip .1em}{[2||4||6]}
\cup
\muu{[1||3||5]\hskip .2em}{\Jarka{[2||4]\hskip .1em}6 
\cup  \bigcup_{\tau\in\Sigma_{2,4}}
  [\tau_2|\tau_4||6]}.
\]
Another choice is the `diagonal image'
\[
U'' := 
\muu{\Jarka{[1||3]\hskip .1em}5 \cup
  \bigcup_{\tau\in\Sigma_{1,3}}[\tau_1|\tau_3||5]\hskip .1em}{[2||4||6]}
\cup
\muu{[1||3||5]\hskip .2em}{\Jarka{2\hskip .1em}{[4||6]} 
\cup  \bigcup_{\tau\in\Sigma_{4,6}}
  [2||\tau_4|\tau_6]}.
\]
Generic points of the corresponding cells are shown in
Figure~\ref{Jaruska-pusinka}.
In both cases, the counterterm $\pasing(g_T)$ is the sum of $6$ terms
corresponding to the six $4$-cells of $U'$ resp.~$U''$. 
 
\subsection{A $4$-dimensional bad cell}
\label{Jarca!!!}
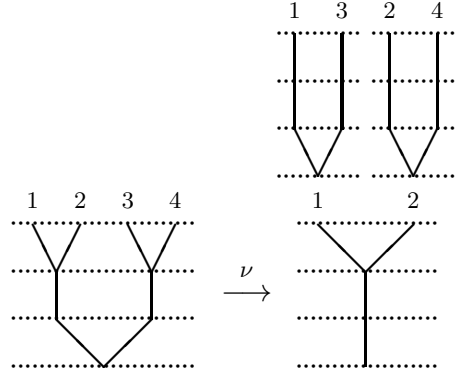
\begin{figure}[ht]
\begin{center}
{
\unitlength=.900000pt
\begin{picture}(180.00,150.00)(0.00,0.00)
\thicklines
\put(180.00,150.00){\makebox(0.00,0.00){\scriptsize $4$}}
\put(160.00,150.00){\makebox(0.00,0.00){\scriptsize $2$}}
\put(140.00,150.00){\makebox(0.00,0.00){\scriptsize $3$}}
\put(120.00,150.00){\makebox(0.00,0.00){\scriptsize $1$}}
\put(170.00,70.00){\makebox(0.00,0.00){\scriptsize $2$}}
\put(130.00,70.00){\makebox(0.00,0.00){\scriptsize $1$}}
\put(70.00,70.00){\makebox(0.00,0.00){\scriptsize $4$}}
\put(50.00,70.00){\makebox(0.00,0.00){\scriptsize $3$}}
\put(30.00,70.00){\makebox(0.00,0.00){\scriptsize $2$}}
\put(10.00,70.00){\makebox(0.00,0.00){\scriptsize $1$}}
\put(100.00,40.00){\makebox(0.00,0.00){\scriptsize $\nu$}}
\put(100.00,30.00){\makebox(0.00,0.00){$\longrightarrow$}}
\multiput(0,0)(40,0){2}{
\multiput(113.5,80)(0,20){4}{
\multiput(0,0)(3,0){12}{\makebox(0,0){$\cdot$}}
}
}
\multiput(2,0)(0,20){4}{
\multiput(0,0)(3,0){26}{\makebox(0,0){$\cdot$}}
}
\multiput(122,0)(0,20){4}{
\multiput(0,0)(3,0){20}{\makebox(0,0){$\cdot$}}
}
\put(150.00,40.00){\line(0,-1){40.00}}
\put(180.00,100.00){\line(-1,-2){10.00}}
\put(180.00,140.00){\line(0,-1){40.00}}
\put(160.00,100.00){\line(0,1){40.00}}
\put(170.00,80.00){\line(-1,2){10.00}}
\put(150.00,40.00){\line(1,1){20.00}}
\put(130.00,60.00){\line(1,-1){20.00}}
\put(140.00,100.00){\line(0,1){40.00}}
\put(130.00,80.00){\line(1,2){10.00}}
\put(120.00,100.00){\line(0,1){40.00}}
\put(130.00,80.00){\line(-1,2){10.00}}
\put(60.00,40.00){\line(1,2){10.00}}
\put(60.00,40.00){\line(-1,2){10.00}}
\put(60.00,20.00){\line(0,1){20.00}}
\put(40.00,0.00){\line(1,1){20.00}}
\put(20.00,40.00){\line(1,2){10.00}}
\put(10.00,60.00){\line(1,-2){10.00}}
\put(20.00,20.00){\line(1,-1){20.00}}
\put(20.00,40.00){\line(0,-1){20.00}}
\end{picture}}
\end{center}
\caption{\label{Mar}
The tree $T$ (left), the tree $S$ (right
bottom), the map $\nu : T \to S$ and its fiber diagram.}
\end{figure}

\begin{figure}[ht]
\begin{center}
{
\unitlength=1pt
\begin{picture}(150.00,120.00)(0.00,0.00)
\put(127.00,83.00){\makebox(0.00,0.00)[rb]{\scriptsize $4$}}
\put(77.00,33.00){\makebox(0.00,0.00)[rb]{\scriptsize $3$}}
\put(57.00,103.00){\makebox(0.00,0.00)[rb]{\scriptsize $2$}}
\put(27.00,73.00){\makebox(0.00,0.00)[rb]{\scriptsize $1$}}
\put(130.00,80.00){\makebox(0.00,0.00){$\bullet$}}
\put(80.00,30.00){\makebox(0.00,0.00){$\bullet$}}
\put(60.00,100.00){\makebox(0.00,0.00){$\bullet$}}
\put(30.00,70.00){\makebox(0.00,0.00){$\bullet$}}
\put(150.00,100.00){\line(-1,-1){10.00}}
\put(50.00,0.00){\line(1,1){100.00}}
\put(0.00,40.00){\line(1,1){80.00}}
\thicklines
\put(20.00,10.00){\vector(1,1){60.00}}
\put(20.00,10.00){\vector(0,1){90.00}}
\put(20.00,10.00){\vector(1,0){90.00}}
\end{picture}}
\end{center}
\caption{\label{Jaruska_ma_padesat.}
A generic point of the cell $\Mar$. The points in the ambient $\RRR^3$
lie on two lines
parallel to the 3rd coordinate. The 2nd coordinate of the points
labeled $1$ and $2$ is less than the 2nd coordinate of the points
labeled $3$ and $4$.
}
\end{figure}
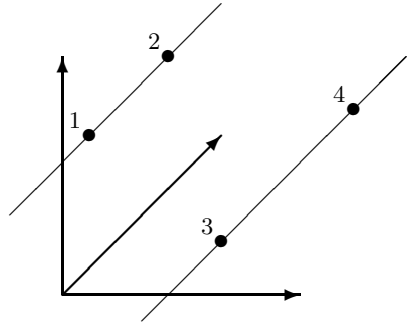

It is the cell \label{mar}$\Mar 
:= \mu[T]$ indexed by the reduced $3$-tree $T :=
[1|2|||3|4]$ shown in Figure~\ref{Mar}.
 A generic point of this cell is presented in
Figure~\ref{Jaruska_ma_padesat.}.
Consider the $3$-tree $S := [1|2]$ and the map $\nu : T \to S$ that
sends the tips of $T$ labeled 1,3 (reps.~2,4) to the tip
1 (resp.~tip 2) of $S$, see again Figure~\ref{Mar}. It is clear
that the corresponding face $\mu[\nu]$ of $\nu$ equals 
$\mu\jarka{[1|||3]\hskip .1em}{[2|||4]}$ and that
$\dim(\mu[\nu]) = \dim(\Mar) = 4$. 

Let us determine the value $\pa(g_T)$ of the differential on the
$4$-dimensional generator $g_T \in G_4(4)$ corresponding to $T$.  As in
Subsection~\ref{Jarca!!}, define
\[
\pareg(g_T) := \textstyle\sum_\sigma \pm \iota(C_\sigma),
\]
with the sum over all faces $\sigma$ such that $\dim(C_\sigma) =
3$. One easily sees that these faces form the union
\begin{eqnarray*}
\pareg \Mar &=& \mu [1|2||3|4] \cup \mu [3|4||1|2] \cup
\mu\JArka{1 \velkasvisla 2}
{[3|4]} \cup \mu\JArka{[1|2]}{3 \velkasvisla 4} 
\\ 
&\cup&
\bigcup_{\tau \in \Sigma_{1,3}}\mu\jarka{\tau_1 \velkasvisla
\tau_3 \hskip .1em}{[2|||4]} \cup
\mu\jarkA{1\hskip .1em}{[2|||3]}4 \cup \mu\jarkA{3\hskip
  .1em}{[1|||4]}2 
\bigcup_{\tau \in
\Sigma_{2,4}}\mu\jarka{[1|||3]\hskip .1em}{\tau_2  \velkasvisla \tau_4}.
\end{eqnarray*}
Observe that the first two terms give the linear
part of $\pal(g_T)$.
As in Subsection~\ref{Jarca!!} we need to describe $2$-dimensional cells in 
the intersection $\pareg \Mar \cap
\pasing \Mar$. The result is:
\begin{eqnarray*}
\pareg \Mar \cap \pasing \Mar &=&
\muu{\bigcup_{\tau \in \Sigma_{1,3}}[\tau_1|\tau_3]}{[2|||4]}
\cup
\muu{[1|||3]\hskip .2em}{\bigcup_{\tau \in \Sigma_{2,4}}[\tau_2|\tau_4]} 
\\
&& \cup\ \mu\jarka{[1||3]\hskip .1em}{[2||4]} 
\cup  \mu\jarka{[3||1]\hskip .1em}{[4||2]},
\end{eqnarray*}
where $\Sigma_{1,3}$ (resp.~$\Sigma_{2,4}$) is the group of permutations of
the set $\{1,3\}$ (resp.~ $\{1,3\}$).
One of the possible choices for the set $U$ of $3$-cells generating
the counterterm $\pasing(g_T)$ is then
\[
U' := \mu\jarka{[1|||3]\hskip .1em}{[4||2]} 
\cup  \mu\jarka{[1||3]\hskip .1em}{[2|||4]}.
\]
The second one is the `diagonal image'
\[
U'' := 
 \mu\jarka{[1|||3]\hskip .1em}{[2||4]} 
\cup  \mu\jarka{[3||1]\hskip .1em}{[2|||4]}.
\]
In both cases, the counterterm $\pasing(g_T)$ has $2$ terms
corresponding to two $3$-cells of $U'$ resp.~$U''$. The differential
$\pa(g_T)$ is the sum of $2$ linear terms, $8$ regular decomposable
terms and $2$ singular terms. As in Subsection~\ref{Jarca!!}, it helps
to represent the cells entering the above calculations by
depicting their generic points.  We leave it to the reader as an
exercise.

\subsection{Proof of Theorem~B}
\label{b}
Let $\bad$ be a set of unlabeled trees indexing (one or more) 
cells in the critical dimension
(such as $\Tam$ or $\Mar$ above). Suppose that, for each $T \in \bad$,
we found an element $\papert(g_T) \in \Free^{\geq 2}(G)$ such that 
\[
\pa(\pal + \papert)(g_T) = 0.
\]  
Since $T$ has the critical dimension,  
$(\pal + \papert)(g_T) \subset \Free(\Greg)$ and
the above formula makes sense. Examples of $\papert(g_T)$ are
given in Subsections~\ref{Jarca!!} and~\ref{Jarca!!!}. 

Now we take, in Lemma~\ref{Zitra_do_VO}, $(\Free(E),\vartheta) :=
\Omega(\Lie')$, $(M,\pal) := (G,\pal)$ and as $\opa = \opal +
\opapert$ we take the differential from
Proposition~\ref{JarunA}. Finally, let $A(n) := \Sigma_n[\bad(n)]$ be,
for $n \geq 4$,
the free $\Sigma_n$-module generated by trees $T \in \bad$ of arity
$n$. With these choices, the assumptions of the lemma are clearly
satisfied and Theorem~B follows.

\section{Bad cells}
\label{jaRcA}

In the first part of this section we analyze the `source-target' conditions
responsible for the existence of bad cells. In the second part we
prove Propositions~\ref{Jarka_byla_u_mne!} and~\ref{JArkA}.  

\subsection{The source-target conditions}
\label{Kremilek}
Assume we are given a pruned unlabeled $h$-tree $T \in \uTree^h(n)$ as
in~(\ref{eq:1}). Suppose that there is an $s \geq 2$ and natural
numbers
\[
1 \leq a_1 < b_1 < a_2 < b_2 \cdots < a_s < b_s \leq k_h = n
\]
and some $1 \leq m < h$ such that
\begin{equation}
\label{v_Mulhouse}
\rho_m \circ \cdots \circ \rho_{h-1}(a_i) = 
\rho_m \circ \cdots \circ \rho_{h-1}(b_i),
\end{equation}
for all $1 \leq i \leq s$.  We also assume that the common values of
the expression in~(\ref{v_Mulhouse}) form a strictly increasing
sequence of $s$ elements of $[k_m]$.  Suppose there is a $h$-tree $S
\in \uTree^h(k)$, $k < n$, and a map $\nu : T \to S$ for which there
exist $1 \leq u < v \leq k$ such that $\nu_h(a_i) = u$ and
$\nu_h(b_i) = v$ for all $1 \leq i \leq s$.

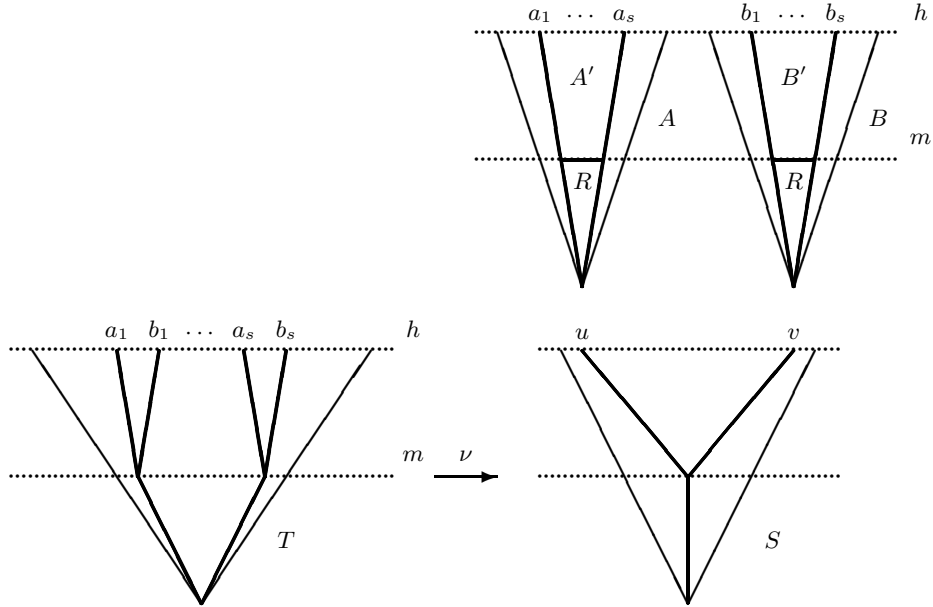
\begin{figure}[ht]
{
\unitlength=.800000pt
\begin{picture}(430.00,280.00)(0.00,0.00)
\thicklines
\put(370.00,200.00){\makebox(0.00,0.00){\scriptsize $R$}}
\put(270.00,200.00){\makebox(0.00,0.00){\scriptsize $R$}}
\put(370.00,250.00){\makebox(0.00,0.00){\scriptsize $B'$}}
\put(270.00,250.00){\makebox(0.00,0.00){\scriptsize $A'$}}
\put(360.00,30.00){\makebox(0.00,0.00){\scriptsize $S$}}
\put(370.00,275){\makebox(0.00,0.00)[b]{\scriptsize $\cdots$}}
\put(270.00,275){\makebox(0.00,0.00)[b]{\scriptsize $\cdots$}}
\put(390.00,275){\makebox(0.00,0.00)[b]{\scriptsize $b_s$}}
\put(350.00,275){\makebox(0.00,0.00)[b]{\scriptsize $b_1$}}
\put(290.00,275){\makebox(0.00,0.00)[b]{\scriptsize $a_s$}}
\put(250.00,275){\makebox(0.00,0.00)[b]{\scriptsize $a_1$}}
\put(430.00,280.00){\makebox(0.00,0.00){\scriptsize $h$}}
\put(430.00,220.00){\makebox(0.00,0.00){\scriptsize $m$}}
\put(410.00,230.00){\makebox(0.00,0.00){\scriptsize $B$}}
\put(310.00,230.00){\makebox(0.00,0.00){\scriptsize $A$}}
\put(370.00,125.00){\makebox(0.00,0.00)[b]{\scriptsize $v$}}
\put(270.00,125.00){\makebox(0.00,0.00)[b]{\scriptsize $u$}}
\put(215.00,70.00){\makebox(0.00,0.00){\scriptsize $\nu$}}
\put(130.00,30.00){\makebox(0.00,0.00){\scriptsize $T$}}
\put(90.00,125){\makebox(0.00,0.00)[b]{\scriptsize $\cdots$}}
\put(130,125){\makebox(0.00,0.00)[b]{\scriptsize $b_s$}}
\put(110.00,125){\makebox(0.00,0.00)[b]{\scriptsize $a_s$}}
\put(70.00,125){\makebox(0.00,0.00)[b]{\scriptsize $b_1$}}
\put(50.00,125){\makebox(0.00,0.00)[b]{\scriptsize $a_1$}}
\put(190.00,70.00){\makebox(0.00,0.00){\scriptsize $m$}}
\put(190.00,130.00){\makebox(0.00,0.00){\scriptsize $h$}}
\put(200.00,60.00){\vector(1,0){30.00}}
\multiput(0,0)(0,60){2}{
\multiput(250.00,60.00)(3,0){48}{\makebox(0,0){$\cdot$}}
}
\multiput(-.25,0)(.050,0){10}{
\put(320.00,60.00){\line(0,-1){60.00}}
\put(320.00,60.00){\line(5,6){50.00}}
\put(320.00,60.00){\line(-5,6){50.00}}
}
\put(260.00,120.00){\line(1,-2){60.00}}
\put(320.00,0.00){\line(1,2){60.00}}
\multiput(0,0)(0,60){2}{
\multiput(220.00,210.00)(3,0){67}{\makebox(0,0){$\cdot$}}
}
\multiput(-.25,0)(.050,0){10}{
\put(370.00,150.00){\line(1,6){20.00}}
\put(350.00,270.00){\line(1,-6){20.00}}
\put(270.00,150.00){\line(1,6){20.00}}
\put(270.00,150.00){\line(-1,6){20.00}}
}
\multiput(0,-.25)(0,.050){10}{
\put(260.00,210.00){\line(1,0){20}}
\put(360.00,210.00){\line(1,0){20}}
}
\put(370.00,150.00){\line(1,3){40.00}}
\put(370.00,150.00){\line(-1,3){40.00}}
\put(270.00,150.00){\line(1,3){40.00}}
\put(270.00,150.00){\line(-1,3){40.00}}
\multiput(0,0)(0,60){2}{
\multiput(0.00,60.00)(3,0){61}{\makebox(0,0){$\cdot$}}
}
\multiput(-.25,0)(.050,0){10}{
\put(90.00,0.00){\line(1,2){30.00}}
\put(60.00,60.00){\line(1,-2){30.00}}
\put(60.00,60.00){\line(-1,6){10.00}}
\put(70.00,120.00){\line(-1,-6){10.00}}
\put(120.00,60.00){\line(1,6){10.00}}
\put(110.00,120.00){\line(1,-6){10.00}}
}
\put(90.00,0.00){\line(-2,3){80.00}}
\put(90.00,0.00){\line(2,3){80.00}}
\end{picture}}
\caption{\label{Jarka_mi_udelala_svacinku!}
The origin of the source-target conditions -- schematic picture.}
\end{figure}

In this situation, denote $A$ (resp.~$B$) the fiber of $\nu$ over $u$
(resp.~$v$) and $A'$ (resp.~$B'$) the maximal pruned subtree
of $A$ (resp.~$B$) with the tips $\Rada a1s$ (resp.~$\Rada
b1s$). Denote finally $R \in \uTree^m(s)$ the pruned unlabeled 
$m$-tree obtained from
$A'$ by amputating everything above level $m$. 
Observe that instead of amputating $A'$ we could have
amputated $B'$ with the same result. 
The situation is visualized in
Figure~\ref{Jarka_mi_udelala_svacinku!}.

The tree $R$ determines a cell $\mu[R] \subset \sfF_m(s)$.
The {\em source and target maps\/} $\pi_s,\pi_t : \mu[C_\nu] \epi
\mu[R]$ are defined as follows.  Since $\mu[C_\nu]$ is the
cartesian product of the cell $\mu[S]$ with the cells indexed by the
fibers of $\nu$, one has the projections $\pi_A :
\mu[C] \epi \mu[A]$ resp.~ $\pi_B : \mu[C] \epi \mu[B]$. One also has
the `forgetful' projections $\pi_A' : \mu[A] \epi \mu[A']$ (resp.~$\pi_B'
: \mu[B] \epi \mu[B']$) given by forgetting all points of the
configurations in $\mu[A]$ (resp.~$\mu[B]$) except those with labels
in $\{\Rada a1s\}$ (resp.~$\{\Rada b1s\}$). Let finally $\pi''_A : \mu[A']
\epi \mu[R]$ (resp.~$\pi''_B : \mu[B'] \epi \mu[R]$) be the projection induced
by the projection $\RRR^h \epi \RRR^m$ to the first $m$ coordinates.
The maps $\pi_s,\pi_t$ are the compositions
\[
\pi_s : \mu[C_\nu] \stackrel{\pi_A}\epi \mu[A] \stackrel{\pi_A'}\epi
\mu[A']  \stackrel{\pi''_A}\epi \mu[R]\
\mbox { and }\
\pi_s : \mu[C_\nu] \stackrel{\pi_B}\epi \mu[B] \stackrel{\pi_B'}\epi
\mu[B']  \stackrel{\pi''_B}\epi \mu[R].
\]

\begin{definition}
The {\em source-target condition\/} is the following equality of 
points of $\mu[R]$
\begin{equation}
\label{Pozitri_jedu_za_Jaruskou.}
\pi_s(\xxx) = \pi_t(\xxx) 
\end{equation}
satisfied by each point $\xxx \in \mu[C_\nu] \cap \pa \mu[T]$.
\end{definition}

\begin{example}
{\rm 
In the Tamarkin case analyzed in Subsection~\ref{Jarca!!}, $h=2$,
$n=6$, $k=2$, $m=1$, the trees $T \in \uTree^2(6)$, $S \in
\uTree^2(2)$ and the map $\nu : T \to S$ are as in
Figure~\ref{Tam}. Moreover, $s=3$, $a_1=1$, $a_2=3$, $a_3=5$, $b_1=2$,
$b_2=4$, $b_3=6$. The amputated tree $R$ equals $[1|2|3]$ so $\mu[R]$
is the open interval $\osfF_1(3)$.

For the $4$-dimensional bad cell of Subsection~\ref{Jarca!!!}, $h=3$,
$n=4$, $k=2$, $m=2$, the trees $T \in \uTree^3(4)$, $S \in
\uTree^3(2)$ and the map $\nu : T \to S$ are as in
Figure~\ref{Mar}. Moreover, $s=2$, $a_1=1$, $a_2=3$, $b_1=2$, and
$b_2=4$. The amputated tree $R$ equals $[1||2]$ so $\mu[R]$ is the
open half-circle in $\sfF_2(2) = \sphere^1$.  
}\end{example}

\subsection{Proofs of Propositions~\ref{Jarka_byla_u_mne!} and~\ref{JArkA}}

Let us introduce the
following terminology. We call a cell of $\sfF$ {\em regular\/} if its
boundary is an union of cells. For subsets $a \subset \sfF(m)$, $b \subset
\sfF(n)$ and $1 \leq i \leq m$ we write, as expected,
\[
a \circ_i b := \{x\circ_i y \in \sfF(m+n-1);\ x \in a,\ y \in b\}. 
\]
We call $a \circ_i b$ the {\em $\circ_i$-composition\/} of the sets
$a$ and $b$.

\begin{lemma}
\label{JaruSka}
Let $\sigma \in \Sigma_m$.  A cell $e \subset \sfF(m)$ is a regular
if and only if the cell $e \nu := \{x\sigma;x \in e\}$ is
regular. The $\circ_i$-composition of cells is regular if and only if and only
if both factors are regular.
\end{lemma}

\begin{proof}
The first part of the lemma is obvious. The second part
follows from the equality $\hbox{$\pa (e' \circ_i e'')$} = \pa e' \circ_i e''
\cup e' \circ_i \pa e''$.
\end{proof}

\begin{theorem}
\label{Prijde_dnes_Jarka?}
For $n \geq 2$ and $h \geq 1$, each cell of dimension $< \dcrit^h(n)$
in the open part $\osfF_h(n)$ is regular and all its faces are
regular, too. Similarly, each cell of $\osfF(n)$ of dimension $<
\dcrit(n)$ is regular along with all its faces.
\end{theorem}

\begin{proof}
The theorem can be proved
by case-studying cells of the indicated dimensions. We, however, prefer a
conceptual approach based on the analysis of the source-target
conditions given in Subsection~\ref{Kremilek}. In particular, we observe
that condition~(\ref{Pozitri_jedu_za_Jaruskou.}) is nontrivial only if the
dimension of the cell $\mu[R]$ is at least one. This implies that
\begin{itemize}
\item[(i)] 
either $s \geq 2$ and $m \geq 2$, or 
\item[(ii)]
$s \geq 3$ and $m \geq 1$.
\end{itemize}

In case (i), one has $n\geq 4$ and $h \geq 3$, while 
in case (ii) one has $n \geq 6$ and $h
\geq 2$. 
Let us prove that each $e \subset \osfF_h(n)$ with $\dim(e) <
\dcrit^h(n)$ is regular. We again distinguish two cases:

\begin{itemize}
\item[(a)] 
$\dcrit^h(n) = \infty$, which means $n = 2,3$, or $n = 4,5$ and $h\leq 2$, or
$n \geq 6$ and $h=1$,
\item[(b)]
$\dcrit^h(n) = n$, which means $n=4,5$ and $h \geq 3$, or $n \geq 6$ and
  $h \geq 2$.
\end{itemize} 

Case (a) is complementary to the cases (i) and (ii) above, therefore
the source-target conditions are trivial and $e$ is a regular
cell. The faces of the cell $e$ are $\circ_i$-compositions of cells $e'$
from $\osfF_{h}(n')$ for some $2 \leq n' \leq n$.
Since, by the definition of the critical dimension, $\dcrit^h(n') \geq
\dcrit^h(n)$,  one has
$\dcrit^{h}(n') = \infty$. As we already established, this implies
that each such an
$e'$ is regular, therefore, by Lemma~\ref{JaruSka}, each face of $e$
is regular, too.

Let us assume~(b), i.e.~$\dcrit^h(n) = n$. Since, by Lemma~\ref{JaruSka},
the symmetric group action preserves the regularity, we may assume that
$e = \mu[T]$ for a reduced
tree $T \in \uTree^h(n)$ as in~(\ref{eq:1}) with $\dim(T) < n$.
By simple combinatorics, 
\begin{itemize}
\item[($\alpha$)] 
either $h = 1$ or
\item[($\beta$)] 
$h = 2$ and $k_1 = 2$.
\end{itemize}
In the first case the cell $\mu[T]$ belongs to the Stasheff polytope
${\sf K}(n) = \sfF_1(n) \subset \sfF(n)$, so $\mu[T]$ and all its
faces are regular cells.  In case ($\beta$), the amputated tree $R$
must be $[1|2]$, thus the corresponding source-target condition is
trivial, so $\mu[T]$ is a regular cell. It is not difficult to verify
that if $T$ satisfies ($\beta$), $S$ is an arbitrary tree and $\nu : T
\to S$ a map, then each reduced fiber of $\nu$ also satisfies
($\alpha$) or ($\beta$) above. This implies that all faces of $\mu[T]$
are regular cells, too.

This finishes the proof of the first part of the theorem.  Since each
cell $e \subset \osfF(n)$ belongs to some $\osfF_h(n)$, $h \geq 1$,
the second part follows from the first part and the inequality
$\dcrit(n) \leq \dcrit^h(n)$.
\end{proof}

\begin{proof}[Proof of Proposition~\ref{Jarka_byla_u_mne!}]
The compatibility of the CW-structures of the $\Sigma$-modules
$\sfF_h$, $h \geq 1$, and $\sfF$ with the operad structures and the
freeness of the symmetric group action on the cells follows from the
very definition of the cell structure reviewed on
page~\pageref{Jarka_na_chalupe.}.
Recall~\cite{batanin:conf,sinha:man} that there is the canonical
embedding
\[
\iota:
\sfF_h(n) \hookrightarrow  
{\mbox {\LARGE $\times$}}_{\substack{1 \leq i,j \leq n \\ i\not= j}} \
\sphere^{h-1} \times {\mbox {\LARGE $\times$}}_{\substack{1 \leq i,j,k
    \leq n \\ i\not= j,\ j\not= k,\ k\not= i}}\ [0,\infty].
\]
It follows from the analysis of the images of the cells of
$\sfF_h(n)$ under $\iota$ given in Section~6 
of~\cite{batanin:conf}, namely from Proposition~6.1 of that
section, that the spaces $\sfF_h(n)$ satisfy condition~(i) of
Definition~\ref{JarunKa} for arbitrary~$n$ and $h$. The analogous
claim for $\sfF(n)$ stems from the fact that each cell of $\sfF(n)$
belongs to the subcomplex $\sfF_h(n)$ for some $h \geq 1$.

Since condition~(i) of Definition~\ref{JarunKa} has already
been established, the space $\sfF_h(n)$ is regular if and only if each
its cell is regular in the sense introduced at the beginning of this
subsection. Since the cells of $\sfF_h(n)$ are iterated
$\circ_i$-compositions of the cells from $\osfF_h(n')$ with $n' \leq
n$, the complex $\sfF_h(n)$ is, by Lemma~\ref{JaruSka}, regular if and
only if all cells of $\osfF_h(n')$ are regular, for each $2 \leq n'
\leq n$. Since clearly $\dcrit^h(n) = \infty$ implies $\dcrit^h(n') =
\infty$ for each $n' \leq n$, the space $\sfF_h(n)$ is, by
Theorem~\ref{Prijde_dnes_Jarka?}, regular if $\dcrit^h(n) = \infty$. 

The non-regularity of $\sfF_h(n)$ if $\dcrit^h(n)$ is finite follows
from Proposition~\ref{JArkA} proved below. This, by the remark following
Proposition~\ref{Jarka_byla_u_mne!}, proves the characterization of
the regularity of the spaces $\sfF_h(n)$.  The similar obvious
analysis applies to $\sfF(n)$ as well.

The regularity of the suboperads $\sfFreg_h$ and $\sfFreg$ follows
from the fact that they are generated by the regular cells in  
$\osfF_h$ resp.~$\osfF$, from Lemma~\ref{JaruSka} and the fact that
cells in the boundary of regular cells are regular established in
Theorem~\ref{Prijde_dnes_Jarka?}. This finishes the proof.
\end{proof}

\begin{proof}[Proof of Proposition~\ref{JArkA}]
We put $e_4^3 := \Mar = \mu[1|2|||3|4]$, the cell introduced in
Subsection~\ref{Jarca!!!}. The cell $e^h_4$ for $h \geq 4$ is defined
as the image of $e_4^3$ under the natural inclusion $\sfF_3(4)
\hookrightarrow \sfF_h(4)$. Likewise, we put  $e_5^3 :=
\mu[1|2|||3|4|5]$ and $e_5^h$ for $h \geq 4$ defined similarly.

The cell $e^2_6$ is the Tamarkin cell $\Tam = \mu[1|2||3|4||5|6]$, and
$e_n^2 := \mu[1|2||3|4||5|6|\cdots|n]$, for $n \geq 7$. The cells
$e_n^h$ for $h \geq 3$ and $n \geq 6$ are the images of $e_n^2$ under
the natural inclusions $\sfF_3(n)
\hookrightarrow \sfF_h(n)$. The cells $e_n$, $n \geq 4$, are the images of
$e_n^3$ under the inclusions $\sfF_3(n)
\hookrightarrow \sfF(n)$. We leave to the reader to verify that the
cells $e^h_n$ defined in this way are bad.
\end{proof}

\section{Free Lie algebras and configuration spaces}
\label{free_Lie}

In this section we prove integral variants of some results whose
characteristic zero versions are known.  Therefore, all algebraic
objects will be considered over the ring $\ZZZ$ of integers. 
The results below easily generalize to an arbitrary 
integral domain with unit.

Let $\TTT(\Rada x1n)$ be the tensor
algebra with generators $\Rada x1n$, $n \geq 1$, and $\LLL(\Rada x1n)$
the free Lie algebra on the same set of generators, considered in the
standard way as a subspace of $\TTT(\Rada x1n)$. We denote by
$\TTT_{\rada 11}(\Rada x1n) \subset \TTT(\Rada x1n)$ the linear
subspace spanned by words containing each of the generators $\Rada
x1n$ precisely once, and $\LLL_{\rada 11}(\Rada x1n) := \LLL(\Rada
x1n) \cap \TTT_{\rada 11}(\Rada x1n)$. We will sometimes simplify the
notation and denote $\Lie(n) := \LLL_{\rada 11}(\Rada x1n)$.

Each space above has a natural right $\Sigma_n$-module action
permuting the generators. Take another set of generators $\Rada \alpha
1n$ and denote by
\begin{equation}
\label{psano_na_letisti}
\Phi : \TTT_{\rada 11}(\Rada \alpha1n) \to \TTT_{\rada 11}(\Rada x1n)', 
\end{equation}
where $\TTT_{\rada 11}(\Rada x1n)'$ is the linear dual of
$\TTT_{\rada 11}(\Rada x1n)$, the isomorphism defined by
\[
\Phi(\alpha_{\rho(1)}\ot \cdots\ot \rho_{\sigma(n)})
(x_{\omega(1)}\ot\cdots\ot x_{\omega(n)}) 
:=
\cases{1}{if $\rho= \omega$, and}{0}{otherwise,} 
\]
where $\rho,\omega \in \Sigma_n$.

For $s,t \geq 1$ denote by $\Sigma_{s,t}$ is the set of all {\em
$(s,t)$-unshuffles\/}, i.e.~permutations 
$\tau \in \Sigma_n$, $n = s+t$, such that 
\[
\tau(1) < \cdots < \tau(s) \mbox { and } \tau(s+1) < \cdots < \tau(s+t).
\]
Let, finally,  $\Ush_{\rada
11}(\Rada \alpha1n) \subset \TTT_{\rada 11}(\Rada \alpha1n)$ be the
linear span of the expressions 
\[
\sum_{\tau \in
\Sigma_{s,t}} \alpha_{\rho\tau(1)}\ot \cdots\ot \alpha_{\rho\tau(n)},
\]
for $\rho \in \Sigma_n$ and $s,t \geq 1$ such that  $s+t=n$.

\begin{theorem}
\label{Jarca}
The map~(\ref{psano_na_letisti}) induces a $\Sigma_n$-equivariant isomorphism 
\[
\underline\Phi : \TTT_{\rada 11}(\Rada \alpha1n)/\Ush_{\rada 11}(\Rada
\alpha1n) \to \LLL_{\rada 11}(\Rada x1n)'.
\]
\end{theorem}

Theorem~\ref{Jarca} will follow from a sequence of claims. The first
one is probably known, but we were unable to find a reference.

\begin{claim}
\label{JARCa}
For each $n \geq 1$, 
$\Lie(n) = \LLL_{\rada 11}(\Rada x1n)$ is the free abelian group
with basis
\begin{equation}
\label{zase}
b_\lambda :=
[x_{\lambda(1)},[x_{\lambda(2)},\ldots,[x_{\lambda(n-1)},x_n]\cdots]],\
\lambda \in \Sigma_{n-1}.
\end{equation}
\end{claim}

\begin{proof}
It is known that $\LLL(\Rada x1n)$ is torsion
free (see, for instance, the overview in~\cite{blessenohl-laue}), so its
subspace  $\Lie(n)$ is torsion-free as well.
Let us prove, by induction on $n$, that the elements in~(\ref{zase})
span $\Lie(n)$. 

This statement is obvious for $n =1,2$. Let $n >
2$. Since $\Lie(n)$ is spanned by elements of the form $[F_1,F_2]$,
$F_1 \in \Lie(s)$, $F_2 \in \Lie(t)$, $s+t=n$, $s,t \geq 1$, it
suffices to prove that each such $[F_1,F_2]$ is a linear combination of
elements of the basis~(\ref{zase}). We may clearly assume that $F_1$
contains $\Rada x1s$ and $F_2$ contains $\Rada x{s+1}n$. 
By induction, $F_2$ is a linear combination of iterated
commutators as in~(\ref{zase}), with $x_n$ at the extreme right
position.

Now we proceed by induction on $s$. If $s=1$,  $[F_1,F_2]$ is an
element of~(\ref{zase}). If $s \geq 2$ we may assume
that $F_1 = [A,B]$, for some $A \in \Lie(a)$,  $B \in \Lie(b)$, with
$a+b = s$, $a,b \geq 1$. In that case
\[
[F_1,F_2] = [A,[B,F_2]] + [B,[F_2,A]].
\]
By induction on $n$, both $[B,F_2]$ and $[F_2,A]$ are combinations of
the commutators as in~(\ref{zase}) with $x_n$ at the rightmost place,
therefore both $[A,[B,F_2]]$ and $[B,[F_2,A]]$ are
combination of basic elements~(\ref{zase}), by induction on $s$.

The linear independence of the elements~(\ref{zase}) 
follows from the well-known fact
that the dimension of $\Lie(n)$ is $(n-1)!$~\cite{blessenohl-laue}, 
which is the
number of elements~(\ref{zase}). 
\end{proof}

There is a straightforward way to verify the linear
independence of the elements~(\ref{zase}) based on Claim~\ref{clm}
below which will be useful also for other purposes.
Recall that $\TTT_{\rada 11}(\Rada x1n)$ is the free abelian group
with basis
\begin{equation}
\label{JarcA}
e_\sigma := x_{\sigma(1)} \ot \cdots \ot  x_{\sigma(n)},\
\sigma \in \Sigma_n. 
\end{equation}

\begin{claim}
\label{clm}
In a (unique) decomposition of \ $b_\lambda$, $\lambda \in \Sigma_{n-1}$,  
into a linear combination of
the basis $\{e_\sigma\}_{\sigma \in \Sigma_n}$, the element
\[
e_{\lambda \times 1} = x_{\lambda(1)} \ot \cdots \ot  x_{\lambda(n-1)}
\ot x_n
\]
appears with coefficient $1$ and $b_\lambda$ is the only basis
element~(\ref{zase}) whose decomposition contains $e_{\lambda \times 1}$.
\end{claim}

\begin{proof}
A simple induction on $n$.
\end{proof}

The following proposition is based on famous Theorem~2.2
of~\cite{ree:AnM69} that, however, assumes the existence of a solution
$\xi \in R$ of the equation $n\xi = \alpha$ for each natural $n \geq
1$ and each $\alpha \in R$, in the ground ring $R$.  
We show that this assumption is not
necessary when this theorem is applied to the subspace $\TTT_{\rada
11}(\Rada x1n)$ not to the whole $\TTT(\Rada x1n)$.

\begin{proposition}
\label{JaRCa}
An element $F = \sum_{\sigma \in \Sigma_n} a_\sigma\cdot x_{\sigma(1)} \ot
\cdots \ot x_{\sigma(n)} \in \TTT_{\rada 11}(\Rada x1n)$, $a_\sigma
\in \ZZZ$, belongs to the subspace $\LLL_{\rada 11}(\Rada x1n)$ if and
only if
\begin{equation}
\label{JaRcA}
\textstyle\sum_{\tau \in \Sigma_{s,t}} a_{\rho\tau} = 0,
\end{equation}
for each permutation $\rho \in \Sigma_n$ and each $s,t \geq 1$ such
that $s+t=n$.
\end{proposition}

\begin{proof}
By analyzing the proof of~\cite[Theorem~2.2]{ree:AnM69}, one sees
that~(\ref{JaRcA}) 
in fact implies $nF \in  \LLL_{\rada 11}(\Rada x1n)$. 
By Claim~\ref{JARCa} this means that $nF = \sum_{\lambda \in
\Sigma_{n-1}} \beta_\lambda \cdot b_\lambda$, for some $\beta_\lambda \in \ZZZ$
and $b_\lambda$ the commutators in~(\ref{zase}).  On the other
hand, it follows from Claim~\ref{clm} that, 
in the expression 
\[
nF = \textstyle\sum_{\sigma \in
\Sigma_n} n a_\sigma \cdot x_{\sigma(1)} \ot \cdots \ot x_{\sigma(n)},
\]
one has
$na_{\lambda \times 1} = \beta_\lambda$, for each $\lambda \in
\Sigma_{n-1}$. This means that $F = \sum_{\lambda \in
\Sigma_{n-1}} a_{\lambda \times 1}\cdot b_\lambda$, so $F$ is a Lie element.
\end{proof}

Another important piece of the proof of Theorem~\ref{Jarca} is

\begin{claim}
\label{jarCa}
The restriction $r :  \TTT_{\rada 11}(\Rada x1n)'
\to \LLL_{\rada 11}(\Rada x1n)'$ is an epimorphism.
\end{claim}

\begin{proof}
We need to show that an arbitrary linear map $\varphi :
\LLL_{\rada 11}(\Rada x1n)\to \ZZZ$ extends into a linear map
$\widetilde\varphi : \TTT_{\rada 11}(\Rada x1n)\to \ZZZ$.  Let
$\{b_\lambda\}_{\lambda \in \Sigma_{n-1}}$ be the basis~(\ref{zase})
of $ \LLL_{\rada 11}(\Rada x1n$), $\{e_\sigma\}_{\sigma \in
\Sigma_{n}}$ the basis~(\ref{JarcA})
of $ \TTT_{\rada 11}(\Rada x1n)$ and put 
\[
\widetilde\varphi(e_\sigma) :=
\cases{\phi(b_\lambda)}{if $\sigma = \lambda \times 1$ for some $\lambda \in
  \Sigma_{n-1}$, and}0{otherwise.}
\]
By Claim~\ref{clm}, $\widetilde\varphi$ defined in this way extends $\varphi$.
\end{proof}

Let $K$ denote the kernel of the composition 
\[
\TTT_{\rada 11}(\Rada \alpha1n) \stackrel\Phi\longrightarrow 
\TTT_{\rada 11}(\Rada x1n)'
\stackrel r\longrightarrow \LLL_{\rada 11}(\Rada x1n)', 
\]
where $\Phi$ is as in~(\ref{psano_na_letisti}) and $r$ the
restriction.

\begin{claim}
\label{jarcA}
An element $x \in \TTT_{\rada 11}(\Rada \alpha1n)$ belongs to the kernel
$K$ if and only if there exists a natural $N$ such that $N \cdot x
\in \Ush_{\rada 11}(\Rada \alpha1n)$.
\end{claim}

\begin{proof}
It is an elementary consequence of Proposition~\ref{JaRCa} that
\[
K= (\Ush_{\rada 11}(\Rada \alpha1n)^\perp)^\perp 
\supset \Ush_{\rada 11}(\Rada \alpha1n), 
\]
where ${}^\perp$ denotes the annihilator in the dual space. It is
another elementary fact that, for any subspace $A$ of a
finite-dimensional vector space $V$, one has $(A^\perp)^\perp \cong
A$, therefore, after extending the scalars to ${\mathbb Q} ={\mathbb
Q} \ot_\ZZZ \ZZZ$, the inclusion in the above display becomes an
isomorphism. The claim follows.
\end{proof}

Claim~\ref{jarcA} implies that, in the composition
\[
\underline\Phi : \TTT_{\rada 11}(\Rada \alpha1n)/\Ush_{\rada 11}(\Rada
\alpha1n) \stackrel\pi\epi \TTT_{\rada 11}(\Rada \alpha1n)/K
\stackrel\cong\to \LLL_{\rada 11}(\Rada x1n)',
\]
the kernel of the projection $\pi$ consists of torsion elements. The
second map, induced by $r$, is an isomorphism by Claim~\ref{jarCa}.
Theorem~\ref{Jarca} will thus be established if we prove

\begin{claim}
\label{jaRCA}
The abelian group $T_{\rada 11}(\Rada \alpha1n)/\Ush_{\rada 11}(\Rada
\alpha1n)$  is
torsion-free. 
\end{claim}

To prove the claim we need some properties of configuration spaces.
We do not know a purely algebraic proof.  As in Section~\ref{JArca},
let $\Conf hn$ be the configuration space of $n$ distinct labeled
points in the $h$-dimensional Euclidean space $\RRR^h$, $n,h \geq 1$.
It is known~\cite{voronov:99} that~(\ref{FNe}) induces a cell
decomposition
\[
\CConf hn = 
\textstyle\bigcup_{\bfT \in \Tree^h(n)} \Ivosek, 
\]
of the one-point compactification $\CConf hn$ of $\Conf
hn$, where $\Ivosek$ is the closure of $[\bfT]$ in $\CConf hn$.  
Let $F_p\CConf hn := 
\bigcup\{\Ivosek;\ \dim(\Ivosek) \leq p\}$. We have an
increasing bounded filtration
\[
\emptyset \subset F_{h+n-1}\CConf hn \subset \cdots \subset
F_{hn}\CConf hn = \CConf hn
\]
which induces a spectral sequence converging to the reduced homology
$\overline H_*(\CConf hn)$.  Let us denote by $(G^h_*(n),\pa^h)$ the
$E^1$-term of this spectral sequence desuspended
$(h+1)$-times,\footnote{The number $h+1$ equals the dimension of the
affine group of ${\RRR^h}$.}
\[
G^h_*(n) := \bigoplus_{p+q=*+h+1} 
\overline H_{p+q} (F_p\CConf hn/F_{p-1}\CConf hn),
\]
with the induced differential. The quotient $F_p\CConf
hn/F_{p-1}\CConf hn$ is isomorphic to the cluster of $p$-dimensional
spheres indexed by trees in $\Tree^h(n)$ with $p$ edges, therefore
$G^h_d(n)$ is spanned by labeled pruned $h$-trees with $d + h + 1$
edges, $G^h_d(n) = \Span(\Tree^h_d(n))$. So $G^h_*(n)$ agrees with the
graded abelian group introduced under the same name on page~\pageref{147}.  Our
spectral sequence degenerates at the $E^2$-level, therefore,
\[
H_*(G^h_*(n),\pa^h) \cong \overline H_{* + h +1}(\CConf hn),
\]
while,
by the Poincar\'e-Lefschetz duality~\cite[Section~13.3]{vassiliev:book-intro},
\[
\overline H_*(\CConf hn) \cong H^{hn-*}(\Conf hn).
\]
The cohomology in the right hand side of the above display is
known~\cite[Theorem~1.6]{cohen-taylor:CM93}; for the purpose of this
paper, it is enough to recall that $H^*(\Conf hn)$ is torsion-free and
nontrivial only in degrees $i(h-1)$, $0\leq i \leq n-1$.  The above
results combine into

\begin{claim}
The homology $H_*(G^h_*(n),\pa^h)$ is torsion-free and concentrated
in degrees $(n-2) + i(h-1)$, $0 \leq i \leq n-1$.
\end{claim}

There is a natural degree zero dg-monomorphism $\iota:
(G^h_*(n),\pa^h) \hookrightarrow (G^{h+1}_*(n),\pa^{h+1})$ that sends
the generator $e_T$ indexed by $T \in \Tree^h(n)$ into the generator
$e_{s(T)}$ indexed by the suspension $s(T) \in \Tree^{h+1}(n)$.  By
simple combinatorics, $G^h_d(n) = G^{h+1}_d(n)$ whenever $d \leq
h(n-1)-1$. Let us denote
\[
(G_*(n),\pa) := \dirlim(G^h_*(n),\pa^h).
\]
It is clear that $G_*(n)$ is the span of the graded set
$\Treeinf_*(n)$ so it coincides with the graded abelian group introduced
in Definition~\ref{JArKa}. It is not difficult to see that the
differential $\pa$ is the differential $\pal$ of
Definition~\ref{jarka} but we will not use this fact.
Observe that $G_d(n) = 0$ for $d < n-2$. The above results imply

\begin{claim}
\label{jarKa}
One has $H_d(G_*(n),\pa)= 0$ for $d\not= n-2$ while $H_{n-2}(G_*(n),\pa)$ is
torsion-free.
\end{claim}

\begin{figure}[ht]
\begin{center}
\unitlength=1.000000pt
\begin{picture}(220.00,140.00)(0.00,10.00)
\put(-120,0){
\put(85.00,55.00){\makebox(0.00,0.00)[lb]{\scriptsize $1$}}
\put(75.00,130.00){\makebox(0.00,0.00)[r]{\scriptsize $n+1$}}
\put(75.00,80.00){\makebox(0.00,0.00)[r]{\scriptsize $3$}}
\put(75.00,30.00){\makebox(0.00,0.00)[r]{\scriptsize $2$}}
\put(94.00,100.00){\makebox(0.00,0.00){\scriptsize $\vdots$}}
\put(94.00,120.00){\makebox(0.00,0.00){\scriptsize $\sigma(n)$}}
\put(94.00,70.00){\makebox(0.00,0.00){\scriptsize $\sigma(2)$}}
\put(94.00,20.00){\makebox(0.00,0.00){\scriptsize $\sigma(1)$}}
\put(80.00,20.00){\makebox(0.00,0.00){\scriptsize $\bullet$}}
\put(80.00,120.00){\makebox(0.00,0.00){\scriptsize $\bullet$}}
\put(80.00,70.00){\makebox(0.00,0.00){\scriptsize $\bullet$}}
\thinlines
\put(80.00,10.00){\line(0,1){135.00}}
\put(20.00,50.00){\line(1,0){180.00}}
\thicklines
\put(80.00,50.00){\vector(1,0){20.00}}
\put(80.00,120.00){\vector(0,1){20.00}}
\put(80.00,70.00){\vector(0,1){20.00}}
\put(80.00,20.00){\vector(0,1){20.00}}
}
\put(100,0){
\put(-20,0){
\put(85.00,55.00){\makebox(0.00,0.00)[lb]{\scriptsize $1$}}
\put(75.00,130.00){\makebox(0.00,0.00)[r]{\scriptsize $s+2$}}
\put(75.00,80.00){\makebox(0.00,0.00)[r]{\scriptsize $4$}}
\put(75.00,30.00){\makebox(0.00,0.00)[r]{\scriptsize $3$}}
\put(94.00,100.00){\makebox(0.00,0.00){\scriptsize $\vdots$}}
\put(94.00,120.00){\makebox(0.00,0.00){\scriptsize $\rho(s)$}}
\put(94.00,70.00){\makebox(0.00,0.00){\scriptsize $\rho(2)$}}
\put(94.00,20.00){\makebox(0.00,0.00){\scriptsize $\rho(1)$}}
\put(80.00,20.00){\makebox(0.00,0.00){\scriptsize $\bullet$}}
\put(80.00,120.00){\makebox(0.00,0.00){\scriptsize $\bullet$}}
\put(80.00,70.00){\makebox(0.00,0.00){\scriptsize $\bullet$}}
\put(80.00,10.00){\line(0,1){135.00}}
\thicklines
\put(80.00,50.00){\vector(1,0){20.00}}
\put(80.00,120.00){\vector(0,1){20.00}}
\put(80.00,70.00){\vector(0,1){20.00}}
\put(80.00,20.00){\vector(0,1){20.00}}
}
\put(50,0){
\put(85.00,55.00){\makebox(0.00,0.00)[lb]{\scriptsize $2$}}
\put(65.00,130.00){\makebox(0.00,0.00){\scriptsize $n+2$}}
\put(65.00,80.00){\makebox(0.00,0.00){\scriptsize $s+4$}}
\put(65.00,30.00){\makebox(0.00,0.00){\scriptsize $s+3$}}
\put(94.00,100.00){\makebox(0.00,0.00){\scriptsize $\vdots$}}
\put(94.00,120.00){\makebox(0.00,0.00){\scriptsize $\rho(n)$}}
\put(102.00,70.00){\makebox(0.00,0.00){\scriptsize $\rho(s+2)$}}
\put(102.00,20.00){\makebox(0.00,0.00){\scriptsize $\rho(s+1)$}}
\put(80.00,20.00){\makebox(0.00,0.00){\scriptsize $\bullet$}}
\put(80.00,120.00){\makebox(0.00,0.00){\scriptsize $\bullet$}}
\put(80.00,70.00){\makebox(0.00,0.00){\scriptsize $\bullet$}}
\put(80.00,10.00){\line(0,1){135.00}}
\thicklines
\put(80.00,50.00){\vector(1,0){20.00}}
\put(80.00,120.00){\vector(0,1){20.00}}
\put(80.00,70.00){\vector(0,1){20.00}}
\put(80.00,20.00){\vector(0,1){20.00}}
}
\thinlines
\put(20.00,50.00){\line(1,0){190.00}}
}
\end{picture}
\end{center}
\caption{\label{JARKa}Generic points of the Fox-Neuwirth cells.}
\end{figure}
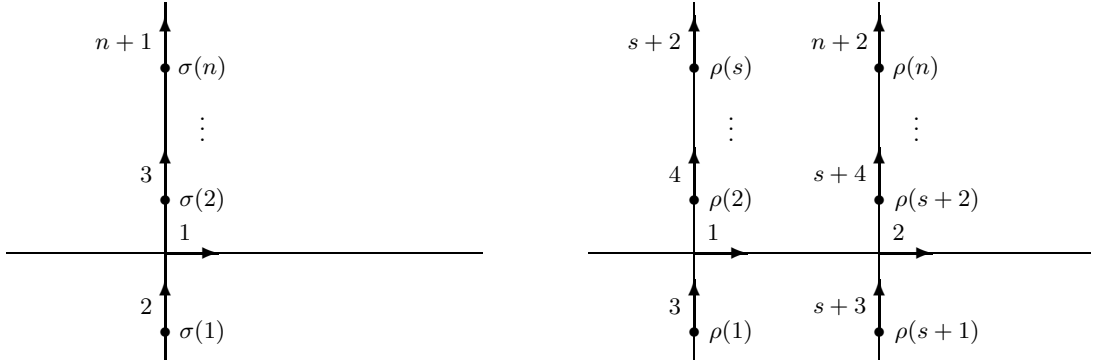

Let us calculate $H_{n-2}(G_*(n),\pa)$. Since $G_d^h(n) = G_d^2(n)$
for $d \leq n-1$, clearly $H_{n-2}(G_*(n),\pa) \cong
H_{n-2}(G_*^2(n),\pa^2)$.
The space $G^2_{n-2}(n)$ is spanned by labeled corollas of
height two
\begin{equation}
\label{JarKa}
\raisebox{-3em}{\rule{0pt}{0pt}}
\unitlength1cm
\begin{picture}(2.00,1.5)(1,0)
\thicklines
\multiput(0,-1)(0,1){3}{\multiput(0.2,0)(0.1,0){27}{\makebox(0,0){$\cdot$}}}
\put(.3,1.3){\makebox(0,0){\scriptsize $\sigma(1)$}}
\put(.95,1.3){\makebox(0,0){\scriptsize $\sigma(2)$}}
\put(2.5,1.3){\makebox(0,0){\scriptsize $\sigma(n)$}}
\put(4.5,.1){\makebox(0,0){,\ $\sigma \in \Sigma_n$,}}
\put(1.6,.8){\makebox(0,0){$\cdots$}}
\put(1.5,0){\line(1,1){1}}
\put(1.5,0){\line(0,-1){1}}
\put(1.5,0){\line(-1,2){.5}}
\put(1.5,0){\line(-1,1){1}}
\end{picture}
\end{equation}
therefore $G_{n-2}(n) \cong \TTT_{\rada 11}(\Rada \alpha 1n)$, with
the $\Sigma_n$-equivariant isomorphism sending the above corolla into
the generator $\alpha_{\sigma(1)} \ot \cdots \ot \alpha_{\sigma(n)}$.
The corolla~(\ref{JarKa}) corresponds to the \label{ukl}cell of
$\CConf 2n$ whose generic point is shown in Figure~\ref{JARKa} (left).

The little arrows numbered $\rada 1{n+1}$ indicate a frame in the
tangent bundle determining the orientation.
Likewise, $G^2_{n-1}(n)$ is spanned by trees of height two
\begin{equation}
\label{Jaruska-milacek}
\raisebox{-1.1cm}{\rule{0em}{0em}}
\unitlength1cm
\begin{picture}(2.00,1.4)(1,0)
\thicklines
\put(-1.5,0){
\multiput(0,-1)(0,1){3}{\multiput(0.2,0)(0.1,0){57}{\makebox(0,0){$\cdot$}}}
\put(.3,1.3){\makebox(0,0){\scriptsize $\rho(1)$}}
\put(.95,1.3){\makebox(0,0){\scriptsize $\rho(2)$}}
\put(2.5,1.3){\makebox(0,0){\scriptsize $\rho(s)$}}
\put(1.6,.8){\makebox(0,0){$\cdots$}}
\put(1.5,0){\line(1,1){1}}
\put(1.5,0){\line(-1,2){.5}}
\put(1.5,0){\line(-1,1){1}}
}
\put(1.5,0){
\put(.35,1.3){\makebox(0,0){\scriptsize $\rho(s\hskip-.3em+\hskip-.3em1)$}}
\put(2.5,1.3){\makebox(0,0){\scriptsize $\rho(n)$}}
\put(1.6,.8){\makebox(0,0){$\cdots$}}
\put(3,0.1){\makebox(0,0)[l]{,\ $\rho \in \Sigma_n,\ 1 \leq s < n$,}}
\put(1.5,0){\line(1,1){1}}
\put(1.5,0){\line(-1,2){.5}}
\put(1.5,0){\line(-1,1){1}}
}
\put(1.5,-1){\line(3,2){1.5}}\put(1.5,-1){\line(-3,2){1.5}}
\end{picture}
\end{equation}
representing the cell whose generic point is shown in
Figure~\ref{JARKa} (right).

By imagining how a generic point of the cell corresponding to the
tree~(\ref{Jaruska-milacek}) moves to the boundary, one sees that
the differential $\pa$ sends this tree into the element that, under
the isomorphism $G_{n-2}(n) \cong \TTT_{\rada 11}(\Rada \alpha 1n)$,
equals 
\begin{equation}
\label{jArKa}
\sum_{\tau \in \Sigma_{s,t}} \sgn(\tau) \cdot
\alpha_{\rho\tau(1)}\ot \cdots\ot \alpha_{\rho\tau(n)}
\end{equation} 
(observe the
$\sgn(\tau)$-factor), therefore
\begin{equation}
\label{JARkA}
H_{n-1}(G_*,\pa) \cong \TTT_{\rada 11}(\Rada \alpha 1n)/
\widetilde{\Ush}_{\rada 11}(\Rada \alpha 1n),
\end{equation}
where $\widetilde{\Ush}_{\rada 11}(\Rada \alpha 1n)$ denotes the span
of elements in~(\ref{jArKa}). We, however, have

\begin{claim}
\label{jarKA}
Let $\sgn$ denote the signum representation.
There is a $\Sigma_n$-equivariant isomorphism
\[
\TTT_{\rada 11}(\Rada \alpha 1n)/
\widetilde{\Ush}_{\rada 11}(\Rada \alpha 1n) \cong \sgn \ot
\TTT_{\rada 11}(\Rada \alpha 1n)/\Ush_{\rada 11}(\Rada \alpha 1n).
\]
\end{claim}

\begin{proof}
The isomorphism $\Psi : \TTT_{\rada 11}(\Rada \alpha 1n) \stackrel\cong\to 
\sgn \ot \TTT_{\rada 11}(\Rada \alpha 1n)$ given by
\[
\Psi(
\alpha_{\sigma(1)} \ot \cdots \ot
\alpha_{\sigma(n)}) := \sgn(\sigma) \ot \alpha_{\sigma(1)} \ot \cdots \ot
\alpha_{\sigma(n)}, \ \sigma \in \Sigma_n,
\]
clearly restricts to an isomorphism 
$\widetilde{\Ush}_{\rada 11}(\Rada \alpha 1n) \cong {\Ush}_{\rada
11}(\Rada \alpha 1n)$ and induces the isomorphism of the claim. 
\end{proof}

By Claim~\ref{jarKa} and isomorphism~(\ref{JARkA}), 
$\TTT_{\rada 11}(\Rada \alpha 1n)/
\widetilde{\Ush}_{\rada 11}(\Rada \alpha 1n)$ is torsion-free,
which, by Claim~\ref{jarKA}, proves Claim~\ref{jaRCA} and therefore
also Theorem~\ref{Jarca}.
Recall that the space $G_*(n)$ inherits a natural free
$\Sigma_n$-action given by relabeling the spanning trees. As a
combination of the above results we get

\begin{theorem}
\label{Stale_chodim_s_Jarkou}
The tree complex $(G_*(n),\pa)$ is a $\Sigma_n$-free resolution of the
$\Sigma_n$-module $\sgn \otimes \Lie(n)'$. 
\end{theorem}

Tensoring $(G_*(n),\pa)$ with the signum representation therefore
leads to a free resolution of $\Lie(n)'$.  The following example shows
that $(G_*(n),\pa)$ is not the smallest possible.

\begin{example}
{ \rm Inspecting Figure~\ref{fig:3}, one sees that $G_1(3)$ has one
$\Sigma_3$-generator $[1|2|3]$, $G_2(3)$ two $\Sigma_3$-generators
$[1|2||3]$ and $[1||2|3]$, and $G_3(3)$ three $\Sigma_3$-generators
$[1||2||3]$,$[1|2|||3]$ and $[1|||2|3]$. In general, the number of
generators of the free $\ZZZ[\Sigma_3]$-module $G_d(3)$ equals $d$.  One
has
\[
\pa ([1|2||3]) = [1|2|3] - [1|3|2] + [3|1|2] = \pa ([3||1|2]) 
\]
and
\[
\pa ([1|2|||3]) = [1|2||3] - [3||1|2].
\]
This shows that one of the generators of $G_2(3)$ is superfluous, so
there is a resolution with only one $\Sigma_3$-generator in degree
$2$. We believe in the existence of a free $\Sigma_3$-resolution
$(\widetilde G_*(3),\pa)$ of
$\Lie(3)'$  in which $\widetilde G_d(3)$ has
$[\frac{d-1}2]+1$ $\Sigma_3$-generators, where $[-]$ denotes the 
integral part, $d \geq 1$.

Notice that the group ring $R[\Sigma_n]$ does not have good properties
even for $R$ a characteristic zero field. For the augmentation
ideal ${\mathcal I}$ of $R[\Sigma_n]$ one has ${\mathcal I}/ {\mathcal
I}^2 = 0$, so there is no good notion of minimality of
$\Sigma_n$-resolutions of $\Sigma_n$-modules.  }
\end{example}

\def\cprime{$'$}

\section*{Glossary}
\label{gloss}
\allowdisplaybreaks
\[
\def\arraystretch{1.2}
\allowdisplaybreaks
\begin{array}{lll}
\allowdisplaybreaks
\gls{\Free(-)}{free operad functor}{jarca!}
\\
\gls{\uFree(-)}{free non-$\Sigma$ operad functor}{jarca!!}
\\
\uTree^h(n),&\mbox {set of pruned $h$-trees with $n$ tips,}
&\mbox{ page~\pageref{prun}}
\\
\Tree^h(n),&\mbox {set of labeled pruned $h$-trees with $n$ tips,}
&\mbox{ page~\pageref{445}}
\\
\uTreeinf(n),&\mbox {set of reduced trees of arbitrary height, 
equals $\dirlim \uTree^h(n)$,} &\mbox{ page~\pageref{098}}
\\
\Treeinf(n),&\mbox {set of reduced labeled trees of arbitrary height, 
equals $\dirlim \Tree^h(n)$,} &\mbox{ page~\pageref{087}}
\\
G^h(n),&\mbox {right $\Sigma_n$-module spanned by $\Tree^h(n)$,} 
&\mbox{ page~\pageref{147}}
\\
G(n),&\mbox {right $\Sigma_n$-module spanned by $\Treeinf(n)$, equals
  $\dirlim G^h$,} 
&\mbox{ page~\pageref{JArKa}}
\\
\sfB(A),&\mbox {bar construction of an associative algebra $A$,}
&\mbox { page~\pageref{jhg}}
\\
\gls{\sfB^h(A)}{$h$th iterate of the bar construction of an
  ass.\ comm.\ algebra}{qwe}
\\
\gls{\hatB_*^h(A)}{desuspension $\downarrow^{h+1}\sfB^h_*(A)$}{ty}
\\
\gls{\hatB^\infty(A)}{direct limit $\dirlim \hatB_*^h(A)$}{3r}
\\
\gls{\Conf hn}{configuration space of distinct labeled 
points in $\RRR^n$}{qaz}
\\
\gls{\osfF_h(n)}{moduli space $\Conf hn/\Aff h$}{wsx}
\\
\gls{\sfF_h(n)}{Fulton-MacPherson compactification of $\osfF_h(n)$}{wsx}
\\
\gls{\sfF(n)}{direct limit $\dirlim \sfF_h(n)$}{tgb}
\\
\gls{[\bfT]}{cell of $\Conf hn$ indexed by a labeled tree
  $\bfT = (T,\ell)$}{yhn} 
\\
\gls{\mu[\bfT]}{quotient $[\bfT]/\Aff h$}{plm}
\\
\gls{\Treeregh(n)}{subset of $\Tree^h(n)$ of trees of dimension
$< \dcrit^h(n)$}{qw}
\\
\gls{\Treereg(n)}{subset of $\Tree(n)$ of trees of dimension
$< \dcrit(n)$,  equals $\dirlim \Treeregh_*$}{bv}
\\
\gls{\sfFreg_h ,\ \sfFreg}{regular skeleton of the configuration
    operad $\sfF_h$ resp.~$\sfF$}{rt}
\\
\gls{\mu[\sigma]}{cell of $\sfF$ corresponding to a face
$\sigma: T\to S$ of $T$}{dnes_oslici}
\\
\gls{(\Rada \sigma1n)}{notation for a 
permutation $\sigma \in \Sigma_n$, $\sigma_i := \sigma^{-1}(i)$, $1
\leq i \leq n$}{Dnes_mi_trhali_zub!}
\\ 
\gls{\Tam}{Tamarkin cell, equals $\mu[1|2||3|4||5|6]$}{tam}
\\
\glsfin{\Mar}{a $4$-dimensional bad cell, equals $\mu[1|2|||3|4]$}{mar}
\end{array}
\]

\end{document}